\documentclass[11pt,a4paper]{scrartcl}
\usepackage[T1]{fontenc}
\usepackage{amsmath,amssymb,amsthm}
\usepackage{authblk}

\usepackage{enumitem}
\usepackage[colorlinks=true,linkcolor=blue,urlcolor=blue,citecolor=blue]{hyperref}
\hypersetup{bookmarksdepth=2} \allowdisplaybreaks
\KOMAoptions{BCOR=0mm}
\areaset{159.2mm}{246.2mm}

\newtheorem{theorem}{Theorem}[section]
\newtheorem{lemma}[theorem]{Lemma}
\newtheorem{proposition}[theorem]{Proposition}
\newtheorem{corollary}[theorem]{Corollary}
\theoremstyle{remark}
\newtheorem{remark}[theorem]{Remark}
\theoremstyle{definition}
\newtheorem{definition}[theorem]{Definition}
\newtheorem{claim}[theorem]{Claim}
\theoremstyle{example}
\newtheorem{example}{Example}

\numberwithin{equation}{section}

\newcommand{\R}{\mathbb{R}}
\newcommand{\dd}{\,\mathrm{d}}

\newcommand{\rs}{\rho_{\sharp}}
\newcommand{\cs}{c_{\sharp}}

\title{\fontfamily{phv}\selectfont Perturbation results for the Choquard equation in the sublinear case}
\author[1]{Simone Secchi\thanks{Supported by the GNAMPA-INdAM 2026 project  ‘‘Nonlinearità e mancanza di compattezza fra ellitticità e dispersione’’. }}
\affil[1]{Università degli Studi di Milano - Bicocca}
\date{\today}

\begin{document}
\maketitle

\begin{flushright}
\textit{To the memory of Antonio Ambrosetti}
\end{flushright}

\begin{abstract}
We study perturbations of the positive ground state of the Choquard equation
\[
-\Delta u+u=(I_\alpha*|u|^p)|u|^{p-2}u\qquad\text{in }\mathbb R^N
\]
by means of a Lyapunov--Schmidt reduction. For $p\ge2$ the energy is
$C^2$ on $H^1$ and a standard machinery applies. The case $\alpha=2$,
$1+2/N<p<2$ is different: the ground state has an algebraic decay and
the energy is not twice differentiable on $H^1$. We introduce a
weighted space adapted to the tail and we perform a finite-dimensional
reduction on it. The outcome is a positive solution of a localized
perturbation problem near every critical point of the first-order
interaction functional with nonzero local degree. For the singularly
perturbed problem with potential $V$ we obtain, under a small relative
oscillation assumption, solutions concentrating near every critical
point of $V$ with nonzero local degree, and a Lusternik--Schnirelman
lower bound near compact sets of local minima. We also consider the
variant in which $V$ multiplies both the linear and nonlocal terms,
and we construct solutions that concentrate near isolated critical
points of nonzero degree.
\end{abstract}

\tableofcontents

\section{Introduction}\label{sec:introduction}

The Choquard equation
\begin{equation}\label{eq:intro-model}
 -\Delta u+u=(I_\alpha*|u|^p)|u|^{p-2}u \qquad\text{in }\mathbb R^N
\end{equation}
is a nonlocal equation that arises afrom the coupling of a (linear) Schr\"odinger equation and a Poisson equation. The 
nonlinearity couples the value of $u$ at a point to the Riesz potential generated by $|u|^p$. 

The equation has a classical physical origin in the Hartree-type description of a quantum particle interacting with its own induced field. In the case $\alpha=2$, it is closely related to the pioneering work of Pekar on the polaron and to the Choquard model for an electron trapped in a polar crystal; in the same spirit, the equation appears in the description of one-component plasmas and other mean-field models with long-range interactions. Mathematically, this is the Euler--Lagrange equation of a nonlocal energy where the potential generated by the density $|u|^p$ is represented by a Riesz convolution kernel. For the physical background we refer to \cite{Pekar,Choquard,Lieb,Liebb} and, for the modern mathematical development, to \cite{Lions,MVS1,MVS2}. 

The existence theory for Choquard-type equations has developed through a series of fundamental contributions. In the classical Hartree setting, Lieb proved the existence of a minimizer for the corresponding variational problem, and Lions developed the concentration-compactness framework that underlies the existence of ground states in the nonlocal setting; see \cite{Lieb,Lions}. For the Choquard equation with general exponents, existence and qualitative properties of positive solutions were later established in the works of Moroz and Van Schaftingen, together with the systematic study of regularity, decay, and radial symmetry; see \cite{MVS1,MVS2}. In the sublinear regime considered here, the recent analysis of decay and nondegeneracy in \cite{XZZ1,XZZ2} is essential, and it is the starting point for the perturbation theory developed in the present paper.

\medskip

We will consider the case $N\ge3$, $I_\alpha(x)=c_{N,\alpha}|x|^{\alpha-N}$ and
\[
 \frac{N+\alpha}{N}<p<\frac{N+\alpha}{N-2},
\]
which is exactly the range in which \eqref{eq:intro-model} has nontrivial finite energy solutions. In this paper we investigate some classes of perturbations of \eqref{eq:intro-model}, with the emphasis on the sublinear regime $p<2$.

Our method is based on a Lyapunov--Schmidt reduction around the translation orbit of a ground state, in the form introduced by Ambrosetti and Badiale \cite{AB1,AB2} (see also \cite{AM1}). Equation \eqref{eq:intro-model} is translation invariant, so if the positive ground state $U$ is unique up to translations, the family $U_\xi=U(\cdot-\xi)$ is an $N$-dimensional manifold of critical points of the energy. The reduction works near this manifold under two conditions. The kernel of the linearized operator at $U$ must be exactly $\operatorname{span}\{\partial_1U,\dots,\partial_NU\}$ (which is usually described in terms of \textit{nondegeneracy} of $U$), and the energy must be regular enough, in a space where the linearization is Fredholm of index zero, for the implicit function theorem to apply.

Nondegeneracy in the sublinear range has just been proved by Xia, Zhang and Zhang \cite{XZZ2}, for $\alpha=2$ and $\frac{N+2}{N}<p<\frac{N+2}{N-2}$ (see Section~\ref{sec:setting}). Regularity is the real issue. For $p\ge2$ the energy is of class $C^2$ on $H^1(\mathbb R^N)$ and the reduction is standard: nevertheless, in Section~\ref{sec:pge2} we describe it as a toy model to which more difficult problems will be compared. When $\alpha=2$ and $p<2$, the ground state decays only algebraically, $U(x)\asymp(1+|x|)^{-\beta}$ with $\beta=(N-2)/(2-p)$ \cite{XZZ1}, and the energy fails to be twice Fr\'echet differentiable at $U$ on $H^1$. The reason is simple: a bump comparable to $U$, placed far from the origin, is small in $H^1$ but it changes the \textit{nonlocal} nonlinearity by a fixed proportion of its size (Section~\ref{sec:psmall}). The $H^1$-based argument therefore breaks down, and we replace $H^1$ by a weighted space that sees the tail of $U$.

\medskip

From now on in this introduction $\alpha=2$ and $1+\frac2N<p<\min \left\{2,\frac{N+2}{N-2} \right\}$. We set
\begin{align*}
 X&=\left\{w\ \text{measurable}\mid \|w\|_X:=\operatorname*{ess\,sup}_{x\in\mathbb R^N}(1+|x|)^\beta|w(x)|<\infty\right\},\\
 Y&=\{w\in X\mid  \nabla w\in L^2(\mathbb R^N)\},
\end{align*}
with norm $\|w\|_Y=\|w\|_X+\|\nabla w\|_{L^2}$. The weighted supremum norm rules out far-away bumps of the kind just described. The gradient term gives the continuous embedding $Y\hookrightarrow H^1$, so we can still use $H^1$-orthogonality to the translation modes.

Let $T_0(u)=(-\Delta+1)^{-1}\left[(I_2*|u|^p)|u|^{p-2}u\right]$, so that the solutions of \eqref{eq:intro-model} are the fixed points of $T_0$. Its derivative at $U$ is a bounded operator on $Y$. Our first result is that $I-T_0'(U)$ is Fredholm of index zero on $Y$. Together with nondegeneracy, this identifies its kernel and range, and makes it invertible on the orthogonal complement of the translation modes.

\begin{theorem}\label{thm:intro-fredholm}
Assume $N\geq3$, $\alpha=2$, and
$1+2/N<p<\min\{2,(N+2)/(N-2)\}$. Let $U$ be the positive radial
ground state of \eqref{eq:intro-model}, and let $Y$ be the weighted space defined
above. Then $T_0'(U)\colon Y\to Y$ is bounded and $I-T_0'(U)$ is Fredholm of index
zero. Moreover,
\[
 \ker(I-T_0'(U))
 =\operatorname{span}\{\partial_{x_1}U,\ldots,\partial_{x_N}U\},
\]
and its range is the $H^1$-orthogonal complement in $Y$ of this
translation space. In particular, the restriction of $I-T_0'(U)$ to
that complement is an isomorphism.
\end{theorem}

As a first application we consider the localized perturbation
\begin{equation}\label{eq:intro-pert}
 -\Delta u+(1+\varepsilon a(x))u
   =(I_2*u^p)u^{p-1}\qquad\text{in }\mathbb R^N,
\end{equation}
where $a$ is bounded and of class $C^1$, has bounded and uniformly continuous gradient, and tends to zero at infinity. For small $\varepsilon$ the reduction produces, for every $\xi\in\mathbb R^N$, a correction of $U_\xi$ and a reduced functional $\Phi_\varepsilon(\xi)$ whose critical points give solutions. The first-order term of $\Phi_\varepsilon$ is the interaction energy
\[
 \Gamma(\xi)=\frac12\int_{\mathbb R^N}
 a(x)U(x-\xi)^2\dd x ,
\]
also called the (first-order expansion of) \textit{Poincaré-Melnikov function} associated to the finite-dimensional reduction.
If $a$ has a sign, $\Gamma$ has a global extremum, since it does not vanish identically and tends to zero at infinity; this gives a solution for every small $\varepsilon$. In general $\varepsilon^{-1}\nabla\Phi_\varepsilon=\nabla\Gamma+O(\varepsilon)$ uniformly, so there are solutions near every critical point of $\Gamma$ with nonzero local degree (Theorem~\ref{thm:general}). The sign-definite case reads as follows.

\begin{theorem}\label{thm:intro-perturbation}
Assume $N\geq3$ and
$1+2/N<p<\min\{2,(N+2)/(N-2)\}$. Let $a$ satisfy assumption
\textup{(A)} of Section~\ref{sec:fixedpoint}, and suppose either
$a\geq0$, $a\not\equiv0$, or $a\leq0$, $a\not\equiv0$. Then there
exists $\varepsilon_0>0$ such that, for every
$0<|\varepsilon|<\varepsilon_0$, equation
\eqref{eq:intro-pert} has a positive solution $u_\varepsilon\in H^1(\mathbb R^N)$.
More precisely, for suitable points $\xi_\varepsilon\in\mathbb R^N$,
\[
 u_\varepsilon=U(\cdot-\xi_\varepsilon)
       +w_\varepsilon(\cdot-\xi_\varepsilon),
 \qquad \|w_\varepsilon\|_Y\leq C|\varepsilon|,
 \qquad u_\varepsilon\geq\frac12U(\cdot-\xi_\varepsilon).
\]
Let $\Phi_\varepsilon$ be the reduced functional, $s=\operatorname{sgn}\varepsilon$, and
$\sigma=+1$ if $a\ge0$, $\sigma=-1$ if $a\le0$. Then $\xi_\varepsilon$ can be chosen as a
maximum point of $\sigma s\,\Phi_\varepsilon$ on a fixed ball; that is, for $a\ge0$ a maximum
point of $\Phi_\varepsilon$ if $\varepsilon>0$ and a minimum point if $\varepsilon<0$, and
the reverse for $a\le0$. Every accumulation point of
$\xi_\varepsilon$ as $\varepsilon\to0$ is a global maximum point of
$\Gamma$ if $a\ge0$ and a global minimum point of $\Gamma$ if $a\le0$.
\end{theorem}

Our second application is the semiclassical problem
\begin{equation}\label{eq:intro-semi}
\begin{cases}
 -\varepsilon^2\Delta u+V(x)u
   =\varepsilon^{-2}(I_2*u^p)u^{p-1}, \\
 u>0,\quad u\in H^1(\mathbb R^N),
\end{cases}
\end{equation}
with a positive potential $V \colon \mathbb{R}^N \to \mathbb{R}$. The algebraic tail of $U$ makes this problem more delicate than its analogue with exponentially decaying ground states, which are typical of the NLS corresponding formally to the choice $I_2=\delta$. After rescaling around a point $\xi$, the potential becomes $V(\xi+\varepsilon\,\cdot)$, which is close to $V(\xi)$ only on bounded sets, whereas the weighted norm measures the perturbation on the whole space, tail included. Local information on $V$ near a prospective concentration point is therefore not enough, and we assume that $V$ has small relative oscillation: for some $\lambda_0>0$,
\[
 V\in C^1(\mathbb R^N),\qquad
 \nabla V\ \text{bounded and uniformly continuous},\qquad
 \sup_{x\in\mathbb R^N}
 \left|\frac{V(x)}{\lambda_0}-1\right|<\eta_*(N,p),
\]
where $\eta_*(N,p)>0$ is the threshold specified in Section~\ref{sec:semiclassical}. Under this assumption the weighted reduction applies uniformly in the center parameter. The reduced gradient converges locally uniformly to a positive multiple of $\nabla V$, which allows the use of the Brouwer degree and gives solutions concentrating near topologically stable critical points of $V$.

\begin{theorem}\label{thm:intro-semi}
Assume $N\geq3$ and
$1+2/N<p<\min\{2,(N+2)/(N-2)\}$. Let $V$ satisfy \textup{(V1)} and
\textup{(V2)} of Section~\ref{sec:semiclassical}. If $\xi_0$ is an
isolated critical point of $V$ with
\[
 \deg(\nabla V,B_r(\xi_0),0)\ne0
\]
for all sufficiently small $r>0$, then, for every sufficiently small
such $r$, equation \eqref{eq:intro-semi} has a positive solution for
all sufficiently small $\varepsilon>0$ whose concentration center
$\xi_\varepsilon$ lies in $B_r(\xi_0)$. Moreover
$u_\varepsilon(\xi_\varepsilon+\varepsilon\,\cdot)-U^{[V(\xi_\varepsilon)]}\to0$
in $L^2(\mathbb R^N)$, where $U^{[\lambda]}$ denotes the ground state of the
frozen equation $-\Delta v+\lambda v=(I_2*v^p)v^{p-1}$.
\end{theorem}

The same convergence of the reduced gradient gives a Lusternik--Schnirelman lower bound for the number of critical points of the reduced functional near a compact set of local minima of $V$. An orbit-separation lemma shows that distinct parameters of the reduction give distinct solutions of the original equation, so the count is also a lower bound for the number of distinct positive solutions (Proposition~\ref{prop:SCcat}).

\bigskip

For $p<2$ the existence theory for the Choquard equation has been proved by  variational methods. Ground states of \eqref{eq:intro-model} exist exactly in the range above \cite{MVS0,MVS3} (see also the survey \cite{MVS2}); in the sublinear range they decay algebraically rather than exponentially, which is the source of the difficulties described above. For the semiclassical problem, Moroz and Van Schaftingen \cite{MVS1} construct by penalization solutions concentrating at local minima of $V$, but their main theorems require $p\ge2$, and the case $p<2$ is only discussed in a closing remark. This gap was filled by Cingolani and Tanaka \cite{CT}. For nonlinearities satisfying their conditions (f1)--(f5), which for the pure power cover the whole interval $\frac{N+\alpha}{N}<p<\frac{N+\alpha}{N-2}$ including $p<2$, and for a continuous potential with $0<\inf V\le\sup V<\infty$ and a bounded open set $\Omega$ with $\inf_\Omega V<\inf_{\partial\Omega}V$, they prove that for small $\varepsilon$ there are at least $\operatorname{cuplength}(K)+1$ positive solutions concentrating in the set $K$ of minima of $V$ in $\Omega$. No smallness of the oscillation of $V$ is needed.

The reduction method, on the other hand, needs nondegeneracy of the ground state, which is delicate for nonlocal equations. For the Hartree case $\alpha=2$, $p=2$ it was proved by Lenzmann \cite{Len} for $N=3$ and by Chen \cite{Chen} for $N=3,4,5$; Chen also carries out the Lyapunov--Schmidt reduction and obtains $\ell(M)$ semiclassical solutions, where $\ell$ is the cup-length, near a nondegenerate critical manifold $M$ of the potential. Xiang \cite{Xia} treated the Coulomb kernel in $\mathbb R^3$ with $p>2$ close to $2$. All these results concern $p\ge2$, where the energy is of class $C^2$ on $H^1$ and the standard reduction applies verbatim. Uniqueness and nondegeneracy for $\alpha=2$ and $\frac{N+2}{N}<p<\frac{N+2}{N-2}$, hence also for $p<2$, were obtained only recently by Xia, Zhang and Zhang \cite{XZZ2} (see \cite{XZZ1} for the decay estimates). As far as we know, no reduction argument had been carried out for $p<2$.

\medskip

In this paper we isolate the obstruction to the standard method when $\alpha=2$ and $p<2$: it is a lack of regularity, not a spectral one, because the effective mass at infinity $2-p$ is positive and nondegeneracy holds. The weighted space $Y$ removes it. We prove that the linearization is Fredholm of index zero on $Y$, that the energy is $C^{2,1}_{\mathrm loc}$ on $Y$, and that the fixed-point map is $C^{1,1}$ from a ball of $X$ into $Y$ (Theorem~\ref{thm:fred}, Propositions~\ref{prop:f0} and \ref{prop:T0}). These are the analytic ingredients that were missing for $p<2$.

The variational theories \cite{MVS1,CT} produce solutions concentrating at minima of $V$. Under the oscillation assumption, the reduction gives a solution near every isolated critical point of $V$ with nonzero local degree, hence also near saddle points (Theorem~\ref{thm:intro-semi}), and near every critical point of $\Gamma$ with nonzero degree in the localized problem (Theorem~\ref{thm:general}). Such solutions cannot be obtained by minimization. The solutions also come with an explicit description. In the localized problem they are $O(\varepsilon)$-close in $Y$ to a translate of $U$, with the center selected by the first-order functional $\Gamma$. In the semiclassical problem they converge to the frozen ground state $U^{[V(\xi_\varepsilon)]}$. In both cases the correction depends $C^1$ on the parameters. Finally, the orbit-separation lemma (Lemma~\ref{lem:orbitsep}) turns a Lusternik--Schnirelman count on the reduced function into a count of distinct positive solutions, in the spirit of \cite{CL,AMS} for local nonlinearities.

The price to pay is that our semiclassical results require a \emph{global} small relative oscillation of $V$ (assumption (V2)), forced by the algebraic tail, whereas \cite{CT} only needs a local well condition. For minima of $V$ the theorem of \cite{CT} is therefore far more general; our category count $\operatorname{cat}(K,K_\delta)$ is of the same type as the classical ones \cite{CL,AMS} but holds under a much more restrictive hypothesis. The reduction is also available only for $\alpha=2$, for the pure power, and for $\frac{N+2}{N}<p<2$ in the sublinear regime, because nondegeneracy is known only there, while \cite{MVS1,CT} cover all $\alpha\in(0,N)$ and, in \cite{CT}, general nonlinearities. Finally, we need $V\in C^1$ with uniformly continuous gradient, whereas the variational approach needs only continuity. The present results are thus complementary to \cite{MVS1,CT}. They do not enlarge the class of potentials admitting a solution concentrating at a minimum, but they give a reduction theory in the sublinear range, together with the additional information described above. 

\medskip

The last section studies the problem with the potential appearing on both sides of the equation,
\begin{equation*}
-\varepsilon^2\Delta u+V(x)u=\varepsilon^{-2}\left(I_2*(Vu^p)\right)Vu^{p-1},	
\end{equation*}
which is the case $K=V$ in the notation of the paper. The key point is that this structure is substantially better behaved than the semiclassical problem considered earlier: the large oscillation of $V$ away from the concentration region is harmless, because it only appears in terms carrying an extra decay. Under the mild hypothesis $V\in C^1$, $0<m\le V\le M<\infty$, and bounded uniformly continuous $\nabla V$, one can rescale around a point $\xi$ and solve a Lyapunov--Schmidt reduction for the profile $U$ in a uniform class of potentials. The reduced functional is then defined by evaluating the natural energy on the corrected profile, and its critical points give positive solutions of the original equation. In the constant-potential case one recovers the dilated ground state $U_\lambda=U(\sqrt\lambda\,\cdot)$, while for general $V$ the reduced equation is driven by the limit profile $\Upsilon\circ V$ and the gradient satisfies
\[
\nabla\Phi_\varepsilon(\xi)\to -c_*\,V(\xi)^{-N/2}\nabla V(\xi),\qquad c_*>0.
\]
Hence the concentration centers are asymptotically located at critical points of $V$, and the reduced problem yields a $C^1$ family of critical points with the natural constraint $\tilde c_\xi=0\iff\nabla\Phi_\varepsilon(\xi)=0$. The main novelty of this section is that, unlike in the semiclassical problem with a fixed coefficient on the nonlocal term, no global small-oscillation assumption on $V$ is needed; the reduction works in the full class $0<m\le V\le M$ as soon as the local potential is controlled near the concentration point.

\medskip

As a concluding remark, in the recent paper \cite{APS}, the authors consider the local NLS $-\Delta u+Vu=V|u|^{p-1}u$, and show that a direct variational approach to such a problem can be troublesome. In particular, this approach depends strongly on the dimension and on the oscillation of $V$ (it can fail for $N\ge3$ when $V$ has isolated bumps. Although a Choquard analogue of \cite{APS} does not exist, it is interesting to compare our results in presence of \emph{the same} potential $V$ in both sides of the equation, and to observe that for the Choquard equation this seems to improve the situation.

\paragraph{Organization.}
Section~\ref{sec:setting} recalls the unperturbed problem and the uniqueness and nondegeneracy input, and  sets up the perturbed equation. Section~\ref{sec:pge2} describes the standard reduction for $p\ge2$, including positivity; it is included for context and is not the main point of the paper. Section~\ref{sec:psmall} develops the weighted-space theory for $\alpha=2$, $p<2$, and the reduction for a localized perturbation. Section~\ref{sec:semiclassical} treats the semiclassical problem under a small-oscillation condition on $V$. Section~\ref{sec:KV} considers the problem with $V$ on both sides of the equation, for which no small-oscillation condition is needed; its reduction and main results are outlined. 

\section{Preliminaries and the unperturbed problem}\label{sec:setting}

We recall that the Ambrosetti--Badiale method \cite{AB1,AB2} (see also
the monograph \cite{AM1}) is a variational finite-dimensional
reduction for perturbed functionals on a Hilbert space $E$. One
considers --- in the most basic case --- a family
\[
f_\varepsilon(u)=f_0(u)+\varepsilon\,\mathcal{P}(u),\qquad \varepsilon\in\R,
\]
where $f_0$ is twice differentiable on $E$ and has a whole manifold
$Z$ of critical points (typically arising from a symmetry, such as
translation invariance), and $\mathcal P$ is a lower-order
perturbation. The unperturbed critical manifold is assumed to be
\emph{nondegenerate} in the transverse directions: the Hessian
$f_0''(z)$ is invertible on the orthogonal complement of $T_zZ$,
uniformly in $z\in Z$. Under this spectral assumption, a
Lyapunov--Schmidt argument produces, for each $z\in Z$, a correction
$w_\varepsilon(z)\perp T_zZ$ solving the auxiliary equation, and a
\emph{reduced functional}
$\Phi_\varepsilon(z):=f_\varepsilon(z+w_\varepsilon(z))$ on the
finite-dimensional manifold $Z$. Critical points of $\Phi_\varepsilon$
then correspond to genuine critical points of $f_\varepsilon$ near
$Z$; this correspondence is often encoded in a \emph{natural
  constraint} property. The first-order expansion
$\Phi_\varepsilon(z)=c_0+\varepsilon\,\Gamma(z)+o(\varepsilon)$
identifies a \emph{Poincar\'e--Melnikov-type function} $\Gamma$ on
$Z$, and critical points of $\Gamma$ (with nonzero local degree, when
degeneracies are present) persist as critical points of
$\Phi_\varepsilon$ for small $\varepsilon$.

\bigskip

Let us come back to our setting. Let $N\ge 3$ be the space dimension,
$0<\alpha <N$, $p>1$, and let $I_\alpha(x)=c_{N,\alpha}|x|^{\alpha-N}$
be the Riesz potential, so that $I_\beta*I_\gamma=I_{\beta+\gamma}$
whenever $\beta,\gamma>0$ and $\beta+\gamma<N$. We consider the
non-local semilinear equation
\begin{equation}\label{eq:choq}
-\Delta u+u=(I_\alpha*|u|^p)|u|^{p-2}u\quad\text{in }\R^N,
\qquad \frac{N+\alpha}{N}<p<\frac{N+\alpha}{N-2}.
\end{equation}
For a positive solution $u$ to \eqref{eq:choq} we put 
\begin{equation} \label{eq:AV}
A=I_\alpha*u^p, \qquad V=Au^{p-2}.
\end{equation}
The linearization of \eqref{eq:choq} at $u$ is associated to the operator
\[
L_{+,u}\varphi=-\Delta\varphi+\varphi-p\,(I_\alpha*(u^{p-1}\varphi))\,u^{p-1}-(p-1)V\varphi .
\]
Before proceeding further, we collect  a few well-known results to be used in the sequel.
\begin{theorem}[{\cite[Theorem 1.1]{XZZ2}}] \label{thm:2.1}
Every positive $H^1$-solution of \eqref{eq:choq} is, up to translation, radial and strictly decreasing in $|x|$.
\end{theorem}
\begin{theorem}[{\cite[Theorem 1.2]{XZZ2}}] \label{thm:2.2}
Assume
\[
\frac{N+2}{N}<p<\frac{N+2}{N-2}\ \text{ if }\alpha=2,\qquad
2\le p<\frac{N+\alpha}{N-2}\ \text{ if }2<\alpha<N .
\]
Then \eqref{eq:choq} has exactly one positive $H^1$-solution $U$ up to translations, and
\[
i_{-,H^1}(L_{+,U})=1,\qquad \ker L_{+,U}=\operatorname{span} \left\{\partial_{x_1}U,\dots,\partial_{x_N}U \right\}.
\]
\end{theorem}
\begin{definition}
    For each $\xi \in\mathbb{R}^N$, the translation operator $\tau_\xi$ acts on a function $u$ as follows:
    \begin{equation*}
    \tau_\xi u\colon x \mapsto u(x-\xi).
    \end{equation*}
\end{definition}
Here and in the rest of the paper, $U$ denotes this ground state centered at the origin and $U_\xi=U(\cdot-\xi)=\tau_\xi U$. It is easy to check that $U$ is a classical solution. The manifold of nontrivial positive critical points of $f_0$ (defined below) is
\[
Z=\left\{U_\xi \mid  \xi\in\R^N\right\},
\]
and nondegeneracy means $\ker f_0''(U)=T_UZ$.
For $\alpha=2$ and $p<2$ put $\beta:=\frac{N-2}{2-p}$. The following properties (F1) and (F4) are proved in \cite{XZZ1,XZZ2}, while (F2) and (F3) are elementary, and (F5) is classical \cite{Stein}.
\begin{itemize}
\item[(F1)] $U(x)\asymp(1+|x|)^{-\beta}$ as $|x| \to +\infty$.
\item[(F2)] $p\beta>N+2$; in particular $\beta>N$. (Indeed $p\beta>N+2\iff p>1+\frac2N$, and $\beta>N\iff 2-p<\frac{N-2}{N}$.)
\item[(F3)] $(p-1)\beta+(N-2)=\beta$.
\item[(F4)] $V$ is bounded, continuous and $V(x)\to1$ as $|x|\to\infty$.
\item[(F5)] The Bessel kernel $G$ of $(-\Delta+1)^{-1}$ is positive, $\int G=1$, decays exponentially, and $G,\nabla G\in L^1$ (with $|\nabla G|\lesssim|x|^{1-N}$ near $0$). Hence for $h\in L^\infty$, $G*h$ and $\nabla(G*h)=\nabla G*h$ are bounded and uniformly continuous, and $\|G*h\|_{C^1}\le C\|h\|_{L^\infty}$.
\end{itemize}
Finally, 
\begin{center}
for $p\ge2$, $V\to0$ at infinity and the ground state $U$ decays exponentially.
\end{center}


\bigskip

Let $a\in L^\infty(\R^N)$ with $a(x)\to0$ as $|x|\to\infty$ and $a\not\equiv0$. Consider
\begin{equation}\label{eq:pert}
-\Delta u+\left(1+\varepsilon a(x)\right)u=(I_\alpha*|u|^p)|u|^{p-2}u\quad\text{in }\R^N,
\end{equation}
whose solutions are critical points of $f_\varepsilon=f_0+\varepsilon\mathcal P$ on $H^1(\R^N)$, where
\[
f_0(u)=\frac12\|u\|_{H^1}^2-\frac1{2p}\!\int(I_\alpha*|u|^p)|u|^p,\qquad
\mathcal P(u)=\frac12\int a\,u^2 .
\]

\section{The superlienar case \texorpdfstring{$p\ge2$ (and $2\le\alpha<N$)}{p >= 2 (and 2 <= alpha < N)}}\label{sec:pge2}

By Theorems \ref{thm:2.1} and \ref{thm:2.2}, $Z$ is an $N$-dimensional, noncompact, $C^2$ manifold of critical points of $f_0$. The results of \cite{XZZ2} yield $\ker L_{+,U}=T_UZ$. Since $V\to0$ and $U^{p-1}$ decays exponentially, multiplication by $V$ and the nonlocal term (via Hardy--Littlewood--Sobolev) are compact. Hence $L_{+,U}=I-K$ with $K$ compact. 

For $p \geq 2$, $f_0 \in C^2(H^1)$ and $\mathcal{P} \in C^2(H^1)$ with $\mathcal{P}''=a$. By translation invariance $f_0''(U_\xi)=\tau_\xi f_0''(U)\tau_{-\xi}$, so all bounds on inverses and on the reduction are uniform in $\xi$; moreover $\|\mathcal P'(U_\xi)\|\le\|a\|_\infty\|U\|$ uniformly. This is the noncompact perturbation setting presented in \cite{AM1}.

A standard Lyapunov--Schmidt reduction yields $w_\varepsilon(\xi)\perp T_{U_\xi}Z$ with $\|w_\varepsilon\|=O(\varepsilon)$ and the reduced functional
\[
\Phi_\varepsilon(\xi)=f_\varepsilon\left(U_\xi+w_\varepsilon(\xi)\right)
=\frac{p-1}{2p}\|U\|_{H^1}^2+\varepsilon\,\Gamma(\xi)+o(\varepsilon),\qquad
\Gamma(\xi)=\frac12\int a(x)\,U^2(x-\xi)\dd x,
\]
with $o(\varepsilon)$ uniform in $\xi$. If $a\ge0$, $a\not\equiv0$, then $\Gamma>0$ and $\Gamma(\xi)\to0$ as $|\xi|\to\infty$, so $\Gamma$ has a global maximum. Let $s=\operatorname{sgn}\varepsilon$. By the uniform expansion, $s\Phi_\varepsilon=sc_0+|\varepsilon|\Gamma+o(\varepsilon)$ has a maximum over a fixed ball attained at an interior point for $\varepsilon$ small; that is, $\Phi_\varepsilon$ has a maximum point if $\varepsilon>0$ and a minimum point if $\varepsilon<0$. In both cases it is a critical point $\xi_\varepsilon$, and $\xi_\varepsilon\to\bar\xi$ (along a sequence), $\bar\xi$ a global maximum point of $\Gamma$. Thus for $|\varepsilon|$ small there is a solution
\[
u_\varepsilon=\tau_{\xi_\varepsilon}\left(U+w_\varepsilon(\xi_\varepsilon)\right).
\]
If $a\le0$ the same argument applies to $-a$ and $-\varepsilon$; the limit points are then global minimum points of $\Gamma$. More generally, every strict (or $C^1$-stable) critical point of $\Gamma$ generates a solution.

\bigskip

We can also show that this solution $u_\varepsilon$ is positive.
For $p<2$ positivity is part of Proposition~\ref{prop:fixed} below, which provides the pointwise estimate $u\ge\frac12U_\xi$. For $p\ge2$ it needs a more involved argument. In this subsection we write $\|\cdot\|=\|\cdot\|_{H^1}$ to ease notation.

\begin{lemma}[regularity of solutions near $U$]\label{lem:reg-near}
Let $N\ge3$, $\alpha\in(0,N)$ and $2\le p<\frac{N+\alpha}{N-2}$ . For every $\eta>0$ there is $\delta>0$ such that every weak solution $u\in H^1(\R^N)$ of \eqref{eq:pert} with $|\varepsilon|\|a\|_\infty\le\delta$ and $\|u-U_\xi\|_{H^1}\le\delta$ for some $\xi$ is continuous and satisfies $\|u-U_\xi\|_{L^\infty}\le\eta$.
\end{lemma}

\begin{proof}
By translation invariance take $\xi=0$. Suppose the claim fails. Then there are $\eta>0$, $\varepsilon_na_n$ with $\varepsilon_n\|a_n\|_\infty\to0$, and solutions $u_n$ with $u_n\to U$ in $H^1$ but $\|u_n-U\|_\infty\ge\eta$. Put $\Phi_n=I_\alpha*|u_n|^p$.

\emph{Preliminaries.} Let $G=(-\Delta+1)^{-1}$, with kernel $0<G(x)\le C|x|^{2-N}$,
$\int G=1$, and exponentially decaying. Put $b_n=\varepsilon_na_n$, so $\|b_n\|_\infty\le\frac12$ for $n$ large, and
$F_n=(I_\alpha*|u_n|^p)|u_n|^{p-2}u_n$; thus $(-\Delta+1+b_n)u_n=F_n$.
We use two standard facts.
\begin{itemize}
\item[(F1)] If $1<\gamma\le s\le\infty$ and $0\le\frac1\gamma-\frac1s\le\frac2N$, with
$s<\infty$ if equality holds on the right, then $G:L^\gamma\to L^s$ is bounded
(this follows from Young's inequality and Hardy--Littlewood--Sobolev at the endpoint $\frac1\gamma-\frac1s=\frac2N$).
\item[(F2)] If $v\in L^2$ solves $(-\Delta+1+b)v=F$ with $\|b\|_\infty\le\frac12$ and $F\in L^\gamma$,
then $v\in L^s$ for $s$ as in (F1) and $\|v\|_{s}\le2\|GF\|_{s}\le C\|F\|_\gamma$.
Indeed $v+G(bv)=GF$ and $\|G(b\,\cdot)\|_{L^t\to L^t}\le\frac12$ for every $t\in[1,\infty]$. Since $GF=v+G(bv)\in L^2$ and
$GF\in L^s$ by (F1), the Neumann series $\sum_j(-Gb)^jGF$ converges in $L^2$ and in $L^s$; its sum is the
unique $L^2$-solution $v$.
\end{itemize}
All constants below depend only on $N,\alpha,p$ and on upper bounds for the norms already known.

\emph{Step 1: $\Phi_n=I_\alpha*|u_n|^p$ is bounded in $L^\infty$.}
Put $q_f=\frac{2N(p-1)}{\alpha+2}$ and $q_\gamma=\frac{(2p-1)N}{N+\alpha}$. Since $p<\frac{N+\alpha}{N-2}$ we have $q_f<2^*$, and $p<2^*$ because $\alpha<N$.
Moreover $q_\gamma<2^*$: this is equivalent to $(2p-1)(N-2)<2(N+\alpha)$, i.e.\ to $p<\frac{3N+2\alpha-2}{2N-4}$, which holds because $\frac{N+\alpha}{N-2}=\frac{2N+2\alpha}{2N-4}$ is smaller.
Hence we can pick $q_0\in(\max\{p,q_f,q_\gamma\},2^*)$, close to $2^*$. Since $u_n\to U$ in $H^1$, $(u_n)$ is bounded in $L^q$ for $q\in[2,2^*]$.
Suppose $(u_n)$ is bounded in $L^2\cap L^q$ with
\begin{equation}\label{eq:qrange}
\max\{p,q_f,q_\gamma\}<q<\frac{pN}{\alpha}.
\end{equation}
Then $|u_n|^p$ is bounded in $L^{q/p}$ with $1<\frac qp<\frac N\alpha$, so by Hardy--Littlewood--Sobolev
$\Phi_n$ is bounded in $L^\theta$, $\frac1\theta=\frac pq-\frac\alpha N$, and
$F_n$ is bounded in $L^\gamma$ with
\[
\frac1\gamma=\frac1\theta+\frac{p-1}q=\frac{2p-1}q-\frac\alpha N .
\]
Since $q>q_\gamma$ we have $\frac{2p-1}q<1+\frac\alpha N$, i.e.\ $\gamma>1$, so (F1) and (F2) apply.
By (F2) and (F1), $(u_n)$ is bounded in $L^{q'}$ for every $q'$ with $\frac1{q'}\in \left[\max\{0,\frac1\gamma-\frac2N\},\frac1\gamma \right]$.
Taking the best choice,
\begin{equation}\label{eq:qmap}
\frac1{q'}=\max\left\{0,\ \frac{2p-1}q-\frac{\alpha+2}N\right\}\quad(\text{any }q'<\infty\text{ if the maximum is }0),
\end{equation}
and $q'>q$ exactly when $\frac{2p-2}q<\frac{\alpha+2}N$, i.e.\ $q>q_f$, which is part of \eqref{eq:qrange}. Since $q'>q$, the new exponent again satisfies the lower bounds in \eqref{eq:qrange}.
Moreover, by interpolation we may replace $q'$ by any smaller number $>q$, which we use to avoid
the value $\frac{pN}\alpha$.
Define $q_{j+1}$ from $q_j$ by \eqref{eq:qmap} as long as $q_j$ satisfies \eqref{eq:qrange}.
Writing $\frac1{q_f}=\frac{\alpha+2}{2N(p-1)}$, \eqref{eq:qmap} reads (when the maximum is positive)
$\frac1{q_{j+1}}-\frac1{q_f}=(2p-1)\left(\frac1{q_j}-\frac1{q_f}\right)$.
Since $\frac1{q_0}-\frac1{q_f}<0$ and $2p-1>1$, this difference is negative and grows geometrically in absolute value,
so after finitely many steps $q_j\ge\frac{pN}\alpha$ and, avoiding equality, $q_j>\frac{pN}\alpha$.
Then $\sup_n\|u_n\|_{q_j}<\infty$ with $\frac{q_j}p>\frac N\alpha$. Splitting $I_\alpha$ at the interface $|x-y|=1$,
\[
\|\Phi_n\|_\infty\le C\left(\|u_n\|_p^p+\|u_n\|_{q_j}^p\right)\le C,
\]
since $|u_n|^p\in L^1$ (because $p\ge2$) and $|u_n|^p\in L^{q_j/p}$ with $q_j/p>N/\alpha$.

\emph{Step 2: bounds in $L^\infty$.} Now $|F_n|\le C|u_n|^{p-1}$ and $(u_n)$ is bounded in $L^q$ for some
$q>\frac{pN}\alpha$. By (F2) with $\gamma=\frac q{p-1}>1$, $(u_n)$ is bounded in $L^{q'}$ with
$\frac1{q'}=\max\{0,\frac{p-1}q-\frac2N\}$ (and in $L^\infty$ if $\frac{p-1}q<\frac2N$). This improves $q$
iff $q>\frac{N(p-2)}2$. That condition holds for $q>\frac{pN}\alpha$: $\frac{N(p-2)}2<\frac{pN}\alpha$ is equivalent to
$p(\alpha-2)<2\alpha$, which is trivial for $\alpha\le2$, and for $\alpha>2$ follows from
$p<\frac{N+\alpha}{N-2}\le\frac{2\alpha}{\alpha-2}$ (the last inequality is equivalent to $(\alpha-N)(\alpha+2)\le0$).
For $p=2$ each step lowers $\frac1q$ by $\frac2N$; for $p>2$ the quantity $\frac1q-\frac2{N(p-2)}$ is multiplied by $p-1>1$ at each step.
In both cases $\sup_n\|u_n\|_\infty<\infty$ after finitely many steps.

\emph{Step 3: convergence.} Interpolating between $u_n\to U$ in $L^2$ and the uniform $L^\infty$ bound gives $u_n\to U$ in every $L^s$, $s<\infty$. Put $\mathcal N(v)=(I_\alpha*|v|^p)|v|^{p-2}v$ and $w_n=u_n-U$. Then $(-\Delta+1)w_n=\mathcal N(u_n)-\mathcal N(U)-\varepsilon_na_nu_n$. The splitting of $I_\alpha$ at $|x-y|=1$, as in Step~1 of the proof of Proposition~\ref{prop:pos} but using only the $L^s$ convergence just obtained, gives $\Phi_n\to I_\alpha*U^p$ in $L^\infty$, and $|u_n|^{p-2}u_n\to|U|^{p-2}U$ in every $L^q$, $q<\infty$ (since $p\ge2$ and the $u_n$ are bounded in $L^\infty$). So the right-hand side tends to $0$ in every $L^q$, $q<\infty$. Hence $w_n\to0$ in $W^{2,q}$, and in $L^\infty$ by Sobolev's embedding for $q>N/2$, contradicting $\|w_n\|_\infty\ge\eta$.
\end{proof}

\begin{remark}
  The proof of Lemma \ref{lem:reg-near} is inspired by 
  \cite[Proposition~3.2]{MVS3}. We have presented an independent proof that takes into account the zeroth-order
term $\varepsilon au$ and that provides some uniformity of the constants.
\end{remark}

We are ready to state and prove the positivity of the perturbed solution.
\begin{proposition}[positivity for $p\ge2$]\label{prop:pos}
Let $p\ge2$ and $a\in L^\infty(\R^N)$. There is $\delta>0$, depending only on $N$, $\alpha$, $p$, such that the following holds. If $|\varepsilon|\,\|a\|_\infty\le\delta$ and $u\in H^1(\R^N)$ is a weak solution of \eqref{eq:pert} with $\|u-U_\xi\|_{H^1}\le\delta$ for some $\xi\in\R^N$, then $u>0$.
\end{proposition}

\begin{proof}
By translation invariance take $\xi=0$; the constants do not depend on $\xi$. Put $M:=(2\|U\|_\infty)^{p-2}$ and $V_u:=(I_\alpha*|u|^p)|u|^{p-2}\ge0$, and take $\delta\le\frac14$.

\emph{Step 1: the potential is small at infinity.} Since $U^p\in L^1\cap L^\infty$ and $\alpha<N$, we have $I_\alpha*U^p(x)\to0$ as $|x|\to\infty$. Fix $R>0$ with $M\,(I_\alpha*U^p)(x)\le\frac14$ for $|x|\ge R$. Next, for $f=|u|^p-U^p$ we split the kernel at $|x-y|=1$ and use Hölder, which gives for some $r>N/\alpha$
\[
\|I_\alpha*f\|_\infty\le C\left(\|f\|_{L^1}+\|f\|_{L^r}\right).
\]
Since $|f|\le p|u-U|(|u|+U)^{p-1}$ and $p\ge2$, both norms are controlled by $\|u-U\|_{L^s}$ with $s=2$ and $s=pr$, times quantities bounded in terms of $\|u\|_{L^{s'}}$, $s'<\infty$. By Lemma~\ref{lem:reg-near} (applied with $\eta\to0$) and interpolation, $\|u-U\|_{L^s}\to0$ as $\delta\to0$ for every $s\in[2,\infty]$. Hence we may choose $\delta$ so small that $M\|I_\alpha*f\|_\infty\le\frac14$ and $\|u\|_\infty\le2\|U\|_\infty$. For $|x|\ge R$ this gives
\[
V_u(x)\le M\left(I_\alpha*U^p+I_\alpha*f\right)(x)\le\frac12 .
\]

\emph{Step 2: $u>0$ on a large ball.} Shrinking $\delta$ further, Lemma~\ref{lem:reg-near} gives $\|u-U\|_\infty<\min_{\overline{B_R}}U$, so $u>0$ on $\overline{B_R}$ and $\{u<0\}\subset\{|x|>R\}$.

\emph{Step 3: $u\ge0$.} Let $u_-=\min\{u,0\}\in H^1$. Test the equation with $u_-$. Since $\nabla u\cdot\nabla u_-=|\nabla u_-|^2$ and $u\,u_-=u_-^2$,
\[
\int|\nabla u_-|^2+(1+\varepsilon a)u_-^2=\int V_u\,u_-^2 .
\]
As $u_-$ vanishes outside $\{u<0\}\subset\{|x|>R\}$, where $V_u\le\frac12$, and $|\varepsilon a|\le\frac14$, we get
\[
\int|\nabla u_-|^2+\frac14\int u_-^2\le0,
\]
so $u_-=0$ and $u\ge0$. Since $u\not\equiv0$ near $B_R$, $u\not\equiv0$.

\emph{Step 4: strict positivity.} We have $u\ge0$, $u\not\equiv0$, and $-\Delta u+(1+\varepsilon a)u=V_uu\ge0$ weakly with bounded coefficient $1+\varepsilon a$. The weak Harnack inequality for nonnegative supersolutions \cite{GT} gives $u>0$ everywhere.
\end{proof}

\begin{remark}
The argument uses $p\ge2$ only to bound $|u|^{p-2}$, and it needs no spectral information and no restriction on $\alpha$ other than $\alpha<N$. For $p<2$ the factor $|u|^{p-2}$ is singular where $u$ vanishes, and this is why the sublinear case is treated in a different manner.
\end{remark}

\begin{remark}
In the setting of Section~\ref{sec:fixedpoint} the hypothesis of Proposition~\ref{prop:pos} holds for the solutions produced by the reduction, because $\|w_\varepsilon(\xi)\|_{H^1}=O(\varepsilon)$ uniformly in $\xi$. 
\end{remark}

\section{The sublinear case \texorpdfstring{$\alpha=2$, $\frac{N+2}{N}<p<2$}{alpha = 2, (N+2)/N < p < 2}}\label{sec:psmall}

Recall that 
\begin{equation*}
-\Delta u+u=A|u|^{p-2}u, \quad -\Delta A=|u|^p
\end{equation*}
and $U\sim c\,r^{-\beta}$, $V\to1$ at infinity. The linearization is a bounded map $H^1\to H^{-1}$ since $V\in L^\infty$, but $u\mapsto A_u|u|^{p-2}u$ is not expected to be Fr\'echet differentiable at $U$ as a map $H^1\to H^{-1}$.

\begin{example}
Take $\varphi_R=\kappa\,U(R)\,b(\cdot-x_R)$ with $|x_R|=R$, $b$ a fixed unit bump and $\kappa>0$ fixed. On its support the remainder of $h(u)=|u|^{p-2}u$ satisfies
\[
A\left[h(U+\varphi)-h(U)-h'(U)\varphi\right]\approx c_\kappa\,V\,U,\qquad
c_\kappa=(1+\kappa)^{p-1}-1-(p-1)\kappa<0
\]
(by the strict concavity of $t\mapsto t^{p-1}$ for $p<2$). Its $H^{-1}$ norm is $\gtrsim U(R)$, while $\|\varphi_R\|_{H^1}\sim\kappa U(R)\to0$. Hence
\[
\frac{\|\text{remainder}\|_{H^{-1}}}{\|\varphi_R\|_{H^1}}\to\frac{|c_\kappa|}{\kappa}\ne0\quad(R\to\infty),
\]
so $f_0$ is not twice Fr\'echet differentiable at $U$. 
\end{example}

Fortunately, the effective mass at infinity is $1-(p-1)V\to2-p>0$, and \cite{XZZ2} proves nondegeneracy there. The obstruction is thus one of regularity only.

\subsection{Reduction in a weighted space \texorpdfstring{$Y\subset H^1$}{Y in H1}}
To overcome this lack of differentiability, we introduce
\[
X=\left\{\varphi\in L^\infty(\R^N) \mid \|\varphi\|_X=\operatorname*{ess\,sup}_{x \in\mathbb{R}^N} (1+|x|)^{\beta}|\varphi(x)|<\infty\right\},
\]
\[
Y=\left\{\varphi\in X \mid \nabla\varphi\in L^2(\R^N)\right\},\qquad \|\varphi\|_Y=\|\varphi\|_X+\|\nabla\varphi\|_{L^2}.
\]
Since $\beta>N$, $X\subset L^q$ for every $q\in[1,\infty]$ with $\|\varphi\|_{L^q}\le C_q\|\varphi\|_X$; in particular $\|\varphi\|_{L^2}\lesssim\|\varphi\|_X$ and $Y\subset H^1$ with $\|\varphi\|_{H^1}\lesssim\|\varphi\|_Y$. Note that elements of $X$ need not belong to $H^1$, which is why we work in $Y$.

\begin{lemma}
The space $Y$ is a Banach space.
\end{lemma}
\begin{proof}
It is easy to prove that
\begin{equation*}
 L^\infty \left(\mathbb{R}^N,1+|\cdot |^\beta \right) = \left\lbrace u \mid \operatorname{ess\, sup}_{x \in\mathbb{R}^N} (1+|x|^\beta) |u(x)| < +\infty \right\rbrace
\end{equation*}
is complete with respect to the norm $\Vert \cdot \Vert_X$. By a general result of Functional Analysis, the cartesian product $D^{1,2}(\mathbb{R}^N) \times L^\infty(\mathbb{R}^N,1+|\cdot |^\beta)$ is a Banach space, and $Y$ is isometric to its diagonal, which is a closed subspace. Hence $Y$ is complete.
\end{proof}

We write the equation \eqref{eq:choq} as a fixed point problem $u=T_\varepsilon(u)$ with 
\begin{equation*}
T_\varepsilon(u)=G*\left(A_u|u|^{p-2}u-\varepsilon a\,u\right), \quad A_u=I_2*|u|^p.
\end{equation*}

The linearization of $T_0$ at $U$ is
\begin{equation}\label{eq:Tprime}
T_0'(U)\varphi=(p-1)G*\varphi+\underbrace{(p-1)G*[(V-1)\varphi]}_{K_1\varphi}
+\underbrace{p\,G*\left[U^{p-1}\,(I_2*(U^{p-1}\varphi))\right]}_{K_2\varphi}.
\end{equation}

\subsection{The Fredholm property on \texorpdfstring{$Y$}{Y}}\label{sec:fredholm}




Define the bounded operator $\mathbb T$ on $Y$ by the formula \eqref{eq:Tprime}:
\[
\mathbb T\varphi=(p-1)G*\varphi+K_1\varphi+K_2\varphi
=G*\left[(p-1)V\varphi+p\,U^{p-1}\left(I_2*(U^{p-1}\varphi)\right)\right].
\]
We shall see in Proposition~\ref{prop:T0} that $\mathbb T=T_0'(U)$. Anyway, the proof of the next theorem only uses the formula.

\begin{theorem}\label{thm:fred}
Let $N\ge3$, $\alpha=2$, $1+\frac2N<p<\min \left\{2,\frac{N+2}{N-2} \right\}$, let $U$ be the positive radial ground state, and $\beta=\frac{N-2}{2-p}$. Then $\mathbb T$ maps $Y$ into $Y$ boundedly, and $I-\mathbb T$ is Fredholm of index $0$ on $Y$.
\end{theorem}

The proof of this result will be accomplished in several steps. We start by collecting a few useful facts.
\begin{enumerate}[label=(\alph*),nosep]
\item $U\le C(1+|x|)^{-\beta}$, $V\in L^\infty$ and $V(x)\to1$ as $|x|\to\infty$ (see \cite[eqs.~(2.9), (2.11), (2.12)]{XZZ2}).
\item $p\beta>N+2$ and $(p-1)\beta+(N-2)=\beta$. Also $\beta>N$, because $p>1+\frac2N$ gives $2-p<\frac{N-2}{N}$.
\item The Bessel kernel $G$ is positive with $\int G=1$. It satisfies $G\lesssim|x|^{2-N}$ and $|\nabla G|\lesssim|x|^{1-N}$ near $0$, and both decay exponentially at infinity. Moreover $\|G*h\|_{W^{2,q}}\le C_q\|h\|_{L^q}$ for $1<q<\infty$ (\cite[Lemma~2.1]{XZZ2}).
\end{enumerate}
We also write $w_s(x)=(1+|x|)^{-s}$, so that $\|\varphi\|_X=\operatorname*{ess\,sup}w_\beta^{-1}|\varphi|$ and $\|\varphi\|_Y=\|\varphi\|_X+\|\nabla\varphi\|_{L^2}$.

Let us begin with some auxiliary lemmas.

%
%

\begin{lemma}[Convolution with finite measures] \label{lem:convfin}
Let $\mu$ be a finite signed Borel measure on $\mathbb R^N$ and
let $s\geq 0$ be such that
\[
  M_s(\mu):=
  \int_{\mathbb R^N}(1+|z|)^s\dd|\mu|(z)<\infty,
\]
where $|\mu|$ denotes the total variation measure of $\mu$.
Then the convolution operator
\[
  T_\mu\phi:=\mu*\phi
\]
is bounded on the weighted space
\[
  X_s=\left\{\phi\in L^\infty_{\mathrm{loc}}(\mathbb R^N)\mid 
  \|\phi\|_{X_s}:=
  \operatorname*{ess\,sup}_{x\in\mathbb R^N}
  (1+|x|)^s|\phi(x)|<\infty\right\},
\]
and
\begin{equation} \label{eq:4.2}
  \|\mu*\phi\|_{X_s}
  \leq M_s(\mu)\|\phi\|_{X_s}.
 \end{equation}
Moreover, if $\phi\in H^1(\mathbb R^N)$, then
\begin{equation} \label{eq:4.3}
  \|\nabla(\mu*\phi)\|_{L^2}
  \leq |\mu|(\mathbb R^N)\|\nabla\phi\|_{L^2}.
\end{equation}
Consequently, convolution with $\mu$ defines a bounded
linear operator on
\[
  Y=X_{\beta-1}\cap
  \{\phi\in H^1(\mathbb R^N) \mid \nabla\phi\in L^2\},
\]
equipped with the norm
\[
  \|\phi\|_Y=\|\phi\|_{X_{\beta-1}}
             +\|\nabla\phi\|_{L^2}.
\]
In particular,
\[
  \|\mu*\phi\|_Y
  \leq
  \left(M_{\beta-1}(\mu)+|\mu|(\mathbb R^N)\right)
  \|\phi\|_Y.
\]
\end{lemma}

\begin{proof}
Since
\[
  1+|x|\leq (1+|x-y|)(1+|y|),
\]
we have
\[
  (1+|x-y|)^{-s}
  \leq (1+|y|)^s(1+|x|)^{-s}.
\]
Therefore,
\begin{align*}
  |\mu*\phi(x)|
  &\leq \int_{\mathbb R^N}|\phi(x-y)|\dd|\mu|(y)\\
  &\leq \|\phi\|_{X_s}(1+|x|)^{-s}
       \int_{\mathbb R^N}(1+|y|)^s\dd|\mu|(y),
\end{align*}
which proves \eqref{eq:4.2}.

For the gradient estimate, differentiation commutes with
convolution in the sense of distributions:
\[
  \nabla(\mu*\phi)=\mu*\nabla\phi.
\]
The latter convolution can be understood as a Bochner integral
with values in $L^2(\mathbb R^N)$. Minkowski's inequality and
the invariance of the $L^2$ norm under translations give
\begin{align*}
  \|\nabla(\mu*\phi)\|_{L^2}
  &\leq
  \int_{\mathbb R^N}
  \|\nabla\phi(\,\cdot-y)\|_{L^2}\dd|\mu|(y)\\
  &=|\mu|(\mathbb R^N)\|\nabla\phi\|_{L^2}.
\end{align*}
Combining the two estimates proves \eqref{eq:4.3}.
\end{proof}

\begin{lemma}[regularity]\label{lem:reg}
Let $s>N/2$ and $|h|\le Mw_s$ a.e. Then $v=G*h$ satisfies
\[
|v|+|\nabla v|\le CMw_s,\qquad \|\nabla v\|_{L^2}\le CM,\qquad \|v\|_{C^{1,\gamma}}\le CM
\]
for each fixed $\gamma\in(0,1)$. If $s\ge\beta$, then $v\in Y$ and $\|v\|_Y\le CM$.
\end{lemma}
\begin{proof}
%
%
By the exponential decay of the Bessel kernel and its gradient,
both $G$ and $\nabla G$ belong to $L^1(\mathbb R^N)$ and have
finite moments of every order. Hence the weighted convolution
estimate of Lemma \ref{lem:convfin} applies to both kernels. In particular,
if $|h|\leq M w_s$ with $s>N/2$, then
\[
  |G*h|+|\nabla G*h|\leq CMw_s.
\]
Moreover, Young's inequality gives
\[
  \|\nabla(G*h)\|_{L^2}
  \leq \|\nabla G\|_{L^1}\|h\|_{L^2}.
\]
Finally, since $(-\Delta+1)(G*h)=h$ and $h\in L^q$
for every $q\geq2$, elliptic regularity yields
\[
  \|G*h\|_{W^{2,q}}\leq C_q\|h\|_{L^q}.
\]
Choosing $q>N$ sufficiently large and applying Morrey's
embedding proves the asserted $C^{1,\gamma}$ estimate.

 Next, $\|\nabla v\|_2\le\|\nabla G\|_1\|h\|_2\le CM\|w_s\|_2$, which is finite since $2s>N$. Finally $h\in L^q$ for all $q\ge2$ (it is bounded and lies in $L^2$). Take $q>N$ with $1-N/q\ge\gamma$: then $\|v\|_{W^{2,q}}\lesssim\|h\|_q\lesssim M$, and Morrey's embedding gives the $C^{1,\gamma}$ bound.
\end{proof}

We can now start the proof of Theorem \ref{thm:fred}.
\paragraph{Step 1: $\mathbb T$ maps $Y$ into $Y$ boundedly.}
Indeed, $(p-1)G*$ is bounded on $Y$ by Lemma \ref{lem:convfin}. 
Then $|(V-1)\varphi|\le(\|V\|_\infty+1)\|\varphi\|_Xw_\beta$, so Lemma~\ref{lem:reg} gives $\|K_1\varphi\|_Y\lesssim\|\varphi\|_X$.

Put now $g_\varphi=U^{p-1}\varphi$. By (a), $|g_\varphi|\lesssim\|\varphi\|_Xw_{p\beta}$ with $p\beta>N$. Splitting $\int|x-y|^{2-N}w_{p\beta}(y)\dd y$ into $\{|y|<|x|/2\}$ and $\{|y|\ge|x|/2\}$ gives $I_2*|g_\varphi|\lesssim\|\varphi\|_Xw_{N-2}$. By (a) and (b), $U^{p-1}w_{N-2}\lesssim w_{(p-1)\beta+N-2}=w_\beta$. So $|U^{p-1}(I_2*g_\varphi)|\lesssim\|\varphi\|_Xw_\beta$ and Lemma~\ref{lem:reg} applies.


\paragraph{Step 2: the principal part is an isomorphism of $Y$.}
Put $\mu=2-p>0$, so $\mu+(p-1)=1$. Using the Fourier transform normalization for which the resolvent kernels have multipliers $\widehat G(\zeta)=(1+|\zeta|^2)^{-1}$ and $\widehat G_\mu(\zeta)=(\mu+|\zeta|^2)^{-1}$, we have
\[
\widehat{G*G_\mu}=\frac1{(1+|\zeta|^2)(\mu+|\zeta|^2)}=\frac{\hat G_\mu-\hat G}{p-1}.
\]
The two sides of this identity are in $L^1$ and have the same Fourier transform; uniqueness of the Fourier transform on $L^1$ therefore gives $(p-1)\,G*G_\mu=G_\mu-G$. Denoting by $\delta$  the Dirac mass at the origin,
\[
\left(\delta-(p-1)G\right)*\left(\delta+(p-1)G_\mu\right)
=\delta+(p-1)G_\mu-(p-1)G-(p-1)(G_\mu-G)=\delta .
\]


Let
\[
  \mu_{\mathcal{D}}:=\delta_0-(p-1)G(x)\dd x,
  \qquad
  \mu_S:=\delta_0+(p-1)G_\mu(x)\dd x.
\]
Since the Bessel kernels $G$ and $G_\mu$ decay exponentially
at infinity, both measures have finite moments of every
nonnegative order. In particular,
\[
  \int_{\mathbb R^N}(1+|x|)^\beta\dd|\mu_D|(x)
  <\infty,
  \qquad
  \int_{\mathbb R^N}(1+|x|)^\beta\dd|\mu_S|(x)
  <\infty.
\]
By the preceding lemma, convolution with either measure
defines a bounded linear operator on $Y$. Thus
\[
  \mathcal{D}
  =I-(p-1)G*,
  \qquad
  S
  =I+(p-1)G_\mu*
\]
are bounded operators on $Y$.

Furthermore, the resolvent identity
\[
  (p-1)G*G_\mu=G_\mu-G
\]
implies, in the convolution algebra of finite signed measures,
\begin{align*}
  \mu_{\mathcal{D}}*\mu_S
  &=\left(\delta_0-(p-1)G\right)
    *\left(\delta_0+(p-1)G_\mu\right)\\
  &=\delta_0+(p-1)G_\mu-(p-1)G
    -(p-1)^2G*G_\mu\\
  &=\delta_0.
\end{align*}
Since convolution of finite measures is associative and
commutative, we also have
\[
  \mu_S*\mu_{\mathcal{D}}=\delta_0.
\]
The convolution operators therefore satisfy
\[
  \mathcal{D}S=S\mathcal{D}=I\quad\text{on }Y.
\]
Hence $\mathcal{D}$ is an isomorphism of $Y$, with bounded inverse
\[
  \mathcal{D}^{-1}=S=I+(p-1)G_\mu*.
\]

\begin{lemma}[compactness]\label{lem:A}
Let $\theta>0$ and $\mathcal H=\{h\text{ measurable}\mid  |h|\le w_{\beta+\theta}\text{ a.e.}\}$. Then $G*\mathcal H$ is precompact in $Y$.
\end{lemma}
\begin{proof}
Let $v_j=G*h_j$ with $h_j\in\mathcal H$. By Lemma~\ref{lem:reg} (with $s=\beta+\theta$), $\|v_j\|_{C^{1,\gamma}}\le C$ and $|v_j|+|\nabla v_j|\le Cw_{\beta+\theta}$. By Arzel\`a--Ascoli on each ball $\overline{B_R}$ and a diagonal argument, a subsequence satisfies $v_j\to v$ and $\nabla v_j\to\nabla v$ locally uniformly.

Given $\epsilon>0$, choose $R$ with $C(1+R)^{-\theta}<\epsilon$ and $C^2\int_{|x|>R}w_{2\beta+2\theta}<\epsilon^2$; the latter choice is possible since $2\beta>N$. Then for all $j,k$,
\[
\sup_{|x|>R}w_\beta^{-1}|v_j-v_k|\le2\epsilon,\qquad \int_{|x|>R}|\nabla v_j-\nabla v_k|^2\le4\epsilon^2,
\]
while on $B_R$ both $\sup|v_j-v_k|$ and $\sup|\nabla v_j-\nabla v_k|$ tend to $0$ as $j,k\to\infty$. So $(v_j)$ is a Cauchy sequence  in $Y$, and the proof is complete.
\end{proof}

\paragraph{Step 3: $K_1$ is compact on $Y$.}
Let $\chi_R$ be a cutoff with $0\le\chi_R\le1$, $\chi_R=1$ on $B_R$ and $\chi_R=0$ outside $B_{2R}$, and split $K_1=K_1^{(R)}+E_R$ with $K_1^{(R)}\varphi=(p-1)G*[(V-1)\chi_R\varphi]$.
\begin{itemize}[nosep]
\item For $\|\varphi\|_Y\le1$, $|(V-1)\chi_R\varphi|\le C_R\mathbf 1_{B_{2R}}\le C_R'\,w_{\beta+1}$. By Lemma~\ref{lem:A}, $K_1^{(R)}$ maps the unit ball into a precompact set, so it is compact.
\item Let $h=(p-1)(V-1)(1-\chi_R)\varphi$ and $\epsilon_R=\operatorname*{ess\,sup}_{|x|>R}|V-1|$, so $\epsilon_R\to0$ by (a). Then $|h|\le(p-1)\epsilon_R\|\varphi\|_Xw_\beta$. By (W), $\|G*h\|_X\le C\epsilon_R\|\varphi\|_X$. By (Y), $\|\nabla G*h\|_2\le\|\nabla G\|_1\|h\|_2\le C\epsilon_R\|\varphi\|_X\|w_\beta\|_2$. So $\|E_R\|_{Y\to Y}\le C\epsilon_R\to0$.
\end{itemize}
Thus $K_1$ is a norm limit of compact operators, hence compact.

\paragraph{Step 4: $K_2$ is compact on $Y$.}
The leading tail of $I_2*g_\varphi$ is $m(\varphi)I_2(x)$ with $m(\varphi)=\int g_\varphi$. By (a), $|m(\varphi)|\le\|\varphi\|_X\int U^{p-1}w_\beta<\infty$, so $m$ is a bounded linear functional on $Y$. Let $\iota\in C^\infty(\R^N)$ with $\iota=I_2$ on $\{|x|\ge1\}$, and put
\[
\rho\varphi:=U^{p-1}\left(I_2*g_\varphi-m(\varphi)\iota\right),\qquad
\psi_0:=G*(U^{p-1}\iota).
\]
Then
\[
K_2\varphi=p\,m(\varphi)\,\psi_0+p\,G*\rho\varphi .
\]

\begin{claim}\label{cl:rho}
$|\rho\varphi|\le C\|\varphi\|_Xw_{\beta+1}$.
\end{claim}
\begin{proof}
Let $|x|\ge2$. Then $\iota(x)=I_2(x)$, and $I_2*g_\varphi-m\iota=\int\left(I_2(x-y)-I_2(x)\right)g_\varphi(y)\dd y$. Split the latter integral.
\begin{enumerate}[label=(\arabic*),nosep]
\item $|y|<|x|/2$. For $t\in[0,1]$, $|x-ty|\ge|x|/2$, so $|I_2(x-y)-I_2(x)|\lesssim|y||x|^{1-N}$. The contribution is $\lesssim\|\varphi\|_X|x|^{1-N}\int|y|w_{p\beta}$, finite since $p\beta>N+1$.
\item $|y|\ge|x|/2$, term $\int I_2(x-y)g_\varphi(y)\dd y$. On $\{|x-y|<|x|/2\}$ we have $|g_\varphi|\lesssim\|\varphi\|_X|x|^{-p\beta}$ and $\int_{|z|<|x|/2}|z|^{2-N}\dd z\lesssim|x|^2$. On $\{|x-y|\ge|x|/2\}$ we have $|I_2(x-y)|\lesssim|x|^{2-N}$ and $\int_{|y|\ge|x|/2}|g_\varphi|\lesssim\|\varphi\|_X|x|^{N-p\beta}$. Both give $|x|^{2-p\beta}$.
\item $|y|\ge|x|/2$, term $I_2(x)\int_{|y|\ge|x|/2}g_\varphi$. This is $\lesssim\|\varphi\|_X|x|^{2-N}|x|^{N-p\beta}=\|\varphi\|_X|x|^{2-p\beta}$.
\end{enumerate}
Since $p\beta\ge N+1$, the total is $\lesssim\|\varphi\|_X|x|^{1-N}$. Multiplying by $U^{p-1}\lesssim|x|^{-(p-1)\beta}$ and using $(p-1)\beta+(N-1)=\beta+1$ gives $|\rho\varphi|\lesssim\|\varphi\|_X|x|^{-\beta-1}$ for $|x|\ge2$. For $|x|<2$ everything is bounded by $C\|\varphi\|_X$, because $g_\varphi\in L^1\cap L^\infty$ makes $I_2*g_\varphi$ bounded.
\end{proof}

By Claim~\ref{cl:rho} and Lemma~\ref{lem:A} (with $\theta=1$), $\varphi\mapsto G*\rho\varphi$ is compact on $Y$. Moreover $|U^{p-1}\iota|\lesssim w_\beta$, so $\psi_0\in Y$ by Lemma~\ref{lem:reg}, and $\varphi\mapsto m(\varphi)\psi_0$ is bounded of rank one. Hence $K_2$ is compact.

\paragraph{Step 5: conclusion.}
By Step 2, $\mathcal D^{-1}=S$ is bounded on $Y$. Hence
\[
I-\mathbb T=\mathcal D-K_1-K_2=\mathcal D\,\left(I-S(K_1+K_2)\right),
\]
and $C:=S(K_1+K_2)$ is compact on $Y$ by Steps 3--4. By the Riesz--Schauder theory, $I-C$ has finite-dimensional kernel and closed range of codimension $\dim\ker(I-C)$, so it is Fredholm of index $0$. Composition with the isomorphism $\mathcal D$ preserves both properties, so $I-\mathbb T$ is Fredholm of index $0$ on $Y$. \qed


\subsection{Regularity of \texorpdfstring{$f_0$ and $T_0$ on $Y$}{f0 and T0 on Y}}\label{sec:regularity}

In the basic Ambrosetti-Badiale scheme, the functional $f_0$ must be at least of class $C^2$ in the whole Hilbert space. In our setting $f_0$ is not $C^2$ on $H^1$, and yet we can prove that it is locally smooth on a whole neighborhood of $U$ in $Y$ with $f_0''$ locally Lipschitz, and that $T_0$ is $C^{1,1}$ from a small $X$-ball into $Y$. In order to perform the reduction, this regularity is sufficient.

Recall that we are assuming
\begin{displaymath}
N\ge3, \quad \alpha=2, \quad 1+\frac2N<p<\min \left\{2,\frac{N+2}{N-2} \right\}, \quad \beta=\frac{N-2}{2-p}.
\end{displaymath}
 The symbol $C$ denotes constants depending only on $N,p$, on $U$ (through the quantities $c_U,\|V\|_\infty,\|U^{p-1}\|_{L^N}$ below) and on $\delta_0$.

\bigskip

Besides (a)--(c) above we observe that:
\begin{enumerate}[label=(\alph*$'$),nosep]
\item $c_U:=\inf_{x \in\mathbb{R}^N}U(x)(1+|x|)^{\beta}>0$. Indeed, by \cite[(2.12)]{XZZ2}, $U(x) \asymp(1+|x|)^{-\beta}$ for $|x|$ large, and $U$ is positive and continuous. Thus $w_\beta\le U/c_U$ and $U\le C w_\beta$.
\item $A:=I_2*U^p$ is finite, $V=AU^{p-2}\in L^\infty$, and $\int AU^p<\infty$ (see \cite[Lemmas~2.3 and 4.1]{XZZ2}).
\item $w_\beta\in L^2$ (since $2\beta>N$), hence $VU\le Cw_\beta\in L^2$; and $U^{p-1}\in L^N$, because $(p-1)\beta=\beta-(N-2)>2$ and so $\int U^{(p-1)N}\lesssim\int w_{(p-1)\beta N}<\infty$.
\end{enumerate}
Fix $\delta_0>0$ with $\delta_0/c_U\le\frac12$ and set
\[
\mathcal W_X=\{u=U+\varphi \mid \varphi\in X,\ \|\varphi\|_X<\delta_0\},\qquad
\mathcal W_Y=\mathcal W_X\cap(U+Y).
\]
Since $\|\cdot\|_X\le\|\cdot\|_Y$, $\mathcal W_Y$ is an open subset of the affine space $U+Y$. The set $\mathcal{W}_X$ is a tubular neighborhood of $U$ in $X$, and $\mathcal{W}_Y$ is a tubular neighborhood of $U$ in $Y$. The following lemma collects some useful estimates that will be used repeatedly in the sequel.
For $u\in\mathcal W_X$ put $A_u=I_2*u^p$ and $V_u=A_uu^{p-2}$.

\begin{lemma}[pointwise ratio estimates]\label{lem:key}
There is a positive constant $C$ such that for all $u,u'\in\mathcal W_X$ and $\psi\in X$ the following hold.
\begin{enumerate}[label=(\roman*),nosep]
\item $\frac12U\le u\le\frac32U$.
\item For $s\in\{p,\,p-1,\,p-2\}$: $|u^s-u'^s|\le C\,U^s\,\|u-u'\|_X$. If moreover $u+\psi\in\mathcal W_X$, then for $q\in\{p,p-1\}$: $\left|(u+\psi)^q-u^q-q\,u^{q-1}\psi\right|\le C\,U^q\,\|\psi\|_X^2$.
\item $A_u\le CA$, $\ |A_u-A_{u'}|\le C\|u-u'\|_X\,A$, and $I_2*(U^{p-1}|\psi|)\le C\|\psi\|_X\,A$.
\item $V_u\le CV$ and $\|V_u-V_{u'}\|_{L^\infty}\le C\|u-u'\|_X$.
\item $AU^{p-1}=VU\le Cw_\beta$; in particular $AU^{p-1}\in L^2$.
\end{enumerate}
\end{lemma}
\begin{proof}
(i) Write $u=U+\varphi$. Since $|\varphi|\le\|\varphi\|_Xw_\beta\le\delta_0U/c_U\le U/2$, we get $|u/U-1|\le\frac12$.

(ii) By (i), $u,u'$ and every point between them lie in $[U/2,3U/2]$ pointwise. For real $s$, $|\frac{d}{dt}t^s|=|s|t^{s-1}\le|s|\kappa_sU^{s-1}$ on $[U/2,3U/2]$, with $\kappa_s=\max\{2^{1-s},(3/2)^{s-1}\}$. So $|u^s-u'^s|\le C U^{s-1}|u-u'|\le CU^{s-1}\|u-u'\|_Xw_\beta\le CU^s\|u-u'\|_X/c_U$. For the second-order bound, Taylor's formula with Lagrange remainder gives $(u+\psi)^q-u^q-qu^{q-1}\psi=\frac12q(q-1)\xi^{q-2}\psi^2$ with $\xi$ between $u$ and $u+\psi$, hence $\xi\in[U/2,3U/2]$ and $\xi^{q-2}\le CU^{q-2}$; and $\psi^2\le\|\psi\|_X^2w_\beta^2\le\|\psi\|_X^2U^2/c_U^2$.

(iii) By (i), $u^p\le CU^p$, so $A_u\le CA$. Applying $I_2*$ to (ii) with $s=p$ gives $|A_u-A_{u'}|\le C\|u-u'\|_XA$. Finally $U^{p-1}|\psi|\le\|\psi\|_XU^{p-1}w_\beta\le\|\psi\|_XU^p/c_U$.

(iv) $V_u\le CA\cdot CU^{p-2}=CV$, where $u^{p-2}\le2^{2-p}U^{p-2}$ because $p-2<0$. Next $V_u-V_{u'}=(A_u-A_{u'})u^{p-2}+A_{u'}(u^{p-2}-u'^{p-2})$. By (iii) and (ii) with $s=p-2$, the first term is $\le C\|u-u'\|_XAU^{p-2}$ and the second is $\le CA\,U^{p-2}\|u-u'\|_X$; both equal $C\|u-u'\|_XV\le C\|V\|_\infty\|u-u'\|_X$.

(v) is immediate from (a$'$)--(c$'$).
\end{proof}

The next auxiliary result will help us prove the regularity of $f_0$. Its proof is standard, and we propose it for the reader's sake.
\begin{lemma}[a criterion for $C^1$]\label{lem:BS}
Let $E$ be a Banach space, $Z$ a normed space, and $\mathcal N\subset Z^*$ a set of functionals with $\|\ell\|\le1$ for $\ell\in\mathcal N$ and $\|z\|=\sup_{\ell\in\mathcal N}|\ell(z)|$ for all $z\in Z$. Let $\mathcal O\subset E$ be open and $\Phi\colon \mathcal O\to Z$. Suppose that for every $u\in\mathcal O$ there is $D\Phi(u)\in L(E,Z)$ such that $u\mapsto D\Phi(u)$ is continuous for the operator norm, and such that for all $\psi\in E$ and $\ell\in\mathcal N$ the function $t\mapsto\ell(\Phi(u+t\psi))$ is differentiable (wherever $u+t\psi\in\mathcal O$) with derivative $\ell\left(D\Phi(u+t\psi)\psi\right)$. Then $\Phi\in C^1(\mathcal O;Z)$ with Fr\'echet derivative $D\Phi$.
\end{lemma}
\begin{proof}
Fix $u\in\mathcal O$ and $h\in E$ small, and $\ell\in\mathcal N$. The scalar function $t\mapsto\ell(\Phi(u+th))$ has the continuous derivative $t\mapsto\ell(D\Phi(u+th)h)$ on $[0,1]$ (continuity from the norm-continuity of $D\Phi$). By the scalar fundamental theorem of calculus,
\[
\ell\left(\Phi(u+h)-\Phi(u)-D\Phi(u)h\right)=\int_0^1\ell\left((D\Phi(u+th)-D\Phi(u))h\right)\dd t,
\]
which is bounded by $\sup_{t}\|D\Phi(u+th)-D\Phi(u)\|\,\|h\|$. Taking the supremum over $\ell\in\mathcal N$ gives $\|\Phi(u+h)-\Phi(u)-D\Phi(u)h\|\le o(1)\|h\|$.
\end{proof}
\noindent We use two cases: $Z=\R$, $\mathcal N=\{\mathrm{id}\}$ (scalar functionals), and $Z=Y^*$ with $\mathcal N=\{\eta\mapsto(\cdot)(\eta):\ \eta\in Y,\ \|\eta\|_Y\le1\}$ (evaluations), which norm $Y^*$.

\begin{proposition}[$f_0\in C^{2,1}_{\mathrm loc}$ on $Y$]\label{prop:f0}
The functional $f_0$ is of class $C^2$ on the open set $\mathcal W_Y\subset U+Y$, with, for $u\in\mathcal W_Y$ and $\psi,\eta\in Y$,
\begin{align*}
f_0'(u)\eta&=\langle u,\eta\rangle_{H^1}-\int A_u\,u^{p-1}\eta,\\
f_0''(u)[\psi,\eta]&=\langle\psi,\eta\rangle_{H^1}-(p-1)\int V_u\,\psi\eta-p\int\left(I_2*(u^{p-1}\psi)\right)u^{p-1}\eta ,
\end{align*}
and for all $u,u'\in\mathcal W_Y$,
\begin{align}
\|f_0'(u)-f_0'(u')\|_{Y^*}&\le\|u-u'\|_{H^1}+C\|u-u'\|_X\ \le\ C\|u-u'\|_Y,\label{eq:f1lip}\\
\|f_0''(u)-f_0''(u')\|_{L(Y,Y^*)}&\le C\|u-u'\|_X\ \le\ C\|u-u'\|_Y .\label{eq:f2lip}
\end{align}
In \eqref{eq:f2lip} the bilinear forms are even bounded on $L^2\times L^2$. In particular $f_0''(U)=L_{+,U}$.
\end{proposition}
\begin{proof}
\emph{Step 1: $f_0$ is well defined and G\^ateaux differentiable.} For $u\in\mathcal W_Y$ we have $u\in H^1$ and $u^p\le CU^p$ by Lemma~\ref{lem:key}(i), so $\mathcal F(u):=\iint I_2(x-y)u^p(x)u^p(y)\dd x\dd y\le C\int AU^p<\infty$; thus $f_0(u)=\frac12\|u\|_{H^1}^2-\frac1{2p}\mathcal F(u)$ is finite. Fix $u\in\mathcal W_Y$ and $\eta\in Y$, and let $v_t=u+t\eta$ for $|t|\le t_0$, with $t_0$ so small that $v_t\in\mathcal W_Y$. The $t$-derivative of the integrand of $\mathcal F(v_t)$ is
\[
I_2(x-y)\left[p\,v_t^{p-1}(x)\eta(x)\,v_t^p(y)+v_t^p(x)\,p\,v_t^{p-1}(y)\eta(y)\right],
\]
whose absolute value is, uniformly in $|t|\le t_0$, at most
\[
C\,I_2(x-y)\left[U^{p-1}|\eta|(x)U^p(y)+U^p(x)U^{p-1}|\eta|(y)\right].
\]
This is integrable on $\mathbb R^N\times\mathbb R^N$, with integral $2C\int AU^{p-1}|\eta|\le2C\|AU^{p-1}\|_{2}\|\eta\|_2<\infty$ by Lemma~\ref{lem:key}(v). Differentiation under the integral sign and the symmetry of $I_2$ give $\frac{d}{dt}\mathcal F(v_t)|_{t=0}=2p\int A_uu^{p-1}\eta$, which is the stated formula for $f_0'(u)\eta$. It is linear in $\eta$, and $|f_0'(u)\eta|\le\|u\|_{H^1}\|\eta\|_{H^1}+C\|VU\|_2\|\eta\|_2\le C(u)\|\eta\|_Y$, so $f_0'(u)\in Y^*$.

\emph{Step 2: estimate \eqref{eq:f1lip} and $f_0\in C^1$.} Write $A_uu^{p-1}-A_{u'}u'^{p-1}=(A_u-A_{u'})u^{p-1}+A_{u'}(u^{p-1}-u'^{p-1})$. By Lemma~\ref{lem:key}(i)--(iii) both terms are bounded by $C\|u-u'\|_XAU^{p-1}=C\|u-u'\|_XVU$. Hence $\left|\int(A_uu^{p-1}-A_{u'}u'^{p-1})\eta\right|\le C\|u-u'\|_X\|VU\|_2\|\eta\|_2$, and adding $|\langle u-u',\eta\rangle_{H^1}|\le\|u-u'\|_{H^1}\|\eta\|_{H^1}$ gives \eqref{eq:f1lip} (recall $\|\cdot\|_{H^1}\lesssim\|\cdot\|_Y$ and $\|\cdot\|_2\lesssim\|\cdot\|_Y$). So $u\mapsto f_0'(u)\in Y^*$ is continuous, and Lemma~\ref{lem:BS} (with $\Phi=f_0$, $E=Y$, $Z=\R$) shows $f_0\in C^1(\mathcal W_Y)$.

\emph{Step 3: derivative of $f_0'$.} Fix $u\in\mathcal W_Y$ and $\psi,\eta\in Y$, and let $v_t=u+t\psi\in\mathcal W_Y$ for $|t|\le t_0$. We differentiate $t\mapsto f_0'(v_t)\eta=\langle v_t,\eta\rangle_{H^1}-\iint I_2(x-y)v_t^p(y)\,v_t^{p-1}(x)\eta(x)\dd y\dd x$. The first term has derivative $\langle\psi,\eta\rangle_{H^1}$. For the second, the $t$-derivative of the integrand is
\[
I_2(x-y)\,p\,v_t^{p-1}(y)\psi(y)\,v_t^{p-1}(x)\eta(x)\;+\;I_2(x-y)v_t^p(y)\,(p-1)v_t^{p-2}(x)\psi(x)\eta(x).
\]
Uniformly in $|t|\le t_0$, the first piece is bounded by $C\,I_2(x-y)U^{p-1}|\psi|(y)U^{p-1}|\eta|(x)$, whose integral is finite by the Hardy--Littlewood--Sobolev inequality (see Step 4); the second piece is bounded by $CA(x)U^{p-2}(x)|\psi||\eta|(x)=CV|\psi||\eta|$ after integrating in $y$, which is in $L^1$ since $V\in L^\infty$ and $\psi,\eta\in L^2$. Differentiating under the integral sign gives, at $t=0$, exactly $f_0''(u)[\psi,\eta]$ as stated (using $A_uu^{p-2}=V_u$). The map $\psi\mapsto f_0''(u)[\psi,\cdot]$ is a bounded linear map $Y\to Y^*$ by the estimates of Step 4. This is the differentiability required in Lemma~\ref{lem:BS} for $\Phi=f_0'$ with the evaluation functionals $\ell_\eta$, $\eta\in Y$.

\emph{Step 4: estimate \eqref{eq:f2lip}.} Subtract the formulas at $u$ and $u'$; the $H^1$ term cancels.
\begin{itemize}[nosep]
\item By Lemma~\ref{lem:key}(iv), $\left|\int(V_u-V_{u'})\psi\eta\right|\le C\|u-u'\|_X\|\psi\|_2\|\eta\|_2$.
\item For the nonlocal term put $\Theta_u(\psi,\eta)=\int(I_2*(u^{p-1}\psi))u^{p-1}\eta$. Then
\[
\Theta_u-\Theta_{u'}=\int\left(I_2*((u^{p-1}-u'^{p-1})\psi)\right)u^{p-1}\eta+\int\left(I_2*(u'^{p-1}\psi)\right)(u^{p-1}-u'^{p-1})\eta .
\]
By Lemma~\ref{lem:key}(i),(ii) with $s=p-1$, $|u^{p-1}-u'^{p-1}|\le C\|u-u'\|_XU^{p-1}$ and $u^{p-1},u'^{p-1}\le CU^{p-1}$. Hence $|\Theta_u-\Theta_{u'}|\le C\|u-u'\|_X\,\Lambda(|\psi|,|\eta|)$ with $\Lambda(f,g)=\int\left(I_2*(U^{p-1}f)\right)U^{p-1}g$. By the Hardy--Littlewood--Sobolev inequality with exponent $r=\frac{2N}{N+2}$ (note $\frac1r+\frac1r+\frac{N-2}{N}=2$) and H\"older's inequality ($\frac1r=\frac1N+\frac12$),
\[
\Lambda(f,g)\le C\|U^{p-1}f\|_r\|U^{p-1}g\|_r\le C\|U^{p-1}\|_{L^N}^2\|f\|_2\|g\|_2 ,
\]
and $U^{p-1}\in L^N$ by (c$'$).
\end{itemize}
This proves \eqref{eq:f2lip} (on $L^2\times L^2$, hence on $Y\times Y$), and also the boundedness of $f_0''(u)$ used in Step 3 (the $H^1$ term is bounded on $Y\times Y$ and the others by the same estimates with $u'=U$). By Step 3 and \eqref{eq:f2lip}, Lemma~\ref{lem:BS} applies to $\Phi=f_0'\colon \mathcal W_Y\to Y^*$, so $f_0'\in C^1$ and $f_0\in C^2(\mathcal W_Y)$. At $u=U$ we have $A_U=A$, $V_U=V$, so the formula is the form of $L_{+,U}$. \qedhere
\end{proof}


\begin{proposition}[$T_0\in C^{1,1}$ from an $X$-ball into $Y$]\label{prop:T0}
Let $\mathsf n(u)=A_uu^{p-1}$ and $T_0(u)=G*\mathsf n(u)$ for $u\in\mathcal W_X$. Then $T_0$ maps $\mathcal W_X$ into $Y$, is Fr\'echet differentiable as a map from the open subset $\mathcal W_X$ of $U+X$ into $Y$, and its derivative is
\[
T_0'(u)\psi=G*\mathcal L_u\psi,\qquad
\mathcal L_u\psi=(p-1)V_u\psi+p\,u^{p-1}\,I_2*(u^{p-1}\psi),
\]
with $\|T_0'(u)\|_{L(X,Y)}\le C$ and
\begin{equation}\label{eq:T0lip}
\|T_0'(u)-T_0'(u')\|_{L(X,Y)}\le C\|u-u'\|_X\qquad(u,u'\in\mathcal W_X).
\end{equation}
Consequently $T_0\in C^{1,1}(\mathcal W_X;Y)$, and the restriction of $T_0$ to $\mathcal W_Y$ is $C^{1,1}$ from $(\mathcal W_Y,\|\cdot\|_Y)$ into $Y$, with $T_0'(U)=\mathbb T$ on $Y$.
\end{proposition}
\begin{proof}
\emph{Mapping and boundedness.} For $u\in\mathcal W_X$, Lemma~\ref{lem:key}(i),(iii),(v) give $0<\mathsf n(u)\le CAU^{p-1}\le Cw_\beta$. By Lemma~\ref{lem:reg} with $s=\beta$, $T_0(u)\in Y$ and $\|T_0(u)\|_Y\le C$. Likewise, for $\psi\in X$: $|V_u\psi|\le CV\|\psi\|_Xw_\beta\le C\|\psi\|_Xw_\beta$, and by Lemma~\ref{lem:key}(i),(iii), $|u^{p-1}I_2*(u^{p-1}\psi)|\le CU^{p-1}\cdot C\|\psi\|_XA=C\|\psi\|_XVU\le C\|\psi\|_Xw_\beta$. So $|\mathcal L_u\psi|\le C\|\psi\|_Xw_\beta$, and Lemma~\ref{lem:reg} gives $\|G*\mathcal L_u\psi\|_Y\le C\|\psi\|_X$.

\emph{Differentiability.} Let $u\in\mathcal W_X$ and $\psi\in X$ with $u+\psi\in\mathcal W_X$. Write $P_1=A_u$, $P_2=u^{p-1}$, $\Delta P_1=A_{u+\psi}-A_u$, $\Delta P_2=(u+\psi)^{p-1}-u^{p-1}$, and note $\mathsf n'(u)\psi=(P_1'\psi)P_2+P_1(P_2'\psi)=\mathcal L_u\psi$ with $P_1'\psi=I_2*(pu^{p-1}\psi)$, $P_2'\psi=(p-1)u^{p-2}\psi$ (recall $A_uu^{p-2}=V_u$). Then
\[
\mathsf n(u+\psi)-\mathsf n(u)-\mathcal L_u\psi=(\Delta P_1-P_1'\psi)\,P_2+P_1\,(\Delta P_2-P_2'\psi)+\Delta P_1\,\Delta P_2 .
\]
By Lemma~\ref{lem:key}(ii),(iii): $|\Delta P_1-P_1'\psi|=\left|I_2*\left[(u+\psi)^p-u^p-pu^{p-1}\psi\right]\right|\le C\|\psi\|_X^2A$; $|\Delta P_2-P_2'\psi|\le C\|\psi\|_X^2U^{p-1}$; $|\Delta P_1|\le I_2*(CU^{p-1}|\psi|)\le C\|\psi\|_XA$; $|\Delta P_2|\le C\|\psi\|_XU^{p-1}$; and $P_2\le CU^{p-1}$, $P_1\le CA$. Hence the remainder is bounded by $C\|\psi\|_X^2\,AU^{p-1}\le C\|\psi\|_X^2\,w_\beta$, and Lemma~\ref{lem:reg} yields
\[
\|T_0(u+\psi)-T_0(u)-T_0'(u)\psi\|_Y\le C\|\psi\|_X^2 .
\]
So $T_0$ is Fr\'echet differentiable $\mathcal W_X\to Y$ with derivative $T_0'(u)$.

\emph{Lipschitz derivative.} For $\psi\in X$,
\begin{multline*}
(\mathcal L_u-\mathcal L_{u'})\psi \\ 
=(p-1)(V_u-V_{u'})\psi+p\,(u^{p-1}-u'^{p-1})\,I_2*(u^{p-1}\psi)+p\,u'^{p-1}\,I_2*\left((u^{p-1}-u'^{p-1})\psi\right).
\end{multline*}
By Lemma~\ref{lem:key}(iv) the first term is $\le C\|u-u'\|_X\|\psi\|_Xw_\beta$. By (i)--(iii) the second is $\le CU^{p-1}\|u-u'\|_X\cdot C\|\psi\|_XA$ and the third is $\le CU^{p-1}\cdot I_2*(CU^{p-1}\|u-u'\|_X|\psi|)\le CU^{p-1}\|u-u'\|_X\|\psi\|_XA$; both are $\le C\|u-u'\|_X\|\psi\|_XVU\le C\|u-u'\|_X\|\psi\|_Xw_\beta$. Lemma~\ref{lem:reg} then gives \eqref{eq:T0lip}.

A Fr\'echet differentiable map with Lipschitz derivative is $C^1$ (indeed $C^{1,1}$). The statement for $\mathcal W_Y$ follows by composing with the bounded inclusion $Y\hookrightarrow X$ ($\|\cdot\|_X\le\|\cdot\|_Y$), and $T_0'(U)=G*\mathcal L_U=\mathbb T$ by the definition \eqref{eq:Tprime}. \qedhere
\end{proof}

\begin{lemma}[gradient identity]\label{lem:grad}
For $u\in\mathcal W_Y$ and $\psi,\eta\in Y$,
\[
f_0'(u)\eta=\langle u-T_0(u),\eta\rangle_{H^1},\qquad
f_0''(u)[\psi,\eta]=\left\langle(I-T_0'(u))\psi,\eta\right\rangle_{H^1}.
\]
Thus $\nabla_{H^1}f_0=\mathrm{id}-T_0$ on $\mathcal W_Y$, and $I-T_0'(u)$ is symmetric for the $H^1$ inner product.
\end{lemma}
\begin{proof}
For $h\in L^2$, $G*h\in H^2$ and $(-\Delta+1)(G*h)=h$ (\cite[Lemma~2.1]{XZZ2} with $r=2$), and $\mathsf n(u),\mathcal L_u\psi\in L^2$ by the bounds above. For $\eta\in H^1$ an integration by parts gives $\langle G*h,\eta\rangle_{H^1}=\int h\eta$. Taking $h=\mathsf n(u)$ gives $\langle T_0(u),\eta\rangle_{H^1}=\int A_uu^{p-1}\eta$, hence the first identity. Taking $h=\mathcal L_u\psi$ gives $\langle T_0'(u)\psi,\eta\rangle_{H^1}=(p-1)\int V_u\psi\eta+p\int\left(I_2*(u^{p-1}\psi)\right)u^{p-1}\eta$, hence the second.
\end{proof}

\begin{remark}\label{rem:unif}
The constants above depend only on $N,p,\delta_0$ and the translation-invariant quantities $c_U,\|V\|_\infty,\|U^{p-1}\|_{L^N}$, $\|w_\beta\|_{L^2}$. Replacing $U$ by $U_\xi$ and the weight by $w_\beta(\cdot-\xi)$ (norms $X_\xi,Y_\xi$), all statements hold with the same constants, since $f_0$ and $T_0$ commute with translations. For $f_\varepsilon=f_0+\varepsilon\mathcal P$ with $\mathcal P(u)=\frac12\int au^2$, $a\in L^\infty$, the functional $\mathcal P$ is $C^2$ on $H^1$ with $\mathcal P''(u)[\psi,\eta]=\int a\psi\eta$ (independent of $u$), so Proposition~\ref{prop:f0} holds for $f_\varepsilon$ with $f_\varepsilon''=f_0''+\varepsilon a$. For $T_\varepsilon=T_0-\varepsilon G*(a\,\cdot)$ the extra term is a bounded linear map $X\to Y$ (by Lemma~\ref{lem:reg}, since $|a\varphi|\le\|a\|_\infty\|\varphi\|_Xw_\beta$), so Proposition~\ref{prop:T0} holds for $T_\varepsilon$ with $T_\varepsilon'(u)=T_0'(u)-\varepsilon G*(a\,\cdot)$.
\end{remark}

\subsection{Invertibility on the orthogonal complement of the tangent space}\label{sec:LS}

Following the idea of the Ambrosetti-Badiale reduction, we want to build an orthogonal deformation of the manifold of critical points for $f_0$. Since we cannot directly refer to the abstract construction, we set up a fixed-point argument that produces, at the end, the desired orthogonal deformation. We isolate some useful ingredients.

\begin{remark} \label{rem:kernel}
 If $\varphi\in Y$ and $\varphi=T_0'(U)\varphi$, applying $-\Delta+1$ gives $L_{+,U}\varphi=0$ weakly with $\varphi\in Y\subset H^1$, so the nondegeneracy theorem of \cite{XZZ2} gives $\varphi\in\operatorname{span}\{\partial_iU\}$. Conversely, $\partial_iU\in Y$: since $U=G*(VU)$ with $|VU|\le Cw_\beta$, Lemma~\ref{lem:reg} gives $|\nabla U|\le Cw_\beta$, i.e.\ $\partial_iU\in X$, and $\nabla\partial_iU\in L^2$ because $U\in H^2$ (indeed $H^3$, Proposition~\ref{prop:H3}).
\end{remark}

\begin{remark}\label{rem:range}
By the gradient identity $\nabla f_0(u)=u-T_0(u)$, the operator $I-T_0'(U)$ is symmetric for the $H^1$ inner product: $\langle(I-T_0'(U))\varphi,\psi\rangle_{H^1}=\langle L_{+,U}\varphi,\psi\rangle$. Hence for $\psi=(I-T_0'(U))\varphi$ with $\varphi\in Y$ one has $\langle\psi,\partial_iU\rangle_{H^1}=\langle\varphi,(I-T_0'(U))\partial_iU\rangle_{H^1}=0$. So the range lies in $W:=\{\psi\in Y:\ \langle\psi,\partial_iU\rangle_{H^1}=0,\ i=1,\dots,N\}$. The $N$ functionals defining $W$ are independent on $Y$ (their Gram matrix on $\partial_1U,\dots,\partial_NU\in Y$ is $\|\partial_1U\|_{H^1}^2\,\mathrm{Id}$, by radial symmetry), so $W$ has codimension $N$. By Theorem~\ref{thm:fred} the operator has index $0$ and kernel of dimension $N$, so its range has codimension $N$ and therefore equals $W$.
\end{remark}
The previous two remarks show that the operator $I-T_0'(U)$ is Fredholm of index $0$ on $Y$, with kernel $T_UZ=\operatorname{span}\{\partial_iU\}$ and range $Y\cap(T_UZ)^\perp$. Hence
\[
\Pi^\perp_U\left(I-T_0'(U)\right) \colon \ Y\cap(T_UZ)^\perp\longrightarrow Y\cap(T_UZ)^\perp
\]
is an isomorphism, which is the invertibility the Lyapunov--Schmidt reduction needs.

With these facts a contraction argument gives $w_\varepsilon(\xi)\perp_{H^1}T_{U_\xi}Z$ with $\|w_\varepsilon\|_Y=O(\varepsilon)$, uniformly in $\xi$ by translation invariance, using Proposition~\ref{prop:T0} for the regularity of $T_0$ and Remark~\ref{rem:unif} for the uniformity in $\xi$.

\subsection{The Lyapunov--Schmidt equation: existence of the correction}\label{sec:fixedpoint}

Throughout this subsection and the following ones, $N\ge3$, $\alpha=2$, $1+\frac2N<p<\min\{2,\frac{N+2}{N-2}\}$, and $X,Y,\mathcal W_X,\mathcal W_Y,\delta_0,c_U$ are as above.

\paragraph{Assumption on the perturbation.}
\begin{itemize}
\item[\textbf{(A)}] $a\in C^1(\R^N)$, with $a$ and $\nabla a$ bounded, $\nabla a$ uniformly continuous, and $a(x)\to0$ as $|x|\to\infty$ (for instance $a\in C^1_c$).
\end{itemize}
The existence of the correction below uses only $a\in L^\infty$; the rest of (A) is used for $C^1$ dependence on $\xi$ (Proposition~\ref{prop:C1}) and in Section~\ref{sec:critpt}. The perturbed equation is
\begin{equation}\label{eq:pertp<2}
-\Delta u+(1+\varepsilon a)u=(I_2*u^p)u^{p-1},
\end{equation}
whose solutions are the fixed points of $T_\varepsilon(u)=T_0(u)-\varepsilon G*(au)$.

\paragraph{Notation.} Put $\tau_\xi v=v(\cdot-\xi)$ and $a_\xi=a(\cdot+\xi)$. Since $\langle\partial_iU,\partial_jU\rangle_{H^1}=\kappa\delta_{ij}$ with $\kappa=\|\partial_1U\|_{H^1}^2>0$ (radial symmetry), define
\[
Y_\perp=\left\{w\in Y\mid \langle w,\partial_iU\rangle_{H^1}=0,\ i=1,\dots,N\right\},\qquad
\Pi^\perp\psi=\psi-\kappa^{-1}\sum_{i=1}^N\langle\psi,\partial_iU\rangle_{H^1}\,\partial_iU .
\]
Since $\partial_iU\in Y$ and $\|\psi\|_{H^1}\lesssim\|\psi\|_Y$, $\Pi^\perp$ is a bounded projection of $Y$ onto the closed subspace $Y_\perp$, which is therefore a Banach space (it is the space $W$ of the paragraph ``Range'' above). Let
\[
L:=I-T_0'(U)=I-\mathbb T\in L(Y),\qquad Q(w):=T_0(U+w)-T_0(U)-T_0'(U)w\quad(w\in Y,\ \|w\|_X<\delta_0).
\]

\begin{lemma}\label{lem:Lperp}
$L$ maps $Y$ into $Y_\perp$, $\ker L=\operatorname{span}\{\partial_iU\}$, and $L_\perp:=L|_{Y_\perp}:Y_\perp\to Y_\perp$ is an isomorphism. Moreover $\Pi^\perp L=L$ on $Y$. Set $M:=\|L_\perp^{-1}\|_{L(Y_\perp)}\,\|\Pi^\perp\|_{L(Y)}$.
\end{lemma}
\begin{proof}
The first two assertions are Remarks \ref{rem:kernel} and \ref{rem:range}. $L_\perp$ is injective since $\ker L\cap Y_\perp=\{0\}$ (its Gram matrix is $\kappa\,\mathrm{Id}$). The operator is also onto: $Y=Y_\perp\oplus\operatorname{span}\{\partial_iU\}$ and $L$ vanishes on the second summand, so $Y_\perp=L(Y)=L(Y_\perp)$. By the open mapping theorem $L_\perp^{-1}$ is bounded. Since $L(Y)\subset Y_\perp$, $\Pi^\perp L=L$.
\end{proof}

\begin{lemma}\label{lem:Q}
There are $C_Q,C_G>0$ such that
\begin{enumerate}[label=(\roman*),nosep]
\item $Q(0)=0$, $Q\in C^1$ on $\{w\in Y \mid \|w\|_X<\delta_0\}$, $\|Q'(w)\|_{L(Y)}\le C_Q\|w\|_X$, and for $\|w\|_X,\|w'\|_X<\delta_0$,
\begin{gather*}
\|Q(w)-Q(w')\|_Y\le C_Q\left(\|w\|_X+\|w'\|_X\right)\|w-w'\|_Y,
\end{gather*}
and in particular $\|Q(w)\|_Y\le C_Q\|w\|_X\|w\|_Y$.
\item for $b\in L^\infty$ and $\varphi\in X$ there results $\|G*(b\varphi)\|_Y\le C_G\|b\|_\infty\|\varphi\|_X$.
\end{enumerate}
\end{lemma}
\begin{proof}
(ii) is Lemma~\ref{lem:reg} with $s=\beta$ and $M=\|b\|_\infty\|\varphi\|_X$, since $|b\varphi|\le Mw_\beta$.
(i) By Proposition~\ref{prop:T0}, $T_0\in C^{1}$ and, for $u=U+w$, $Q'(w)=T_0'(U+w)-T_0'(U)$, whose norm in $L(X,Y)$ is $\le C\|w\|_X$ by \eqref{eq:T0lip}; restricting to $Y\subset X$ ($\|\cdot\|_X\le\|\cdot\|_Y$) gives the bound on $\|Q'(w)\|_{L(Y)}$ (this defines $C_Q$). For the Lipschitz estimate write $Q(w)-Q(w')=\int_0^1Q'(w'+t(w-w'))(w-w')\dd t$ and use $\|w'+t(w-w')\|_X\le\|w\|_X+\|w'\|_X$ (the set $\{\|w\|_X<\delta_0\}$ is convex).
\end{proof}

\bigskip

We can now solve the fixed-point equation on which the finite-dimensional reduction is based. Because $f_0$, $T_0$ commute with translations and $a\,\tau_\xi\tilde v=\tau_\xi(a_\xi\tilde v)$, we have $T_\varepsilon(\tau_\xi\tilde v)=\tau_\xi\left[T_0(\tilde v)-\varepsilon G*(a_\xi\tilde v)\right]$. Write $u=\tau_\xi(U+w)$ with $w\in Y_\perp$. As the $H^1$ inner product is translation invariant, the Lyapunov--Schmidt equation 
\begin{displaymath}
  \Pi^\perp_\xi(u-T_\varepsilon u)=0
\end{displaymath} 
becomes, in the fixed space $Y_\perp$,
\begin{equation}\label{eq:LS}
\mathfrak S(\varepsilon,\xi,w):=\Pi^\perp\left[w-\left(T_0(U+w)-T_0(U)\right)+\varepsilon\,G*\left(a_\xi(U+w)\right)\right]=0,
\end{equation}
using $T_0(U)=U$. The parameter $\xi$ enters only through $a_\xi$. Since $w-(T_0(U+w)-T_0(U))=Lw-Q(w)$ and $\Pi^\perp L=L$, equation \eqref{eq:LS} is equivalent to the fixed point problem
\begin{equation}\label{eq:fixed}
w=\mathcal N_{\varepsilon,\xi}(w):=L_\perp^{-1}\Pi^\perp\left[Q(w)-\varepsilon G*(a_\xi U)-\varepsilon G*(a_\xi w)\right].
\end{equation}

\begin{proposition}[existence of the correction, uniform in $\xi$]\label{prop:fixed}
Let $a\in L^\infty$ and $c_1:=C_G\|a\|_\infty(1+\|U\|_X)$. Fix $\rho_0\in(0,\delta_0)$ with $MC_Q\rho_0\le\frac18$. There is $\varepsilon_0>0$, depending only on $M,C_Q,c_1,\rho_0$, such that for every $|\varepsilon|\le\varepsilon_0$ and \emph{every} $\xi\in\R^N$, $\mathcal N_{\varepsilon,\xi}$ has a unique fixed point $\tilde w_\varepsilon(\xi)$ in $B_{\rho_0}=\{w\in Y_\perp:\|w\|_Y\le\rho_0\}$, and
\[
\|\tilde w_\varepsilon(\xi)\|_Y\le2Mc_1|\varepsilon|\le\rho_0/2 .
\]
In addition $U+\tilde w_\varepsilon(\xi)\ge\frac12U>0$.
\end{proposition}
\begin{proof}
Choose $\varepsilon_0$ so that $Mc_1\varepsilon_0\le\frac18$, $Mc_1\varepsilon_0(1+\rho_0)\le\frac{\rho_0}2$ and $2Mc_1\varepsilon_0\le\frac{\rho_0}{2}$. By Lemma~\ref{lem:Q}(ii), $\|G*(a_\xi U)\|_Y\le c_1$ and $\|G*(a_\xi w)\|_Y\le c_1\|w\|_Y$ (the constant does not depend on $\xi$ because $\|a_\xi\|_\infty=\|a\|_\infty$), and $\|w\|_X\le\|w\|_Y\le\rho_0<\delta_0$ on $B_{\rho_0}$.

\emph{Self-map.} For $w\in B_{\rho_0}$, by Lemma~\ref{lem:Q}(i) with $w'=0$,
\[
\|\mathcal N(w)\|_Y\le M\left[C_Q\rho_0^2+|\varepsilon|c_1(1+\rho_0)\right]\le\frac{\rho_0}8+\frac{\rho_0}2<\rho_0 .
\]
\emph{Contraction.} For $w,w'\in B_{\rho_0}$,
\[
\|\mathcal N(w)-\mathcal N(w')\|_Y\le M\left[2C_Q\rho_0+|\varepsilon|c_1\right]\|w-w'\|_Y\le\left(\frac14+\frac18\right)\|w-w'\|_Y .
\]
Banach's theorem on the complete metric space $B_{\rho_0}$ gives existence and uniqueness. For the size: with $\tilde w=\mathcal N(\tilde w)$ and $\|\tilde w\|_Y\le\rho_0$,
\[
\|\tilde w\|_Y\le M\left[C_Q\|\tilde w\|_Y^2+|\varepsilon|c_1(1+\|\tilde w\|_Y)\right]\le\frac18\|\tilde w\|_Y+\frac18\|\tilde w\|_Y+Mc_1|\varepsilon|,
\]
so $\|\tilde w\|_Y\le\frac43Mc_1|\varepsilon|\le2Mc_1|\varepsilon|$. Finally $|\tilde w|\le\|\tilde w\|_Xw_\beta\le(\|\tilde w\|_X/c_U)\,U\le(\delta_0/c_U)U\le\frac12U$.
\end{proof}

\begin{remark}
For $\varepsilon=0$ the fixed point is $\tilde w_0=0$, i.e.\ $U$ itself. The $O(\varepsilon)$ bound is the quantitative form of ``the perturbed solution is a small deformation of $U_\xi$'', and the uniformity in $\xi$ is exactly what the noncompactness of $Z$ requires: it comes from translation invariance and from the fact that the perturbation enters \eqref{eq:fixed} only through $\|a_\xi\|_\infty$.
\end{remark}

\subsection{\texorpdfstring{$C^1$ dependence on $(\varepsilon,\xi)$}{C1 dependence on (epsilon, xi)}}\label{sec:C1dep}

\begin{lemma}\label{lem:xi}
Under (A), the map $\xi\mapsto\left(\varphi\mapsto G*(a_\xi\varphi)\right)$ is of class $C^1$ from $\R^N$ into $L(X,Y)$, with $\partial_{\xi_j}G*(a_\xi\varphi)=G*((\partial_ja)_\xi\varphi)$.
\end{lemma}
\begin{proof}
For $h\in\R^N$, $a_{\xi+h}-a_\xi-h\cdot(\nabla a)_\xi=\int_0^1h\cdot\left((\nabla a)_{\xi+th}-(\nabla a)_\xi\right)\dd t$, whose sup norm is $o(|h|)$ uniformly in $\xi$ because $\nabla a$ is uniformly continuous. So $\xi\mapsto a_\xi\in L^\infty$ is $C^1$ with derivative $(\nabla a)_\xi$ (continuous in $\xi$ by the same uniform continuity). Compose with the bounded linear map $b\mapsto(\varphi\mapsto G*(b\varphi))$ from $L^\infty$ to $L(X,Y)$ (Lemma~\ref{lem:Q}(ii)).
\end{proof}

\begin{proposition}\label{prop:C1}
Assume (A). After decreasing $\varepsilon_0$, the map $(\varepsilon,\xi)\mapsto\tilde w_\varepsilon(\xi)\in Y_\perp$ is of class $C^1$ on $(-\varepsilon_0,\varepsilon_0)\times\R^N$.
\end{proposition}
\begin{proof}
$\mathfrak S$ is $C^1$ on $(-\varepsilon_0,\varepsilon_0)\times\R^N\times\{w\in Y_\perp:\|w\|_Y<\rho_0\}$: in $w$ by Proposition~\ref{prop:T0} (and linearity of the last term), in $\varepsilon$ because it is affine, in $\xi$ by Lemma~\ref{lem:xi}; the partial derivatives are jointly continuous. At a solution $(\varepsilon,\xi,\tilde w)$, since $\Pi^\perp L=L$,
\[
D_w\mathfrak S=L_\perp-\Pi^\perp Q'(\tilde w)+\varepsilon\,\Pi^\perp G*(a_\xi\,\cdot)
=L_\perp\left[I-L_\perp^{-1}\Pi^\perp\left(Q'(\tilde w)-\varepsilon G*(a_\xi\,\cdot)\right)\right].
\]
By Lemma~\ref{lem:Q} and $\|\tilde w\|_Y\le2Mc_1|\varepsilon|$, the operator inside the bracket has norm smaller than $Mc_1|\varepsilon|(2MC_Q+1)$, which in turn is smaller than $\frac12$ once $\varepsilon_0$ is small (depending only on $M,C_Q,c_1$). Hence $D_w\mathfrak S$ is invertible, and the implicit function theorem gives a $C^1$ branch of solutions near each $(\varepsilon,\xi)$. Since $\|\tilde w_\varepsilon(\xi)\|_Y\le\rho_0/2$, this branch stays in $B_{\rho_0}$ by continuity, so by the uniqueness in Proposition~\ref{prop:fixed} it coincides with $\tilde w_\varepsilon(\xi)$.
\end{proof}

\subsection{The Melnikov function and its expansion}\label{sec:reduced}

For $|\varepsilon|<\varepsilon_0$ define
\[
\Phi_\varepsilon(\xi):=f_\varepsilon\left(\tau_\xi(U+\tilde w_\varepsilon(\xi))\right)=\Psi_\varepsilon\left(\xi,\tilde w_\varepsilon(\xi)\right),\qquad
\Psi_\varepsilon(\xi,w)=f_0(U+w)+\frac\varepsilon2\int a_\xi(U+w)^2 ,
\]
where we used the translation invariance of $f_0$ and the change of variables $x\mapsto x+\xi$ in the last integral. By Proposition~\ref{prop:f0} and Lemma~\ref{lem:xi}, $\Psi_\varepsilon$ is $C^1$ on $\R^N\times\{w\in Y:\|w\|_Y<\delta_0\}$, hence by Proposition~\ref{prop:C1} $\Phi_\varepsilon\in C^1(\R^N)$.

\begin{proposition}[uniform expansion]\label{prop:expansion}
Under (A) with $a\in L^\infty$, for $|\varepsilon|<\varepsilon_0$ and all $\xi\in\R^N$,
\[
\Phi_\varepsilon(\xi)=\frac{p-1}{2p}\|U\|_{H^1}^2+\varepsilon\,\Gamma(\xi)+r_\varepsilon(\xi),\quad
\Gamma(\xi)=\frac12\int a(x)\,U^2(x-\xi)\dd x,\quad \sup_{\xi \in \mathbb{R}^N} |r_\varepsilon(\xi)|\le C\varepsilon^2 .
\]
\end{proposition}
\begin{proof}
Write $w=\tilde w_\varepsilon(\xi)$, $\|w\|_Y\le2Mc_1|\varepsilon|$. By Proposition~\ref{prop:f0}, Taylor's formula with integral remainder, $f_0'(U)=0$ and Lemma~\ref{lem:grad},
\[
f_0(U+w)=f_0(U)+\frac12\langle Lw,w\rangle_{H^1}+\mathrm{Rem}(w),
\]
\[
|\mathrm{Rem}(w)|\le\int_0^1(1-t)\left\|f_0''(U+tw)-f_0''(U)\right\|_{L(Y,Y^*)}\|w\|_Y^2\dd t\le C\|w\|_X\|w\|_Y^2\le C\|w\|_Y^3,
\]
where the middle bound is \eqref{eq:f2lip}. Further $f_0(U)=\frac12\|U\|_{H^1}^2-\frac1{2p}\int AU^p=\frac{p-1}{2p}\|U\|_{H^1}^2$, because $\|U\|_{H^1}^2=\int AU^p$ (test the equation with $U$). Finally
\[
\frac\varepsilon2\int a_\xi(U+w)^2=\varepsilon\Gamma(\xi)+\varepsilon\int a_\xi Uw+\frac\varepsilon2\int a_\xi w^2 ,
\]
and $\left|\int a_\xi Uw\right|\le\|a\|_\infty\|U\|_{L^2}\|w\|_{L^2}\le C\|w\|_Y$, $\int a_\xi w^2\le C\|w\|_Y^2$. Collecting: $r_\varepsilon=\frac12\langle Lw,w\rangle+\mathrm{Rem}(w)+\varepsilon\int a_\xi Uw+\frac\varepsilon2\int a_\xi w^2=O(\varepsilon^2)$ uniformly in $\xi$.
\end{proof}

\subsection{Critical points of \texorpdfstring{$\Phi_\varepsilon$}{Phi} solve the equation (natural constraint)}\label{sec:natural}

Fix $|\varepsilon|<\varepsilon_0$ and $\xi$, and write $v=U+\tilde w_\varepsilon(\xi)$ and
\[
g:=v-T_0(v)+\varepsilon\,G*(a_\xi v)\in Y .
\]
Equation \eqref{eq:LS} says $\Pi^\perp g=0$, so
\[
g=\sum_{i=1}^Nc_i\,\partial_iU,\qquad c=c(\varepsilon,\xi)\in\R^N .
\]
The function $u=\tau_\xi v$ solves \eqref{eq:pertp<2} if and only if $g=0$, i.e.\ $c=0$.

The natural-constraint argument we are going to prove 
differentiates the equation twice in the direction of translations, and for this it needs $\partial_iU\in H^2$, i.e.\ $U\in H^3$. We could not find an explicit proof in the literature, so we prove it from (F1)--(F5).

\begin{proposition}\label{prop:H3}
Let $N\ge3$, $\alpha=2$, $1+\frac2N<p<\min \left\{2,\frac{N+2}{N-2} \right\}$. Then $U\in C^1(\R^N)\cap H^3(\R^N)$. In particular $\partial_iU\in H^2$ and $\partial_i\partial_jU\in H^1$ for all $i,j$.
\end{proposition}
\begin{proof}
  Put $f:=AU^{p-1}=VU$, so that $-\Delta U+U=f$ weakly with
  $U\in H^1$. Throughout, $w_s(x)=(1+|x|)^{-s}$.

\emph{Step 1: $U\in C^1\cap H^2$.} By (F1) and (F4),
$0<f=VU\le Cw_\beta$; thus $f$ is bounded and lies in $L^2$ (as
$2\beta>N$). Since $U\in H^1$ solves $-\Delta U+U=f$, we get
$\hat U=\hat f/(1+|\zeta|^2)$ and so $U\in H^2$. Also $U=G*f$, and
(F5) gives $U\in C^1$. (Alternatively use Lemma~\ref{lem:reg}.)
Moreover $U$ is bounded and $U>0$.

\emph{Step 2: $A\in C^1$ and $|\nabla A(x)|\le C(1+|x|)^{-(N-1)}$.} Write $I_2=c|x|^{2-N}$. Since $N\ge3$, $I_2\in W^{1,1}_{\mathrm loc}$, its distributional gradient is the function $\nabla I_2(x)=c(2-N)x|x|^{-N}$ (there is no contribution from the origin: the flux through $\partial B_r$ is $O(r^{N-1}r^{2-N})\to0$), and $|\nabla I_2|=c'|x|^{1-N}$. Let $h=U^p$. By (F1), (F2), $0<h\le Cw_{p\beta}$ with $p\beta>N+2$, so $h\in L^1\cap L^\infty$ and
\[
\nabla A=\nabla I_2*h .
\]
Split $\nabla I_2=k_1+k_2$ with $k_1=\nabla I_2\mathbf 1_{B_1}\in L^1$ and $k_2=\nabla I_2\mathbf 1_{B_1^c}\in L^\infty$. Then $k_1*h$ is uniformly continuous ($L^1*L^\infty$) and so is $k_2*h$ ($L^\infty*L^1$); hence $A\in C^1$ and $\nabla A$ is bounded.
For the decay, let $s>N$; we claim
\begin{equation}\label{eq:riesz-decay}
\int|x-y|^{1-N}w_s(y)\dd y\le C(1+|x|)^{1-N}.
\end{equation}
For $|x|\le2$ the left side is bounded by $\int_{|x-y|\le1}|x-y|^{1-N}\dd y+\|w_s\|_{L^1}\le C$. For $|x|\ge2$ split into $\{|y-x|\le|x|/2\}$, where $|y|\ge|x|/2$ and so $w_s(y)\le C|x|^{-s}$, giving at most $C|x|^{-s}\int_{|z|\le|x|/2}|z|^{1-N}\dd z=C|x|^{1-s}\le C|x|^{1-N}$; and $\{|y-x|>|x|/2\}$, where $|x-y|^{1-N}\le C|x|^{1-N}$, giving at most $C|x|^{1-N}\|w_s\|_{L^1}$. With $s=p\beta>N$ we obtain $|\nabla A|\le C(1+|x|)^{1-N}$.

\emph{Step 3: $f\in H^1$.} Since $U\in C^1$, $U>0$ and $A\in C^1$, the function $f=AU^{p-1}$ belongs to $C^1(\R^N)$, with
\[
\nabla f=\nabla A\,U^{p-1}+(p-1)AU^{p-2}\nabla U=\nabla A\,U^{p-1}+(p-1)V\nabla U .
\]
The first term is bounded by $C\|U\|_\infty^{p-1}(1+|x|)^{1-N}\in L^2$, because $2(N-1)>N$ for $N\ge3$. The second is in $L^2$ because $V\in L^\infty$ and $\nabla U\in L^2$. A $C^1$ function that is in $L^2$ together with its classical gradient lies in $H^1$; thus $f\in H^1$.

\emph{Step 4.} From $\hat U=\hat f/(1+|\zeta|^2)$,
\[
(1+|\zeta|^2)^{3/2}|\hat U|=(1+|\zeta|^2)^{1/2}|\hat f|\in L^2,
\]
so $U\in H^3$. \qedhere
\end{proof}

\begin{remark}
The only place where $p<2$ could have hurt is the factor $U^{p-2}$ in $\nabla(U^{p-1})$, which is singular where $U$ is small; it is absorbed by $A$, via the bounded function $V=AU^{p-2}$. This is the same mechanism that makes the weighted estimates of Lemma~\ref{lem:key} work.
\end{remark}

\begin{proposition}[natural constraint]\label{prop:natural}
Assume (A) and $|\varepsilon|<\varepsilon_0$, with $\varepsilon_0$ small enough that $2Mc_1\varepsilon_0\,C_*\,N\max_{i,j}\|\partial_i\partial_jU\|_{H^1}<\kappa$, where $\|\psi\|_{H^1}\le C_*\|\psi\|_Y$. Then
\[
\nabla\Phi_\varepsilon(\xi)=-(\kappa I-B)\,c,\qquad B_{ij}=\langle\partial_i\partial_jU,\tilde w_\varepsilon(\xi)\rangle_{H^1},\quad\|B\|\le C|\varepsilon| ,
\]
and $\kappa I-B$ is invertible. Moreover
\begin{equation}\label{eq:gradPhi-pert}
\partial_{\xi_j}\Phi_\varepsilon(\xi)=\frac\varepsilon2\int(\partial_ja)_\xi\,v^2 .
\end{equation}
Consequently, if $\nabla\Phi_\varepsilon(\xi)=0$, then $c=0$, and $u=\tau_\xi\left(U+\tilde w_\varepsilon(\xi)\right)$ is a positive solution of \eqref{eq:pertp<2}.
\end{proposition}
\begin{proof}
\emph{(1) Regularity.} By Proposition~\ref{prop:H3}, $\partial_iU\in H^2$. Rewriting the definition of $g$,
\[
v=T_0(v)-\varepsilon G*(a_\xi v)+\sum_ic_i\partial_iU .
\]
For $v\in\mathcal W_Y$ we have $0<\mathsf n(v)\le Cw_\beta\in L^2$, so $T_0(v)=G*\mathsf n(v)\in H^2$; also $a_\xi v\in L^2$, so $G*(a_\xi v)\in H^2$. Hence $v\in H^2$.

\emph{(2) Derivative of $\Phi_\varepsilon$.} By the chain rule, with $\partial_w\Psi_\varepsilon(\xi,w)\eta=\langle v-T_0(v),\eta\rangle_{H^1}+\varepsilon\int a_\xi v\,\eta=\langle g,\eta\rangle_{H^1}$ (Lemma~\ref{lem:grad}, and $\langle G*h,\eta\rangle_{H^1}=\int h\eta$ for $h\in L^2$, $\eta\in H^1$),
\[
\partial_{\xi_j}\Phi_\varepsilon(\xi)=\partial_{\xi_j}\Psi_\varepsilon(\xi,\tilde w)+\langle g,\partial_{\xi_j}\tilde w\rangle_{H^1}=\frac\varepsilon2\int(\partial_ja)_\xi v^2 .
\]
The last term of the chain rule vanishes because $\partial_{\xi_j}\tilde w\in Y_\perp$ (the space $Y_\perp$ is closed and $\tilde w$ is $C^1$ into it) while $g\in\operatorname{span}\{\partial_iU\}$. This is \eqref{eq:gradPhi-pert}.

\emph{(3) Translation identity.} The functional $f_0$ is $C^1$ on $H^1$ by a standard application of the  Hardy--Littlewood--Sobolev inequality, since $\frac{N+2}N<p<\frac{N+2}{N-2}$ (see \cite{MVS2}), and $f_0'(v)\eta=\langle v,\eta\rangle_{H^1}-\int A_vv^{p-1}\eta=\langle v-T_0v,\eta\rangle_{H^1}$ for all $\eta\in H^1$, because $\mathsf n(v)\in L^2$. As $v\in H^2$, $s\mapsto\tau_{se_j}v$ is $C^1$ into $H^1$ with derivative $-\partial_jv$ at $s=0$, and $f_0(\tau_{se_j}v)\equiv f_0(v)$. Differentiating at $s=0$,
\[
\langle v-T_0v,\partial_jv\rangle_{H^1}=0 .
\]
Therefore
\[
\langle g,\partial_jv\rangle_{H^1}=\varepsilon\int a_\xi v\,\partial_jv=-\frac\varepsilon2\int(\partial_ja)_\xi v^2=-\partial_{\xi_j}\Phi_\varepsilon(\xi),
\]
by an integration by parts.

\emph{(4) Expansion of the left side.} Since $g=\sum_ic_i\partial_iU$ and $\partial_jv=\partial_jU+\partial_j\tilde w$,
\[
\langle g,\partial_jv\rangle_{H^1}=\sum_ic_i\left[\kappa\delta_{ij}+\langle\partial_iU,\partial_j\tilde w\rangle_{H^1}\right],\qquad
\langle\partial_iU,\partial_j\tilde w\rangle_{H^1}=-\langle\partial_i\partial_jU,\tilde w\rangle_{H^1}=-B_{ij}.
\]
The integration by parts needs $\partial_iU\in H^2$ and $\partial_i\partial_jU\in H^1$, i.e.\ $U\in H^3$ (Proposition~\ref{prop:H3}), and $\tilde w\in H^1$. Hence $\partial_{\xi_j}\Phi_\varepsilon=-\sum_i(\kappa\delta_{ij}-B_{ij})c_i$, and $B$ is symmetric. Finally $|B_{ij}|\le\|\partial_i\partial_jU\|_{H^1}\|\tilde w\|_{H^1}\le C_*\|\partial_i\partial_jU\|_{H^1}\,2Mc_1|\varepsilon|$, so $\|B\|_{\mathrm op}\le N\max|B_{ij}|<\kappa$ by the choice of $\varepsilon_0$, and $\kappa I-B$ is invertible. The last claims follow: $c=0$ gives $g=0$, i.e.\ $v=T_0(v)-\varepsilon G*(a_\xi v)$, and applying $(-\Delta+1)$ and translating gives \eqref{eq:pertp<2}; positivity is Proposition~\ref{prop:fixed}.
\end{proof}

\subsection{Existence of a critical point and conclusion}\label{sec:critpt}

\begin{theorem}[perturbed solutions for $\alpha=2$, $p<2$]\label{thm:main}
Let $N\ge3$, $1+\frac2N<p<\min \left\{2,\frac{N+2}{N-2} \right\}$, and let $a$ satisfy (A) with $a\ge0$, $a\not\equiv0$. There is $\varepsilon_1>0$ such that for every $0<|\varepsilon|<\varepsilon_1$ equation \eqref{eq:pertp<2} has a positive solution $u_\varepsilon\in H^1$ of the form
\[
u_\varepsilon=\tau_{\xi_\varepsilon}\left(U+\tilde w_\varepsilon(\xi_\varepsilon)\right),\qquad\left\|\tau_{-\xi_\varepsilon}u_\varepsilon-U\right\|_Y=\|\tilde w_\varepsilon(\xi_\varepsilon)\|_Y\le C|\varepsilon|,\qquad u_\varepsilon\ge\frac12U(\cdot-\xi_\varepsilon),
\]
where $\xi_\varepsilon$ is a maximum point of $\Phi_\varepsilon$ over a ball $\overline{B_R}$ if $\varepsilon>0$, a minimum point if $\varepsilon<0$, and every accumulation point of $\xi_\varepsilon$ as $\varepsilon\to0$ is a global maximum point of $\Gamma$. If $a\le0$, $a\not\equiv0$, the same holds with $\xi_\varepsilon$ a minimum point of $\Phi_\varepsilon$ over $\overline{B_R}$ if $\varepsilon>0$, a maximum point if $\varepsilon<0$, and accumulation points are global minimum points of $\Gamma$. 
\end{theorem}
\begin{proof}
Let first $a\ge0$. $\Gamma$ is continuous and, since $U>0$ and $a\ge0$, $a\not\equiv0$, we have $\Gamma>0$. Also $\Gamma(\xi)\to0$ as $|\xi|\to\infty$: for $\mu>0$ choose $r$ with $a\le\mu$ on $\{|x|\ge r\}$; then
\[
\int aU_\xi^2\le\mu\|U\|_{L^2}^2+\|a\|_\infty\int_{|x|<r}U_\xi^2=\mu\|U\|_2^2+\|a\|_\infty\int_{|y+\xi|<r}U^2(y)\dd y,
\]
and the last integral tends to $0$ as $|\xi|\to\infty$. Hence $\Gamma$ attains its maximum at some $\bar\xi$. Take $\mu_0:=\Gamma(\bar\xi)/2>0$ and $R>|\bar\xi|$ with $\sup_{|\xi|\ge R}\Gamma\le\Gamma(\bar\xi)-\mu_0$.

Let $C_0$ be the constant of Proposition~\ref{prop:expansion}, $c_0=\frac{p-1}{2p}\|U\|_{H^1}^2$, $s=\operatorname{sgn}\varepsilon$, and $\varepsilon_1:=\min\{\varepsilon_0,\mu_0/(3C_0)\}$. Then $s\Phi_\varepsilon=sc_0+|\varepsilon|\Gamma+sr_\varepsilon$, and for $0<|\varepsilon|<\varepsilon_1$
\[
\max_{|\xi|=R}\left(s\Phi_\varepsilon-sc_0\right)\le|\varepsilon|\left(\Gamma(\bar\xi)-\mu_0\right)+C_0\varepsilon^2<|\varepsilon|\Gamma(\bar\xi)-C_0\varepsilon^2\le s\Phi_\varepsilon(\bar\xi)-sc_0 .
\]
So the maximum of the $C^1$ function $s\Phi_\varepsilon$ over $\overline{B_R}$ is attained at an interior point $\xi_\varepsilon$ (a maximum of $\Phi_\varepsilon$ if $\varepsilon>0$, a minimum if $\varepsilon<0$), which is a critical point of $\Phi_\varepsilon$. By Proposition~\ref{prop:natural}, $u_\varepsilon=\tau_{\xi_\varepsilon}(U+\tilde w_\varepsilon(\xi_\varepsilon))$ solves \eqref{eq:pertp<2}, and the estimates follow from Proposition~\ref{prop:fixed}. Finally $s\Phi_\varepsilon(\xi_\varepsilon)\ge s\Phi_\varepsilon(\bar\xi)$ gives $\Gamma(\xi_\varepsilon)\ge\Gamma(\bar\xi)-2C_0|\varepsilon|$, so any limit point $\xi_*\in\overline{B_R}$ satisfies $\Gamma(\xi_*)=\max\Gamma$.

If $a\le0$, the equation \eqref{eq:pertp<2} with data $(\varepsilon,a)$ coincides with the one with data $(-\varepsilon,-a)$, and $-a\ge0$ satisfies (A); the reduced functional $\Phi_\varepsilon$ is the same function, and the interaction functional is $-\Gamma$. Applying the first case with $\varepsilon'=-\varepsilon$ gives the stated maximum/minimum behaviour and that limit points are global minimum points of $\Gamma$. \qedhere
\end{proof}

As we have anticipated, stable critical points of the potential \(V\) give rise to solutions for \(\varepsilon\) small.
\begin{theorem}[general perturbation, nonzero local degree]\label{thm:general}
Assume (A) and let $\xi_0$ be an isolated critical point of $\Gamma$ with $\deg(\nabla\Gamma,B_r(\xi_0),0)\ne0$ for all small $r>0$ (for example a strict local extremum or a nondegenerate critical point). Then for every small $r>0$ there is $\varepsilon_r>0$ such that for $0<|\varepsilon|<\varepsilon_r$ equation \eqref{eq:pertp<2} has a positive solution
\[
u_\varepsilon=\tau_{\xi_\varepsilon}\left(U+\tilde w_\varepsilon(\xi_\varepsilon)\right),\qquad|\xi_\varepsilon-\xi_0|<r,\quad\|\tilde w_\varepsilon(\xi_\varepsilon)\|_Y\le C|\varepsilon|,\quad u_\varepsilon\ge\frac12U(\cdot-\xi_\varepsilon).
\]
\end{theorem}
\begin{proof}
By \eqref{eq:gradPhi-pert} and $\partial_j\Gamma(\xi)=\frac12\int(\partial_ja)_\xi U^2$ (differentiate $\Gamma(\xi)=\frac12\int a(y+\xi)U(y)^2\dd y$),
\[
\varepsilon^{-1}\partial_j\Phi_\varepsilon(\xi)-\partial_j\Gamma(\xi)=\frac12\int(\partial_ja)_\xi\left(2U\tilde w+\tilde w^2\right),
\]
so that, with $\tilde w=\tilde w_\varepsilon(\xi)$ and $\|\tilde w\|_{L^2}\le C\|\tilde w\|_Y\le C|\varepsilon|$,
\[
\sup_\xi\left|\varepsilon^{-1}\nabla\Phi_\varepsilon(\xi)-\nabla\Gamma(\xi)\right|\le\frac12\|\nabla a\|_\infty\|\tilde w\|_2\left(2\|U\|_2+\|\tilde w\|_2\right)\le C|\varepsilon| .
\]
Take $r$ with $\xi_0$ the only zero of $\nabla\Gamma$ in $\overline{B_r(\xi_0)}$ and $\sigma=\min_{\partial B_r(\xi_0)}|\nabla\Gamma|>0$. For $C|\varepsilon|<\sigma$ the vector fields $\varepsilon^{-1}\nabla\Phi_\varepsilon$ and $\nabla\Gamma$ are homotopic without zeros on $\partial B_r(\xi_0)$, so $\deg(\varepsilon^{-1}\nabla\Phi_\varepsilon,B_r(\xi_0),0)=\deg(\nabla\Gamma,B_r(\xi_0),0)\ne0$ ($\nabla\Phi_\varepsilon$ is continuous by Proposition~\ref{prop:C1}). Hence $\Phi_\varepsilon$ has a critical point $\xi_\varepsilon\in B_r(\xi_0)$, and Proposition~\ref{prop:natural} and Proposition~\ref{prop:fixed} give the conclusion. \qedhere
\end{proof}

\section{The semiclassical problem for \texorpdfstring{$\alpha=2$, $p<2$}{alpha = 2, p < 2} under small oscillation of \texorpdfstring{$V$}{V}}\label{sec:semiclassical}

In this section we extend the analysis already performed to the semiclassical problem associated to the Choquard equation.

In the sequel, $N\ge3$, $\alpha=2$, $1+\frac2N<p<\min \left\{2,\frac{N+2}{N-2} \right\}$, $\gamma:=\frac1{p-1}$, and we keep using the notation already introduced in the previous sections. We study
\begin{equation}\label{eq:SC}
  \begin{cases}
-\varepsilon^2\Delta u+V(x)u=\varepsilon^{-2}(I_2*u^p)\,u^{p-1},\\
 u>0,\ u\in H^1(\R^N).  
  \end{cases}
\end{equation}

This singularly perturbed model presents an obstruction with respect to the theory developed earlier. 
This is due to the fact  that the true solution \(v\) has a far-field tail $v^{2-p}\approx A/V(\xi+\varepsilon z)$, so that its ratio to the frozen ground state $U^{[V(\xi)]}$ is $\left(V(\xi)/V(\xi+\varepsilon z)\right)^{1/(2-p)}$, which may be far from $1$ if $V$ is allowed to oscillate wildly. If $V$ has \emph{small oscillation}, this ratio is uniformly close to $1$, and the whole problem becomes an instance of the linear perturbation problem of Section~\ref{sec:fixedpoint}, with $\varepsilon a_\xi$ replaced by the bounded function $b_{\xi,\varepsilon}(z)=V(\xi+\varepsilon z)-1$ (after normalizing $\lambda_0=1$). Smallness of $\|b\|_\infty$ is all that the reduction uses; it does not need $\varepsilon\to0$. The limit $\varepsilon\to0$ is used only to identify the reduced function.

\subsection{Scaling, hypotheses, and the rescaled problem}

\paragraph{Hypotheses.} Let $\lambda_0>0$ be a reference value (we will choose $\lambda_0=V(\xi_0)$ or $\lambda_0 = \min_{\bar\Omega}V$).
\begin{itemize}
\item[\textbf{(V1)}] $V\in C^1(\R^N)$, with $\nabla V$ bounded and uniformly continuous.
\item[\textbf{(V2)}] $\sup_{x \in\mathbb{R}^N}\left|V(x)/\lambda_0-1\right|\le\eta_0$, with $\eta_0<\eta_*$, where $\eta_*=\eta_*(N,p)>0$ is the constant fixed in \eqref{eq:eta} below.
\end{itemize}

As usual, we first prove the main scaling property of our problem.
\begin{lemma}[scaling]\label{lem:scaling}
Let $\lambda_0>0$. If $u$ solves \eqref{eq:SC} with potential $V$, then $\tilde u(x')=\lambda_0^{-\gamma}u(x'/\sqrt{\lambda_0})$ solves \eqref{eq:SC} with the same $\varepsilon$ and potential $\tilde V(x')=V(x'/\sqrt{\lambda_0})/\lambda_0$. In particular $\tilde V(\sqrt{\lambda_0}\,\xi_0)=1$ if $V(\xi_0)=\lambda_0$, and the frozen problem $-\Delta v+\lambda v=(I_2*v^p)v^{p-1}$ has the solution $U^{[\lambda]}(z)=\lambda^\gamma U(\sqrt\lambda z)$.
\end{lemma}
\begin{proof}
Let $u(x)=\lambda^\gamma\tilde u(\sqrt\lambda x)$. Then $-\varepsilon^2\Delta u=\lambda^{\gamma+1}\varepsilon^2(-\Delta\tilde u)(\sqrt\lambda x)$, and $Vu=\lambda^\gamma V\tilde u$. In the Riesz term, the substitution $y'=\sqrt\lambda y$ gives $(I_2*u^p)(x)=\lambda^{\gamma p-1}(I_2*\tilde u^p)(\sqrt\lambda x)$, and $u^{p-1}=\lambda\,\tilde u^{p-1}$ since $\gamma(p-1)=1$. So the right side is $\varepsilon^{-2}\lambda^{\gamma p}(I_2*\tilde u^p)\tilde u^{p-1}$ and $\gamma p=\gamma+1$. Dividing by $\lambda^{\gamma+1}$ gives the claim. The last statement is the case $V\equiv\lambda$, $\varepsilon=1$ applied to $U$.
\end{proof}

By Lemma~\ref{lem:scaling} (which preserves (V1) and turns (V2) into $|\tilde V-1|\le\eta_0$), \emph{we assume from now on $\lambda_0=1$} whenever convenient.

\paragraph{The rescaled problem.} For $\xi\in\R^N$, $\varepsilon>0$, and $u(x)=v((x-\xi)/\varepsilon)$, equation \eqref{eq:SC} is equivalent to
\begin{equation}\label{eq:SCxi}
-\Delta v+\left(1+b_{\xi,\varepsilon}(z)\right)v=(I_2*v^p)v^{p-1},\qquad b_{\xi,\varepsilon}(z):=V(\xi+\varepsilon z)-1,\qquad\|b_{\xi,\varepsilon}\|_\infty\le\eta_0 ,
\end{equation}
and the energy 
\[
I_\varepsilon(u)=\frac12\int(\varepsilon^2|\nabla u|^2+Vu^2)-\frac1{2p}\varepsilon^{-2}\int(I_2*u^p)u^p
\]
equals $\varepsilon^N\mathcal J(b_{\xi,\varepsilon},v)$, where
\[
\mathcal J(b,v):=f_0(v)+\frac12\int b\,v^2 .
\]
In the notation of Section~\ref{sec:fixedpoint}, \eqref{eq:SCxi} is equivalent to 
\[
v=T_0(v)-G*(b\,v).
\]
Choose first $r_*>0$ so small that, for the map $\mathcal F_U(\zeta)_i:=\langle U(\cdot+\zeta)-U,\partial_iU\rangle_{H^1}$ (a $C^1$ map from $\R^N$ to itself, see Lemma~\ref{lem:orbitsep}),
\[
\sup_{|\zeta|\le r_*} \left\|D\mathcal F_U(\zeta)-\kappa I_N \right\|_{\mathrm op}<\frac{\kappa}{4}.
\]
This is possible since $\mathcal{F}_U$ is a $C^1$ map $\R^N\to\R^N$, see Lemma~\ref{lem:orbitsep}, $D\mathcal F_U(0)=\kappa I_N$ and $D\mathcal F_U$ is continuous. Recalling that $U$ is strictly radially decreasing, $d_*:=U(0)-U(r_*)>0$. Choose then $\rs\in(0,\delta_0)$ so small that
\[
MC_Q\rs\le\frac18 \hbox{ and } \rs<d_*.
\]
Set $\cs:=C_G(1+\|U\|_X+\rs)$, let $C_*$ satisfy $\|\psi\|_{H^1}\le C_*\|\psi\|_Y$, and put $m_3:=\max_{i,j}\|\partial_i\partial_jU\|_{H^1}<\infty$ (by Proposition~\ref{prop:H3}). Finally, fix $\eta_*>0$ such that
\begin{equation}\label{eq:eta}
M\cs\eta_*\le\min\left\{\frac18,\frac{\rs}4\right\},\qquad M\cs\eta_*(2MC_Q+1)\le\frac12,\qquad 2M\cs\eta_*\,C_*\,N\,m_3<\frac{\kappa}{4} .
\end{equation}
These constants depend only on $N$ and $p$ through $U$.

\subsection{Finite-dimensional reduction for an arbitrary small bounded potential}

As anticipated above, the finite-dimensional reduction can be performed for an arbitrary small bounded potential $b$. So,
fix $b\in L^\infty(\R^N)$ and put
\begin{align*}
\mathfrak S(b,w)&:=\Pi^\perp\left[w-\left(T_0(U+w)-T_0(U)\right)+G*\left(b(U+w)\right)\right],\\
\mathcal N_b(w)&:=L_\perp^{-1}\Pi^\perp\left[Q(w)-G*\left(b(U+w)\right)\right].
\end{align*}
As in \eqref{eq:LS}--\eqref{eq:fixed}, $\mathfrak S(b,w)=0$ if and only if $w=\mathcal N_b(w)$. 

\begin{proposition}\label{prop:redb}
Let $\|b\|_\infty<\eta_*$. Then:
\begin{enumerate}[label=(\roman*),nosep]
\item $\mathcal N_b$ has a unique fixed point $\tilde w(b)$ in $B_{\rs}=\{w\in Y_\perp:\|w\|_Y\le\rs\}$, and $\|\tilde w(b)\|_Y\le2M\cs\|b\|_\infty\le\rs/2$. Moreover $U+\tilde w(b)\ge\frac12U$.
\item The map $b\mapsto\tilde w(b)$ is of class $C^1$ from $\{b\in L^\infty:\|b\|_\infty<\eta_*\}$ into $Y_\perp$.
\item If $b\equiv\lambda-1$ is constant, then $\tilde w(b)=U^{[\lambda]}-U$; in particular $\|U^{[\lambda]}-U\|_Y\le2M\cs|\lambda-1|$ for $|\lambda-1|<\eta_*$.
\end{enumerate}
\end{proposition}
\begin{proof}
(i) This is the proof of Proposition~\ref{prop:fixed} with $\varepsilon a_\xi$ replaced by $b$. For $w\in B_{\rs}$, Lemma~\ref{lem:Q} gives $\|Q(w)\|_Y\le C_Q\rs^2$ and $\|G*(b(U+w))\|_Y\le C_G\|b\|_\infty(\|U\|_X+\|w\|_X)\le\cs\|b\|_\infty$. Hence $\|\mathcal N_b(w)\|_Y\le M[C_Q\rs^2+\cs\|b\|_\infty]\le\frac{\rs}8+\frac{\rs}4<\rs$. For $w,w'\in B_{\rs}$,
\[
\|\mathcal N_b(w)-\mathcal N_b(w')\|_Y\le M\left[2C_Q\rs+C_G\|b\|_\infty\right]\|w-w'\|_Y\le\left(\frac14+\frac18\right)\|w-w'\|_Y .
\]
Banach's theorem gives existence and uniqueness. For the size, $\|\tilde w\|_Y\le M[C_Q\|\tilde w\|_Y^2+\cs\|b\|_\infty]\le\frac18\|\tilde w\|_Y+M\cs\|b\|_\infty$, so $\|\tilde w\|_Y\le\frac87M\cs\|b\|_\infty\le2M\cs\|b\|_\infty\le\frac{\rs}2$. Finally $|\tilde w|\le\|\tilde w\|_Xw_\beta\le(\delta_0/c_U)U\le\frac12U$.

(ii) $\mathfrak S$ is $C^1$ on $\{\|b\|_\infty<\eta_*\}\times\{w\in Y_\perp:\|w\|_Y<\rs\}$: it is affine in $b$, with $\partial_b\mathfrak S[\hat b]=\Pi^\perp G*(\hat b(U+w))$, continuous in $w$ in $L(L^\infty,Y_\perp)$ by Lemma~\ref{lem:Q}(ii); and $C^1$ in $w$ by Proposition~\ref{prop:T0}. At a solution, as in Proposition~\ref{prop:C1},
\[
D_w\mathfrak S=L_\perp\left[I-L_\perp^{-1}\Pi^\perp\left(Q'(\tilde w)-G*(b\,\cdot)\right)\right],
\]
and the bracket differs from $I$ by an operator of norm $\le M\left(C_Q\cdot2M\cs\|b\|_\infty+C_G\|b\|_\infty\right)\le M\cs\eta_*(2MC_Q+1)\le\frac12$. So $D_w\mathfrak S$ is invertible, and the implicit function theorem gives a $C^1$ branch through each solution. Since $\|\tilde w(b)\|_Y\le\rs/2<\rs$, the branch stays in $B_{\rs}$ and coincides with $\tilde w(b)$ by uniqueness.

(iii) Let $b\equiv\lambda-1$ and $v=U+\tilde w(b)$, $g=v-T_0v+(\lambda-1)G*v=\sum c_i\partial_iU$. As in the proof of Proposition~\ref{prop:natural}, $v\in H^2$ and $\langle v-T_0v,\partial_jv\rangle_{H^1}=0$ by translation invariance, while $\langle(\lambda-1)G*v,\partial_jv\rangle_{H^1}=(\lambda-1)\int v\,\partial_jv=0$. So $\langle g,\partial_jv\rangle_{H^1}=0$, and the computation of step (4) of that proof gives $(\kappa I-B)c=0$ with $\|B\|<\kappa$ by \eqref{eq:eta}; thus $c=0$ and $v$ solves $-\Delta v+\lambda v=(I_2*v^p)v^{p-1}$. It is positive ($v\ge U/2$) and in $H^1$, hence by uniqueness (Theorem~1.2 of \cite{XZZ2}) and Lemma~\ref{lem:scaling}, $v=U^{[\lambda]}(\cdot-\zeta)$ for some $\zeta$. Next, $v$ is radial: for $\mathcal R\in O(N)$, $T_0$, $Q$, $L_\perp$, $\Pi^\perp$ and $b$ commute with $w\mapsto w\circ \mathcal R^{-1}$ (because $U$ is radial, $\operatorname{span}\{\partial_iU\}$ is $\mathcal R$-invariant, and the $H^1$ product is $\mathcal R$-invariant), so $w\circ\mathcal R^{-1}$ is another fixed point in $B_{\rs}$, and $w\circ \mathcal R^{-1}=w$ by uniqueness. A radial translate of the radial strictly decreasing function $U^{[\lambda]}$ has $\zeta=0$. \qedhere
\end{proof}

\subsection{The reduced function and the natural constraint}

Assume (V1), (V2) with $\lambda_0=1$ and $\eta_0<\eta_*$. For $\xi\in\R^N$, $\varepsilon>0$:
\[
b_{\xi,\varepsilon}=V(\xi+\varepsilon\,\cdot)-1,\qquad \tilde w_{\xi,\varepsilon}:=\tilde w(b_{\xi,\varepsilon}),\qquad v_{\xi,\varepsilon}:=U+\tilde w_{\xi,\varepsilon},\qquad u_{\xi,\varepsilon}(x):=v_{\xi,\varepsilon}\left(\frac{x-\xi}\varepsilon\right),
\]
\[
\Phi_\varepsilon(\xi):=\mathcal J(b_{\xi,\varepsilon},v_{\xi,\varepsilon})=\varepsilon^{-N}I_\varepsilon(u_{\xi,\varepsilon}).
\]
\begin{lemma}\label{lem:bxi}
Under assumption (V1), $\xi\mapsto b_{\xi,\varepsilon}\in L^\infty$ is $C^1$ with $\partial_{\xi_j}b_{\xi,\varepsilon}=(\partial_jV)(\xi+\varepsilon\,\cdot)$, and $\partial_{z_j}b_{\xi,\varepsilon}=\varepsilon\,\partial_{\xi_j}b_{\xi,\varepsilon}$.
\end{lemma}
\begin{proof}
For any $h\in\R^N$,
\begin{equation*}
		b_{\xi+h,\varepsilon}-b_{\xi,\varepsilon}-h\cdot(\nabla V)(\xi+\varepsilon\,\cdot)=\int_0^1h\cdot\left[(\nabla V)(\xi+th+\varepsilon z)-(\nabla V)(\xi+\varepsilon z)\right]\dd t,
\end{equation*}
 with sup norm $o(|h|)$ uniformly in $z$ by uniform continuity of $\nabla V$. Continuity of the derivative in $\xi$ follows from the same property. The last identity is the chain rule.
\end{proof}

\begin{proposition}\label{prop:natSC}
Assume (V1), (V2). Then $\Phi_\varepsilon\in C^1(\R^N)$ and
\begin{equation}\label{eq:gradPhi}
\nabla\Phi_\varepsilon(\xi)=\frac12\int(\nabla V)(\xi+\varepsilon z)\,v_{\xi,\varepsilon}(z)^2\dd z .
\end{equation}
Moreover, writing $g=v-T_0v+G*(b\,v)$ with $v=v_{\xi,\varepsilon}$, $b=b_{\xi,\varepsilon}$, we have $g=\sum_ic_i\partial_iU$ with
\[
(\kappa I-B)\,c=-\varepsilon\,\nabla\Phi_\varepsilon(\xi),\qquad B_{ij}=\langle\partial_i\partial_jU,\tilde w_{\xi,\varepsilon}\rangle_{H^1},\quad\|B\|_{\mathrm op}<\kappa .
\]
Consequently $|c|\le C\varepsilon$, and $\nabla\Phi_\varepsilon(\xi)=0$ if and only if $c=0$, if and only if $u_{\xi,\varepsilon}$ is a positive solution of \eqref{eq:SC}.
\end{proposition}
\begin{proof}
$\Phi_\varepsilon(\xi)=\mathcal J(b_{\xi,\varepsilon},U+\tilde w(b_{\xi,\varepsilon}))$ is a composition of $C^1$ maps: $\mathcal J$ is $C^1$ on $L^\infty\times\mathcal W_Y$ by Proposition~\ref{prop:f0}, $b\mapsto\tilde w(b)$ by Proposition~\ref{prop:redb}(ii), and $\xi\mapsto b_{\xi,\varepsilon}$ by Lemma~\ref{lem:bxi}. By the chain rule, with $\partial_w\mathcal J(b,v)[\eta]=\langle v-T_0v,\eta\rangle_{H^1}+\int b\,v\eta=\langle g,\eta\rangle_{H^1}$ (Lemma~\ref{lem:grad}),
\[
\partial_{\xi_j}\Phi_\varepsilon=\frac12\int\partial_{\xi_j}b\;v^2+\langle g,\partial_{\xi_j}\tilde w\rangle_{H^1}=\frac12\int(\partial_jV)(\xi+\varepsilon z)v^2\dd z ,
\]
since $\partial_{\xi_j}\tilde w\in Y_\perp\perp g$. This is \eqref{eq:gradPhi}. As in Proposition~\ref{prop:natural}: $v\in H^2$ (since $bv\in L^2$, $\partial_iU\in H^2$ by Proposition~\ref{prop:H3}), and $\langle v-T_0v,\partial_jv\rangle_{H^1}=0$ by translation invariance of $f_0$. By Lemma~\ref{lem:bxi} and an integration by parts,
\[
\langle g,\partial_jv\rangle_{H^1}=\int b\,v\,\partial_jv=-\frac12\int\partial_{z_j}b\;v^2=-\frac\varepsilon2\int(\partial_jV)(\xi+\varepsilon z)v^2=-\varepsilon\,\partial_{\xi_j}\Phi_\varepsilon .
\]
Expanding the left side exactly as in step (4) of Proposition~\ref{prop:natural} gives $\langle g,\partial_jv\rangle_{H^1}=\sum_i(\kappa\delta_{ij}-B_{ij})c_i$, and $\|B\|_{\mathrm op}\le N\,m_3C_*\|\tilde w\|_Y\le2M\cs\eta_*C_*Nm_3<\kappa$ by \eqref{eq:eta}. The bound $|c|\le C\varepsilon$ follows from $|\nabla\Phi_\varepsilon|\le\frac12\|\nabla V\|_\infty\|v\|_{L^2}^2\le C$. Finally $g=0$ gives $v=T_0(v)-G*(bv)$, i.e.\ \eqref{eq:SCxi}, hence \eqref{eq:SC} for $u_{\xi,\varepsilon}$; positivity: $v\ge U/2$.
\end{proof}

\subsection{The limit \texorpdfstring{$\varepsilon\to0$}{epsilon -> 0}}

Write $\lambda=V(\xi)$ for the frozen value and $h(\xi):=\frac12\|U\|_{L^2}^2\,V(\xi)^{2\gamma-N/2}\,\nabla V(\xi)$. Since $\|U^{[\lambda]}\|_{L^2}^2=\lambda^{2\gamma-N/2}\|U\|_{L^2}^2$, we have $h(\xi)=\frac12\|U^{[V(\xi)]}\|_{L^2}^2\,\nabla V(\xi)$.

\begin{lemma}\label{lem:limitSC}
Assume (V1), (V2). Let $\xi_n\to\bar\xi$ and $\varepsilon_n\to0$. Then $v_{\xi_n,\varepsilon_n}\to U^{[V(\bar\xi)]}$ locally uniformly and in $L^2(\R^N)$, and
\[
\nabla\Phi_{\varepsilon_n}(\xi_n)\longrightarrow h(\bar\xi).
\]
\end{lemma}
\begin{proof}
Write $v_n=v_{\xi_n,\varepsilon_n}=U+w_n$, $b_n=b_{\xi_n,\varepsilon_n}$, $g_n=v_n-T_0v_n+G*(b_nv_n)=\sum_ic_{n,i}\partial_iU$. We have $\|w_n\|_Y\le\rs/2$, hence $0<v_n\le Cw_\beta$, and $|c_n|\le C\varepsilon_n\to0$ by Proposition~\ref{prop:natSC}.

\emph{Compactness.} $v_n=T_0(v_n)-G*(b_nv_n)+\sum_ic_{n,i}\partial_iU$. By Lemma~\ref{lem:reg}, $T_0(v_n)$ and $G*(b_nv_n)$ are bounded in $C^{1,\gamma'}$ (right-hand sides bounded by $Cw_\beta$), and $\partial_iU\in C^{0,\gamma'}$ since $U=G*(VU)\in C^{1,\gamma'}$. Thus $(v_n)$ is bounded in $C^{0,\gamma'}(\R^N)$, and by Arzel\`a--Ascoli a subsequence converges locally uniformly to some $v^*$. Also $(w_n)$ is bounded in $Y\subset H^1$, so $w_n\rightharpoonup w^*:=v^*-U$ weakly in $H^1$ (along a further subsequence), and $\|w^*\|_X\le\liminf\|w_n\|_X$, $\|\nabla w^*\|_{L^2}\le\liminf\|\nabla w_n\|_{L^2}$, so $\|w^*\|_Y\le\liminf\|w_n\|_Y\le\rs/2$; moreover $w^*\in Y_\perp$ by weak convergence.

\emph{Limit equation.} By dominated convergence (all integrands are dominated by multiples of $w_\beta$ or $w_{p\beta}$, with $p\beta>N+2$), $A_{v_n}=I_2*v_n^p\to A_{v^*}$ pointwise, hence $T_0(v_n)\to T_0(v^*)$ pointwise; and $G*(b_nv_n)\to(\lambda-1)G*v^*$ pointwise because $b_n(z)=V(\xi_n+\varepsilon_nz)-1\to V(\bar\xi)-1=:\lambda-1$ pointwise and boundedly. Since $c_n\to0$, we obtain
\[
v^*=T_0(v^*)-(\lambda-1)G*v^* ,
\]
that is, $\mathfrak S(b_c,w^*)=0$ for the constant $b_c\equiv\lambda-1$, with $w^*\in B_{\rs}$. By the uniqueness in Proposition~\ref{prop:redb}(i) and (iii), $w^*=U^{[\lambda]}-U$, i.e. $v^*=U^{[V(\bar\xi)]}$. The limit does not depend on the subsequence, so the whole sequence converges locally uniformly; since $|v_n|\le Cw_\beta\in L^2$, dominated convergence gives convergence in $L^2$.

\emph{The gradient.} By \eqref{eq:gradPhi}, $\nabla\Phi_{\varepsilon_n}(\xi_n)=\frac12\int(\nabla V)(\xi_n+\varepsilon_nz)v_n^2\dd z$. The integrand converges pointwise to $\nabla V(\bar\xi)(U^{[\lambda]})^2$ and is dominated by $\|\nabla V\|_\infty Cw_{2\beta}\in L^1$. Hence the limit is $\frac12\nabla V(\bar\xi)\|U^{[\lambda]}\|_{L^2}^2=h(\bar\xi)$.
\end{proof}

\begin{corollary}\label{cor:C1conv}
For every $R>0$, $\sup_{|\xi|\le R}|\nabla\Phi_\varepsilon(\xi)-h(\xi)|\to0$ as $\varepsilon\to0$. Moreover, with $\Upsilon(s):=\frac12\|U\|_{L^2}^2\int_1^ss'^{\,2\gamma-N/2}\dd s'$ (so that $\Upsilon'>0$ and $\nabla(\Upsilon\circ V)=h$),
\[
\sup_{|\xi|\le R}\left|\Phi_\varepsilon(\xi)-\Upsilon(V(\xi))-\left(\Phi_\varepsilon(0)-\Upsilon(V(0))\right)\right|\to0 .
\]
\end{corollary}
\begin{proof}
The first claim follows from Lemma~\ref{lem:limitSC} by contradiction: if it failed, there would be $\varepsilon_n\to0$, $|\xi_n|\le R$ with $|\nabla\Phi_{\varepsilon_n}(\xi_n)-h(\xi_n)|\ge\sigma>0$; extracting $\xi_n\to\bar\xi$ and using continuity of $h$ contradicts the lemma. For the second, $\psi_\varepsilon=\Phi_\varepsilon-\Upsilon\circ V$ is $C^1$ with $\sup_{B_R}|\nabla\psi_\varepsilon|\to0$, hence $|\psi_\varepsilon(\xi)-\psi_\varepsilon(0)|\le R\sup_{B_R}|\nabla\psi_\varepsilon|\to0$ on the convex set $B_R$.
\end{proof}

\begin{remark}
The corollary says that, up to an additive constant depending on $\varepsilon$, $\Phi_\varepsilon$ is uniformly $C^1_{\mathrm loc}$-close to $\Upsilon\circ V$, with $\Upsilon'>0$. Thus the reduced function reproduces the critical-point structure of $V$ itself, as in the usual semiclassical theory, although the correction $\tilde w_{\xi,\varepsilon}$ is only $O(\eta_0)$ in $Y$, not $O(\varepsilon)$. Explicitly, $\Upsilon(s)=\frac{\|U\|_2^2}{2\theta}(s^\theta-1)$ with $\theta=2\gamma+1-\frac N2$, so $\Upsilon\circ V=c\,V^\theta+\text{const}$, which is the expected leading term of the reduced function. Here $\theta>0$: indeed $\theta>0\iff\frac2{p-1}>\frac{N-2}2\iff p<\frac{N+2}{N-2}$, which holds in the whole range considered.
\end{remark}

\subsection{Solutions concentrating at stable critical points and multiplicity}

\begin{theorem}[concentration at stable critical points of $V$]\label{thm:SCstable}
Let $V$ satisfy (V1) and let $\xi_0$ be an isolated critical point of $V$ with $\deg\left(\nabla V,B_r(\xi_0),0\right)\ne0$ for small $r$ (for example a nondegenerate critical point, or a strict local extremum). Assume (V2) with some $\lambda_0>0$:
\[
\sup_x\left|\frac{V(x)}{\lambda_0}-1\right|\le\eta_0<\eta_*(N,p).
\]
Then for every small $r>0$ there is $\varepsilon_r>0$ such that for all $\varepsilon\in(0,\varepsilon_r)$ equation \eqref{eq:SC} has a positive solution $u_\varepsilon$ of the form
\[
u_\varepsilon(x)=\lambda_0^{\gamma}\,v_\varepsilon\left(\sqrt{\lambda_0}\,\frac{x-\xi_\varepsilon}{\varepsilon}\right),\qquad|\xi_\varepsilon-\xi_0|<r,\quad v_\varepsilon\ge\frac12U,\quad\|v_\varepsilon-U^{[V(\xi_\varepsilon)/\lambda_0]}\|_{L^2}\to0 \ \ (\varepsilon\to0).
\]
In particular $u_\varepsilon$ concentrates at $\xi_0$ (take $r\to0$ along a sequence of $\varepsilon\to0$), and $u_\varepsilon(\xi_\varepsilon+\varepsilon\,\cdot)-U^{[V(\xi_\varepsilon)]}\to0$ in $L^2(\R^N)$.
\end{theorem}
\begin{proof}
By Lemma~\ref{lem:scaling} we may assume $\lambda_0=1$ (the change of variables $x'=\sqrt{\lambda_0}x$ is linear, so it preserves (V1), the isolation of the critical point and the value of the degree; $\xi_0$ is replaced by $\sqrt{\lambda_0}\xi_0$). Choose $r$ with $\xi_0$ the only zero of $\nabla V$ in $\overline{B_r(\xi_0)}$. With $h=\omega(\xi)\nabla V(\xi)$, $\omega=\frac12\|U\|_2^2V^{2\gamma-N/2}>0$ continuous, the homotopy $((1-t)\omega+t)\nabla V$ has no zero on $\partial B_r(\xi_0)$, so $\deg(h,B_r(\xi_0),0)=\deg(\nabla V,B_r(\xi_0),0)\ne0$. Let $\sigma=\min_{\partial B_r(\xi_0)}|h|>0$. By Corollary~\ref{cor:C1conv}, for $\varepsilon$ small, $|\nabla\Phi_\varepsilon-h|<\sigma$ on $\partial B_r(\xi_0)$, and $\nabla\Phi_\varepsilon$ is continuous (Proposition~\ref{prop:natSC}); by homotopy invariance of the Brouwer degree, $\deg(\nabla\Phi_\varepsilon,B_r(\xi_0),0)=\deg(h,B_r(\xi_0),0)\ne0$, so there is $\xi_\varepsilon\in B_r(\xi_0)$ with $\nabla\Phi_\varepsilon(\xi_\varepsilon)=0$. By Proposition~\ref{prop:natSC}, $u_{\xi_\varepsilon,\varepsilon}$ solves \eqref{eq:SC}. The statement on $v_\varepsilon$ follows from Proposition~\ref{prop:redb}(i) and Lemma~\ref{lem:limitSC} (applied along any sequence $\varepsilon_n\to0$, $\xi_{\varepsilon_n}\to\bar\xi$, using $V(\xi_{\varepsilon_n})\to V(\bar\xi)$). \end{proof}

We will use the following topological multiplicity result.
\begin{lemma}[Lusternik--Schnirelman]\label{lem:LS}
Let $\Omega\subset\R^N$ be a bounded open set, $g\in C(\bar\Omega)$ with $m:=\min_{\bar\Omega}g<\min_{\partial\Omega}g$, $K:=\{g=m\}$, and $\delta>0$ with $K_\delta:=\{\operatorname{dist}(\cdot,K)\le\delta\}\subset\Omega$. Put $m_\delta:=\min\{g(x):x\in\bar\Omega,\ \operatorname{dist}(x,K)\ge\delta\}$, so that $m_\delta>m$, and let $0<\sigma<(m_\delta-m)/4$. If $\Phi\in C^1(\Omega)\cap C(\bar\Omega)$ and $\sup_{\bar\Omega}|\Phi-g-c|<\sigma$ for some constant $c$, then $\Phi$ has at least $\operatorname{cat}(K,K_\delta)$ critical points in $\Omega$, all contained in $\{g\le m+2\sigma\}$.
\end{lemma}
\begin{proof}
Here $\operatorname{cat}(K,K_\delta)$ is the Lusternik--Schnirelman category of $K$ relative to $K_\delta$. Note that $m_\delta>m$ because $\{x\in\bar\Omega:\operatorname{dist}(x,K)\ge\delta\}$ is a nonempty compact set (it contains $\partial\Omega$, since $K_\delta\subset\Omega$) on which $g>m$; and $m_\delta\le\min_{\partial\Omega}g$.

Let $\tilde\Phi=\Phi-c$ and $\tilde\Phi^t=\{x\in\bar\Omega:\tilde\Phi(x)\le t\}$. For $t\le m+2\sigma$ and $x\in\tilde\Phi^t$ we have $g(x)<\tilde\Phi(x)+\sigma\le m+3\sigma<m_\delta$, hence $\operatorname{dist}(x,K)<\delta$. Thus $\tilde\Phi^t$ is a compact subset of $K_\delta\subset\Omega$ (closed in $\bar\Omega$ by continuity, and it does not meet $\partial\Omega$). Moreover $K\subset\tilde\Phi^{m+\sigma}$ because $\tilde\Phi<g+\sigma=m+\sigma$ on $K$, and $\tilde\Phi^t=\emptyset$ for $t<m-\sigma$ because $\tilde\Phi>g-\sigma\ge m-\sigma$.

For $1\le k\le k_0:=\operatorname{cat}(K,K_\delta)$ set $c_k:=\inf\{t\le m+2\sigma:\ \operatorname{cat}(\tilde\Phi^t,K_\delta)\ge k\}$. By monotonicity of the category, $\operatorname{cat}(\tilde\Phi^{m+\sigma},K_\delta)\ge k_0$, so $m-\sigma\le c_1\le\dots\le c_{k_0}\le m+\sigma$. The Palais--Smale condition holds on the compact set $\tilde\Phi^{m+2\sigma}\subset\Omega$, and a (locally Lipschitz) pseudo-gradient flow of the $C^1$ function $\tilde\Phi$ decreases $\tilde\Phi$, hence leaves $\tilde\Phi^{m+2\sigma}$ invariant, stays in $K_\delta$, and is homotopic to the identity there. The standard deformation argument therefore applies: if $c_k$ were not a critical value, the flow would deform $\tilde\Phi^{c_k+\epsilon}$ into $\tilde\Phi^{c_k-\epsilon}$, contradicting the definition of $c_k$; and if $c_k=c_{k+1}=\dots=c_{k+j}=:\bar c$, the set of critical points at level $\bar c$ has category (relative to $K_\delta$) at least $j+1$, so it contains more than $j$ points, possibly infinitely many. In all cases $\Phi$ has at least $k_0$ critical points in $\tilde\Phi^{m+\sigma}\cap\Omega$, and at each of them $g\le\tilde\Phi+\sigma\le m+2\sigma$. \end{proof}

\begin{lemma}[injectivity of the reconstruction map]\label{lem:orbitsep}
Assume (V1)--(V2) with $\lambda_0=1$, with the constants chosen in \eqref{eq:eta}. For every $\varepsilon>0$ and all $\xi,\xi'\in\R^N$,
\[
 u_{\xi,\varepsilon}=u_{\xi',\varepsilon}
 \quad\Longrightarrow\quad \xi=\xi'.
\]
In particular, the reconstruction map $\xi\mapsto u_{\xi,\varepsilon}$ is injective on $\R^N$, and therefore on the set of critical points of $\Phi_\varepsilon$.
\end{lemma}
\begin{proof}
Write $v_\xi=v_{\xi,\varepsilon}=U+w_\xi$ and
$v_{\xi'}=v_{\xi',\varepsilon}=U+w_{\xi'}$, where
$w_\xi,w_{\xi'}\in Y_\perp$ and
\[
\|w_\xi\|_Y,\|w_{\xi'}\|_Y
 \le 2M\cs\|b_{\xi,\varepsilon}\|_\infty
 \le 2M\cs\eta_0<2M\cs\eta_*\le\frac{\rs}{2}.
\]
For every $\xi$ the profile $v_\xi$ belongs to $H^2(\R^N)$: indeed $v_\xi=T_0(v_\xi)-G*(b\,v_\xi)+\sum_ic_i\partial_iU$ with every term in $H^2$, as in the proof of Proposition~\ref{prop:natSC} (this does not use criticality of $\xi$).

Suppose that $u_{\xi,\varepsilon}=u_{\xi',\varepsilon}$ and set
$s=(\xi'-\xi)/\varepsilon$. By the definition of the rescaled
profiles,
\begin{equation}\label{eq:profile-translate}
 v_{\xi'}(z)=v_\xi(z+s)\qquad(z\in\mathbb R^N).
\end{equation}
We first show that $|s|<r_*$. Since
$\|w_\xi\|_X,\|w_{\xi'}\|_X\le\rs/2$ and
$w_\beta(x)=(1+|x|)^{-\beta}\le1$, evaluation of
\eqref{eq:profile-translate} at $z=0$ gives
\[
 U(0)-\frac{\rs}{2}
 \le v_{\xi'}(0)=v_\xi(s)
 \le U(s)+\frac{\rs}{2}.
\]
If $|s|\ge r_*$, radial monotonicity of $U$ yields
$U(s)\le U(r_*)$, and hence
$U(0)-U(r_*)\le\rs$, contradicting the choice
$\rs<d_*=U(0)-U(r_*)$. Thus $|s|<r_*$. 
Here we use that $v_\xi=T_0(v_\xi)-G*(b\,v_\xi)+\sum_ic_i\partial_iU$, where the first two terms are in $C^{1,\gamma}$ by Lemma~\ref{lem:reg} (their arguments are bounded by $Cw_\beta$) and $\partial_iU$ is continuous; hence $v_\xi$ is continuous, and the identity \eqref{eq:profile-translate} and the bound $|w_\xi|\le\|w_\xi\|_Xw_\beta$ hold pointwise everywhere.

For $v\in H^2(\mathbb R^N)$ define
\[
 \mathcal F_v(\zeta)_i
 :=\langle v(\cdot+\zeta)-U,\partial_iU\rangle_{H^1},
 \qquad i=1,\ldots,N.
\]
The map $\mathcal F_v:\mathbb R^N\to\mathbb R^N$ is $C^1$.
Orthogonality of $w_\xi$ and $w_{\xi'}$ to the translation modes,
together with \eqref{eq:profile-translate}, gives
$\mathcal F_{v_\xi}(0)=\mathcal F_{v_\xi}(s)=0$.
Differentiating translations in $H^1$ and integrating by parts,
\[
 [D\mathcal F_v(\zeta)]_{ij}
 =\langle\partial_jv(\cdot+\zeta),\partial_iU\rangle_{H^1}
 =-\langle v(\cdot+\zeta),\partial_{ij}U\rangle_{H^1}.
\]
Consequently,
\[
 \left|[D\mathcal F_{v_\xi}(\zeta)-D\mathcal F_U(\zeta)]_{ij}\right|
 \le \|w_\xi\|_{H^1}\,\|\partial_{ij}U\|_{H^1}
 \le C_*\|w_\xi\|_Y m_3.
\]
It follows, using the standard estimate of the operator norm by $N$
times the largest matrix entry, that
\[
\sup_{\zeta\in\mathbb R^N}
\|D\mathcal F_{v_\xi}(\zeta)-D\mathcal F_U(\zeta)\|_{\mathrm op}
\le 2M\cs\eta_* C_*N m_3<\frac{\kappa}{4}.
\]
By the choice of $r_*$, for $|\zeta|\le r_*$ we therefore have
\[
\|D\mathcal F_{v_\xi}(\zeta)-\kappa I_N\|_{\mathrm op}
<\frac{\kappa}{2}.
\]
For any $a,b\in B_{r_*}(0)$, integration along the segment from $b$
to $a$ gives
\[
\left(\mathcal F_{v_\xi}(a)-\mathcal F_{v_\xi}(b)\right)
\cdot(a-b)
\ge\frac{\kappa}{2}|a-b|^2.
\]
Thus $\mathcal F_{v_\xi}$ is injective on $B_{r_*}(0)$. Since
$0,s\in B_{r_*}(0)$ and both are zeros, $s=0$. Hence
$\xi=\xi'$.
\end{proof}

We are now ready to prove a multiplicity result for equation \eqref{eq:SC}.
\begin{proposition}[multiplicity of distinct positive solutions]
\label{prop:SCcat}
Let $V$ satisfy (V1), let $\Omega$ be a bounded open set with
$m_0:=\min_{\bar\Omega}V<\min_{\partial\Omega}V$, and put
$K:=\left\{\xi\in\bar\Omega:V(\xi)=m_0\right\}$. Assume (V2) with
$\lambda_0=m_0$. For every $\delta>0$ such that $K_\delta\subset\Omega$,
the equation \eqref{eq:SC} has at least
$\operatorname{cat}(K,K_\delta)$ pairwise distinct positive
solutions for all sufficiently small $\varepsilon>0$. Their
concentration centers are critical points of $\Phi_\varepsilon$
in $\Omega$.
\end{proposition}
\begin{proof}
Normalize $\lambda_0=m_0=1$ by Lemma~\ref{lem:scaling}, so that $\Upsilon\circ V$ has its minimum over $\bar\Omega$ exactly on $K$ ($\Upsilon'>0$) and $\min_{\partial\Omega}\Upsilon\circ V>\min_{\bar\Omega}\Upsilon\circ V$. Take $R$ with $\bar\Omega\subset B_R$. By Corollary~\ref{cor:C1conv}, $\sup_{\bar\Omega}|\Phi_\varepsilon-\Upsilon\circ V-c_\varepsilon|\to0$ with $c_\varepsilon=\Phi_\varepsilon(0)-\Upsilon(V(0))$. Lemma~\ref{lem:LS} with $g=\Upsilon\circ V$ (applied with any admissible $\sigma$, and $\varepsilon$ so small that the sup above is $<\sigma$) gives $\operatorname{cat}(K,K_\delta)$ critical points of $\Phi_\varepsilon$ in $\Omega$, all in $\{g\le\min g+2\sigma\}$; by Proposition~\ref{prop:natSC}, each gives a positive solution. Since $\sigma$ can be taken arbitrarily small as $\varepsilon\to0$, their centers satisfy $V(\xi)\to m_0$. By Lemma~\ref{lem:orbitsep}, distinct critical points give distinct reconstructed solutions. The Lusternik--Schnirelman lower bound is therefore a lower bound for the number of pairwise distinct positive solutions of \eqref{eq:SC}.
\end{proof}

\begin{remark}
Smallness of the oscillation is used only to guarantee $\|b_{\xi,\varepsilon}\|_{L^\infty}\le\eta_0<\eta_*$, so that the reduction of Section~\ref{sec:fixedpoint} applies with $b$ in place of $\varepsilon a_\xi$. It is a global condition on $V$, which seems hard to avoid with this method because of the algebraic tail of $U$; we do not know whether a local condition near $\xi_0$ would suffice for the equation of this section. For the equation with the potential on both sides (Section~\ref{sec:KV}), the global condition can be dispensed with. The threshold $\eta_*$ depends only on $N$ and $p$, through $M,C_Q,C_G,\cs,\rs$, and is not explicit.

Assumption (V1) enters through the $C^1$ regularity of $\Phi_\varepsilon$ and the exact formula \eqref{eq:gradPhi} for its gradient. For the minimization-type result of Proposition~\ref{prop:SCcat} it can probably be relaxed, but we have not tried. The $C^1(H^1)$ regularity of $f_0$ used in Section~\ref{sec:natural} is standard \cite{MVS2}, and the deformation argument in Lemma~\ref{lem:LS} is the classical one \cite{AMS,CL}.
%
\end{remark}

\section{The singularly perturbed problem with a potential on both sides}\label{sec:KV}

In this section we consider a particular model in which the potential $V$ multiplies the nonlocal term as well:
\begin{equation}\label{eq:KV}
  \begin{cases}
  -\varepsilon^2\Delta u+V(x)u=\varepsilon^{-2}\left(I_2*(Vu^p)\right)\,Vu^{p-1},\\
  u>0,\ u\in H^1(\R^N),
  \end{cases}
\end{equation}
with $N\ge3$, $\alpha=2$ and $1+\frac2N<p<\min\left\{2,\frac{N+2}{N-2}\right\}$ as in Sections~\ref{sec:psmall} and \ref{sec:semiclassical}. The equation is the Euler--Lagrange equation of
\[
u \mapsto \frac12\int\left(\varepsilon^2|\nabla u|^2+Vu^2\right)-\frac1{2p}\,\varepsilon^{-2}\int\left(I_2*(Vu^p)\right)Vu^p ,
\]
so the variational structure used in Section~\ref{sec:semiclassical} is preserved.

As we shall see, the presence of the \emph{same} potential function $V$ \emph{both} in the left-hand and in the right-hand side of \eqref{eq:KV} makes the finite-dimensional reduction more flexible. Roughly speaking, with $A=I_2*(Vu^p)$ the equation reads $-\varepsilon^2\Delta u=V\,(\varepsilon^{-2}Au^{p-1}-u)$. The potential multiplies only the Thomas--Fermi combination $\varepsilon^{-2}Au^{p-1}-u$, which vanishes to leading order in the far field, because there $u^{2-p}\approx \varepsilon^{-2}A$. For the frozen profile this combination is $-\Delta U/\lambda$ and decays two orders faster than $U$. As a consequence the large oscillation of $V$ at distance $\gtrsim1/\varepsilon$ from the bump is harmless: it only enters through terms that carry the extra decay (Lemma~\ref{lem:KVres}), and it only changes the far-field linearization from $-\Delta+(2-p)$ to $-\Delta+(2-p)V$, which is still coercive (Lemma~\ref{lem:KVprinc}). For a potential $K\ne V$ on the right the first of these facts fails, and one needs the oscillation of $K/V$ to be small.

Our assumption on the potential $V$ can be relaxed as follows:
\begin{itemize}
\item[\textbf{(W)}] $V\in C^1(\R^N)$, $0<m:=\inf_{\mathbb{R}^N} V\le\sup_{\mathbb{R}^N} V=:M<\infty$, and $\nabla V$ is bounded, by $L$ say, and uniformly continuous.
\end{itemize}
We emphasize that no condition on $M/m$ is imposed.


\subsection{Dilation and the class of rescaled potentials}

As usual, we investigate the scaling properties of our problem.
\begin{lemma}[dilation]\label{lem:dil}
Let $W$ be measurable with $0<\inf_{\mathbb{R}^N} W\le\sup_{\mathbb{R}^N} W<\infty$, and let $v$ solve
\begin{equation}\label{eq:KVresc}
-\Delta v+Wv=\left(I_2*(Wv^p)\right)Wv^{p-1}.
\end{equation}
For $\lambda>0$ the function $\hat v(y)=v(y/\sqrt\lambda)$ solves \eqref{eq:KVresc} with the potential $\hat W(y)=W(y/\sqrt\lambda)/\lambda$. In particular $U_\lambda(y):=U(\sqrt\lambda\,y)$ solves \eqref{eq:KVresc} with $W\equiv\lambda$.
\end{lemma}
\begin{proof}
Write $v(z)=\hat v(\sqrt\lambda z)$. Then $-\Delta v=\lambda(-\Delta\hat v)(\sqrt\lambda z)$ and $Wv=\lambda\hat W\hat v$. In the nonlocal term, the substitution $y=\sqrt\lambda z'$ and $I_2(z-y/\sqrt\lambda)=\lambda^{(N-2)/2}I_2(\sqrt\lambda z-y)$ give
\[
\left(I_2*(Wv^p)\right)(z)=\lambda\cdot\lambda^{\frac{N-2}2-\frac N2}\left(I_2*(\hat W\hat v^p)\right)(\sqrt\lambda z)=\left(I_2*(\hat W\hat v^p)\right)(\sqrt\lambda z).
\]
Hence the equation for $v$, evaluated at $z=y/\sqrt\lambda$, reads
\[
\lambda\left[-\Delta\hat v+\hat W\hat v\right]=\lambda\left(I_2*(\hat W\hat v^p)\right)\hat W\hat v^{p-1}.
\]
The last claim is the case $v=U$, $W\equiv1$, $\lambda$ arbitrary.
\end{proof}

Note that, unlike Lemma~\ref{lem:scaling}, the amplitude does not scale: the frozen ground state for the potential $\lambda$ is $U(\sqrt\lambda\,\cdot)$.

For $\xi\in\R^N$ and $\varepsilon>0$ the rescaling $v(z)=u(\xi+\varepsilon z)$ turns \eqref{eq:KV} into \eqref{eq:KVresc} with $W=W_{\xi,\varepsilon}=V(\xi+\varepsilon\,\cdot)$. We always write $W=1+b$.

\begin{definition}[The class $\mathfrak B$]
 For $\hat m\in(0,1]$, $\hat L>0$, $\eta\in[0,\frac12]$ and $\varepsilon>0$ let $\mathfrak B=\mathfrak B_{\hat m,\hat L}(\eta,\varepsilon)$ be the set of functions $b$ on $\R^N$ such that
\[
\hat m\le1+b\le\hat m^{-1},\qquad |b(0)|\le\eta,\qquad |b(y)-b(y')|\le\hat L\varepsilon|y-y'| .
\]
In particular $|b(y)|\le\eta+\hat L\varepsilon|y|$. 
\end{definition}

\begin{remark}
The class is star-shaped: $tb\in\mathfrak B$ for $t\in[0,1]$, since $1+tb$ is a convex combination of $1$ and $1+b$.
\end{remark}
Two properties of the class $\mathfrak B$ will be useful:
\begin{itemize}
\item[(a)] \emph{Dilation at the center} ($\eta=0$). For any $\xi$ put $\lambda=V(\xi)$. By Lemma~\ref{lem:dil}, $\hat b(y)=V(\xi+\varepsilon y/\sqrt\lambda)/\lambda-1$ belongs to $\mathfrak B_{\hat m,\hat L}(0,\varepsilon)$ with $\hat m=m/M$ and $\hat L=L\,m^{-3/2}$.
\item[(b)] \emph{Patches} ($\eta>0$). After one global dilation we may assume $V(\xi_0)=1$. If $|V-1|\le\eta$ on a set $\Omega_0$ and $\xi\in\Omega_0$, then $b_{\xi,\varepsilon}=V(\xi+\varepsilon\,\cdot)-1\in\mathfrak B_{\hat m,\hat L}(\eta,\varepsilon)$ with $\hat m=\min\{m,1/M\}$ and $\hat L=L$ (the function $z\mapsto V(\xi+\varepsilon z)$ is $L\varepsilon$-Lipschitz).
\end{itemize}
In both cases the \emph{global} oscillation of $V$ is only reflected in the constant $\hat m$.

We use $w_s=(1+|x|)^{-s}$, $X$, $Y$, $c_U$, $G$, $A=I_2*U^p$ from Section~\ref{sec:psmall}, and write $\varrho:=AU^{p-2}$ for the function called $V$ there, to avoid confusion with the potential. Recall that $0<U\le Cw_\beta$, $U\ge c_Uw_\beta$, $\varrho\in C(\R^N)$, $\varrho\to1$ at infinity, $\varrho U=AU^{p-1}$, $p\beta>N+2$, $(p-1)\beta+N-2=\beta$, and $-\Delta U+U=\varrho U$. In particular
\begin{equation}\label{eq:DeltaU}
\Delta U=U-AU^{p-1}=(1-\varrho)\,U\quad\text{a.e.}
\end{equation}
For $b\in L^\infty$ with $\inf_{\mathbb{R}^N} (1+b)>0$ put $W=1+b$,
\begin{align*}
&A_{b}(v):=I_2*(Wv^p), \\
&T_b(v):=G*\left[W\,A_b(v)\,v^{p-1}\right], \\
&F_b(w):=(U+w)-T_b(U+w)+G*\left(b(U+w)\right).
\end{align*}
Then $v=U+w$ solves \eqref{eq:KVresc} if and only if $F_b(w)=0$, and $T_0$ of Section~\ref{sec:psmall} is $T_b$ with $b=0$.

\subsection{The residual estimate}

\begin{lemma}\label{lem:KVres}
Let $b\in\mathfrak B_{\hat m,\hat L}(\eta,\varepsilon)$. Then $F_b(0)=G*\mathcal R_b$ with
\[
\mathcal R_b=b\,\Delta U-(1+b)\left(I_2*(bU^p)\right)U^{p-1},
\]
and
\[
\|F_b(0)\|_Y\le C_1\left(\eta+\chi(\varepsilon)+\varepsilon\right),\qquad \chi(\varepsilon):=\sup_{y}\min\{\hat L\varepsilon|y|,\hat m^{-1}\}\,|1-\varrho(y)|,
\]
where $C_1$ depends only on $N,p,\hat m,\hat L$. Moreover $\chi(\varepsilon)\to0$ as $\varepsilon\to0$. In particular, for $\eta=0$, $\|F_b(0)\|_Y\to0$ as $\varepsilon\to0$, uniformly in $b\in\mathfrak B_{\hat m,\hat L}(0,\varepsilon)$, although $\|b\|_\infty$ need not be small.
\end{lemma}
\begin{proof}
Since $U\in H^2$, $G*(-\Delta+1)U=U$, so $F_b(0)=G*\mathcal R_b$ with $\mathcal R_b=-\Delta U+(1+b)U-(1+b)A_b(U)U^{p-1}$. Write $A_b(U)=A+I_2*(bU^p)$ and use $-\Delta U=AU^{p-1}-U$:
\[
\mathcal R_b=bU-bAU^{p-1}-(1+b)\left(I_2*(bU^p)\right)U^{p-1}=b\,(U-AU^{p-1})-(1+b)\left(I_2*(bU^p)\right)U^{p-1},
\]
and $U-AU^{p-1}=\Delta U$. We now prove the estimates.

By Lemma~\ref{lem:reg} (with $s=\beta$), $\|G*h\|_Y\le C\,\|h\|_X$ whenever $h\in X$, so it suffices to bound $\|\mathcal R_b\|_X$.

First term. By \eqref{eq:DeltaU}, $|b\Delta U|=|b|\,|1-\varrho|\,U\le C|b|\,|1-\varrho|\,w_\beta$. Since $|b|\le\eta+\min\{\hat L\varepsilon|y|,\hat m^{-1}\}$,
\[
\|b\Delta U\|_X\le C\left(\eta\,\|1-\varrho\|_\infty+\chi(\varepsilon)\right).
\]
Second term. We have $|b|U^p\le\eta U^p+\hat L\varepsilon|y|U^p$ and $|y|U^p\le Cw_{p\beta-1}$ with $p\beta-1>N$. By the estimate used in Step~1 of the proof of Theorem~\ref{thm:fred} (or Lemma~\ref{lem:KVtail} below), $I_2*w_s\le Cw_{N-2}$ for $s>N$. Hence $|I_2*(bU^p)|\le\eta A+C\hat L\varepsilon\,w_{N-2}$, and, using $AU^{p-1}=\varrho U\le Cw_\beta$ and $w_{N-2}U^{p-1}\le Cw_{N-2+(p-1)\beta}=Cw_\beta$,
\[
\left\|(1+b)\left(I_2*(bU^p)\right)U^{p-1}\right\|_X\le C\hat m^{-1}(\eta+\hat L\varepsilon).
\]

Finally, fix $R>0$. For $|y|\le R$ the quantity is at most $\hat L\varepsilon R\|1-\varrho\|_\infty$; for $|y|>R$ it is at most $\hat m^{-1}\sup_{|y|>R}|1-\varrho|$, which tends to $0$ as $R\to\infty$ because $\varrho\to1$. So $\limsup_{\varepsilon\to0}\chi(\varepsilon)\le\hat m^{-1}\sup_{|y|>R}|1-\varrho|$ for every $R$, i.e.\ $\chi(\varepsilon)\to0$.
\end{proof}

The next technical estimates will be used later.
\begin{lemma}\label{lem:KVtail}
Let $s>N$. There is $C_s$ such that for all $z\in\R^N$ and $R\ge1$,
\[
\int_{\R^N}|z-y|^{2-N}w_s(y)\dd y\le C_s\,w_{N-2}(z),\qquad
\int_{|y|>R}|z-y|^{2-N}w_s(y)\dd y\le C_sR^{N-s}\,w_{N-2}(z).
\]
\end{lemma}
\begin{proof}
The first inequality is the splitting used in Step~1 of the proof of Theorem~\ref{thm:fred}; it also follows from the second with $R=1$ and the bound on $\{|y|\le1\}$. For the second, let $\Theta_R(z)$ denote the integral.

If $|z|\le R/2$, then $|z-y|\ge|y|/2$ on $\{|y|>R\}$, so $\Theta_R(z)\le C\int_{|y|>R}|y|^{2-N-s}\dd y=CR^{2-s}\le CR^{N-s}\,w_{N-2}(z)$, since $w_{N-2}(z)\ge c\,R^{2-N}$.

If $|z|>R/2$, split at $|y-z|=|z|/2$. On $\{|y-z|\le|z|/2\}$ we have $1+|y|\ge\frac12(1+|z|)$, so that part is at most $C(1+|z|)^{-s}\int_{|y-z|\le|z|/2}|z-y|^{2-N}\dd y\le C(1+|z|)^{2-s}\le C(R/2)^{N-s}w_{N-2}(z)$. On the complement $|z-y|^{2-N}\le C(1+|z|)^{2-N}$ and $\int_{|y|>R}w_s\le CR^{N-s}$.
\end{proof}

\subsection{Estimates for \texorpdfstring{$T_b$}{T\_b}}\label{sec:KVQ}

Throughout this subsection $\hat m\in(0,1]$ is fixed and
\[
\mathcal O=\mathcal O(\hat m):=\left\{b\in L^\infty(\R^N) \mid \frac12\hat m\le1+b\le2\hat m^{-1}\right\}.
\]
The class $\mathfrak B_{\hat m,\hat L}(\eta,\varepsilon)$ is contained in $\mathcal O$, and in the interior of $\mathcal O$ with respect to $\|\cdot\|_\infty$ (this is why $\hat m/2$ and $2/\hat m$ appear). We write $W=1+b$, $W'=1+b'$, and $\delta_0,\mathcal W_X,c_U$ are as in Section~\ref{sec:psmall} (they do not depend on $b$). The letter $C$ denotes constants depending only on $N,p,\hat m$ and on $U$ (through $c_U$, $\|\varrho\|_\infty$ and the constants in Lemmas~\ref{lem:reg} and \ref{lem:key}). For $u\in\mathcal W_X$ put
\[
A_{b,u}:=I_2*(Wu^p),\qquad \mathsf n_b(u):=W\,A_{b,u}\,u^{p-1},\qquad T_b(u)=G*\mathsf n_b(u),
\]
and, for $\psi\in X$,
\[
a_{b,u}:=A_{b,u}u^{p-2},\qquad J_{b,u}\psi:=I_2*\left(Wu^{p-1}\psi\right),\qquad
\mathcal L_{b,u}\psi:=(p-1)W\,a_{b,u}\,\psi+p\,W\,u^{p-1}\,J_{b,u}\psi .
\]
For $b=0$ these are $A_u$, $\mathsf n$, $T_0$, $V_u$ and $\mathcal L_u$ of Proposition~\ref{prop:T0}.

\begin{lemma}[estimates for $T_b$]\label{lem:KVQ}
There are $C>0$ depending only on $N,p,\hat m$ and $U$ such that, with $\|\cdot\|_X$ the norm of $X$ and $\mathcal W_X=\{U+\varphi:\|\varphi\|_X<\delta_0\}$:
\begin{enumerate}[label=(\roman*),nosep]
\item For $b,b'\in\mathcal O$ and $u,u'\in\mathcal W_X$: $0\le\mathsf n_b(u)\le Cw_\beta$, $T_b(u)\in Y$ with $\|T_b(u)\|_Y\le C$, and
\[
\|T_b(u)-T_{b'}(u')\|_Y\le C\left(\|u-u'\|_X+\|b-b'\|_\infty\right).
\]
\item For $b\in\mathcal O$, $u\in\mathcal W_X$ and $\psi\in X$ with $u+\psi\in\mathcal W_X$,
\[
\left\|T_b(u+\psi)-T_b(u)-G*\mathcal L_{b,u}\psi\right\|_Y\le C\|\psi\|_X^2 .
\]
Hence $T_b:\mathcal W_X\to Y$ is Fr\'echet differentiable (with $X$ as the domain norm) with $T_b'(u)=G*\mathcal L_{b,u}$.
\item For $b,b'\in\mathcal O$ and $u,u'\in\mathcal W_X$,
\[
\|T_b'(u)\|_{L(X,Y)}\le C,\qquad \|T_b'(u)-T_{b'}'(u')\|_{L(X,Y)}\le C\left(\|u-u'\|_X+\|b-b'\|_\infty\right).
\]
In particular $T_b\in C^{1,1}(\mathcal W_X;Y)$ for each $b\in\mathcal O$ (constants uniform in $b$), and $\|L_b-L_{b'}\|_{L(Y)}\le C_3\|b-b'\|_\infty$, where $L_b=I-T_b'(U)+G*(b\,\cdot)$ is the operator of Section~\ref{sec:KVlin}.
\item The map $(b,u)\mapsto T_b(u)$ is of class $C^1$ from $\operatorname{int}\mathcal O\times\mathcal W_X$ (with the norms $\|\cdot\|_\infty$ and $\|\cdot\|_X$) into $Y$. In particular $(b,w)\mapsto F_b(w)$ is $C^1$ on $\operatorname{int}\mathcal O\times\{w\in Y:\|w\|_X<\delta_0\}$.
\item The map $Q_b(w):=T_b(U+w)-T_b(U)-T_b'(U)w$ satisfies, for $w,w'\in Y$ with $\|w\|_X,\|w'\|_X<\delta_0$,
\[
\|Q_b(w)\|_Y\le C_Q'\|w\|_X\|w\|_Y,\qquad \|Q_b(w)-Q_b(w')\|_Y\le C_Q'\left(\|w\|_X+\|w'\|_X\right)\|w-w'\|_Y ,
\]
and $Q_b'(w)=T_b'(U+w)-T_b'(U)$ has norm $\|Q_b'(w)\|_{L(Y)}\le C_Q'\|w\|_X$.
\end{enumerate}
\end{lemma}

\begin{proof}
We use the pointwise facts of Lemma~\ref{lem:key} and its proof, which do not involve $b$: for $u,u'\in\mathcal W_X$, $\psi\in X$ with $u+\psi\in\mathcal W_X$,
\begin{enumerate}[label=(K\arabic*),nosep]
\item $\frac12U\le u\le\frac32U$;
\item $|u^s-u'^s|\le CU^s\|u-u'\|_X$ for $s\in\{p,p-1,p-2\}$, and $|(u+\psi)^q-u^q-qu^{q-1}\psi|\le CU^q\|\psi\|_X^2$ for $q\in\{p,p-1\}$;
\item $I_2*U^p=A$, $\ I_2*(U^{p-1}|\psi|)\le C\|\psi\|_XA$ (because $U^{p-1}w_\beta\le U^p/c_U$), $\ AU^{p-2}=\varrho\le\|\varrho\|_\infty$ and $AU^{p-1}=\varrho U\le Cw_\beta$.
\end{enumerate}
Since $\frac12\hat m\le W,W'\le2\hat m^{-1}$ and $|W-W'|=|b-b'|$, we deduce, for $\delta:=\|u-u'\|_X+\|b-b'\|_\infty$:
\begin{enumerate}[label=(\alph*),nosep]
\item $0<A_{b,u}\le CA$.
\item $|A_{b,u}-A_{b',u'}|\le C\delta A$. Indeed $A_{b,u}-A_{b',u'}=I_2*\left((W-W')u^p\right)+I_2*\left(W'(u^p-u'^p)\right)$, and we use (K1), (K2), (K3).
\item $\left|I_2*\left(W[(u+\psi)^p-u^p-pu^{p-1}\psi]\right)\right|\le C\|\psi\|_X^2A$, by (K2) with $q=p$ and (K3).
\item $I_2*(Wu^{p-1}|\psi|)\le C\|\psi\|_XA$, by (K1) and (K3); more generally, $|J_{b,u}\psi-J_{b',u'}\psi|\le C\delta\|\psi\|_XA$, since $|Wu^{p-1}-W'u'^{p-1}|\le|W-W'|u^{p-1}+W'|u^{p-1}-u'^{p-1}|\le C\delta U^{p-1}$.
\item $0<a_{b,u}\le C\varrho$ and $|a_{b,u}-a_{b',u'}|\le C\delta\varrho$. Indeed $u^{p-2}\le2^{2-p}U^{p-2}$, so $a_{b,u}\le CAU^{p-2}=C\varrho$; and
$a_{b,u}-a_{b',u'}=(A_{b,u}-A_{b',u'})u^{p-2}+A_{b',u'}(u^{p-2}-u'^{p-2})$ with (b), (K2) for $s=p-2$.
\end{enumerate}

\emph{(i).} By (a) and (K1), $0\le\mathsf n_b(u)=WA_{b,u}u^{p-1}\le C\,AU^{p-1}\le Cw_\beta$. Lemma~\ref{lem:reg} (with $s=\beta$) gives $T_b(u)\in Y$, $\|T_b(u)\|_Y\le C$. For the difference,
\[
\mathsf n_b(u)-\mathsf n_{b'}(u')=(W-W')A_{b,u}u^{p-1}+W'(A_{b,u}-A_{b',u'})u^{p-1}+W'A_{b',u'}(u^{p-1}-u'^{p-1}),
\]
and by (a), (b), (K2) with $s=p-1$ each term is bounded by $C\delta\,AU^{p-1}\le C\delta w_\beta$. Lemma~\ref{lem:reg} gives the estimate.

\emph{(ii).} Put $P_1=A_{b,u}$, $P_2=u^{p-1}$, $\Delta P_1=A_{b,u+\psi}-A_{b,u}$, $\Delta P_2=(u+\psi)^{p-1}-u^{p-1}$, $P_1'\psi=I_2*(pWu^{p-1}\psi)=pJ_{b,u}\psi$, $P_2'\psi=(p-1)u^{p-2}\psi$. Then $W\left[(P_1'\psi)P_2+P_1(P_2'\psi)\right]=\mathcal L_{b,u}\psi$ and
\[
\mathsf n_b(u+\psi)-\mathsf n_b(u)-\mathcal L_{b,u}\psi=W\left[(\Delta P_1-P_1'\psi)P_2+P_1(\Delta P_2-P_2'\psi)+\Delta P_1\Delta P_2\right].
\]
By (c), $|\Delta P_1-P_1'\psi|\le C\|\psi\|_X^2A$; by (K2) with $q=p-1$, $|\Delta P_2-P_2'\psi|\le C\|\psi\|_X^2U^{p-1}$; by (b), $|\Delta P_1|\le C\|\psi\|_XA$; by (K2), $|\Delta P_2|\le C\|\psi\|_XU^{p-1}$; and $P_1\le CA$, $P_2\le CU^{p-1}$. Hence the right side is bounded by $C\|\psi\|_X^2AU^{p-1}\le C\|\psi\|_X^2w_\beta$, and Lemma~\ref{lem:reg} gives the claim, provided $G*\mathcal L_{b,u}\psi$ makes sense: by (d), (e) and (K3), $|\mathcal L_{b,u}\psi|\le C\|\psi\|_X\left(\varrho\,w_\beta+U^{p-1}A\right)\le C\|\psi\|_Xw_\beta$, so $\|G*\mathcal L_{b,u}\psi\|_Y\le C\|\psi\|_X$.

\emph{(iii).} The bound on $\|T_b'(u)\|_{L(X,Y)}$ was just shown. For the difference, write
\begin{align*}
(\mathcal L_{b,u}-\mathcal L_{b',u'})\psi&=(p-1)(W-W')\,a_{b,u}\psi+(p-1)W'(a_{b,u}-a_{b',u'})\psi\\
&\quad+p\,(Wu^{p-1}-W'u'^{p-1})\,J_{b,u}\psi+p\,W'u'^{p-1}\,(J_{b,u}\psi-J_{b',u'}\psi).
\end{align*}
By (e) the first two terms are bounded by $C\delta\|\psi\|_X\varrho\,w_\beta$; by (d) and $|Wu^{p-1}-W'u'^{p-1}|\le C\delta U^{p-1}$ the third is bounded by $C\delta\|\psi\|_XU^{p-1}A$, and by (d) again so is the fourth. With (K3), $|(\mathcal L_{b,u}-\mathcal L_{b',u'})\psi|\le C\delta\|\psi\|_Xw_\beta$, and Lemma~\ref{lem:reg} gives the estimate. Since a Fr\'echet differentiable map with a locally Lipschitz derivative is $C^{1,1}$, $T_b\in C^{1,1}(\mathcal W_X;Y)$. Finally, with $u=u'=U$,
\[
\|L_b-L_{b'}\|_{L(Y)}\le\|T_b'(U)-T_{b'}'(U)\|_{L(X,Y)}+\|G*((b-b')\,\cdot)\|_{L(X,Y)}\le C\|b-b'\|_\infty ,
\]
by Lemma~\ref{lem:Q}(ii) for the second term and since $\|\cdot\|_X\le\|\cdot\|_Y$.

\emph{(iv).} The partial derivative of $(b,u)\mapsto T_b(u)$ in $u$ is $T_b'(u)$, continuous in $(b,u)$ with values in $L(X,Y)$ by (iii). In $b$, the map $b\mapsto T_b(u)$ is a quadratic polynomial in $b$ (it is $G*$ of $(1+b)\,I_2*((1+b)u^p)\,u^{p-1}$), so it is Fr\'echet differentiable on $L^\infty$ with
\[
\partial_bT_b(u)\hat b=G*\left[\hat b\,A_{b,u}u^{p-1}+W\,\left(I_2*(\hat bu^p)\right)u^{p-1}\right]\in Y,\qquad\|\partial_bT_b(u)\hat b\|_Y\le C\|\hat b\|_\infty ,
\]
by (a), (K1), (K3) and Lemma~\ref{lem:reg}. Moreover, by (b), (K2) and the same decomposition as in (i),
\begin{multline*}
\left|\hat b\,(A_{b,u}u^{p-1}-A_{b',u'}u'^{p-1})\right|+\left|\left(Wu^{p-1}-W'u'^{p-1}\right)I_2*(\hat bu^p)\right|\\
+\left|W'u'^{p-1}I_2*\left(\hat b(u^p-u'^p)\right)\right|\le C\delta\|\hat b\|_\infty w_\beta ,
\end{multline*}
so $\|\partial_bT_b(u)-\partial_bT_{b'}(u')\|_{L(L^\infty,Y)}\le C\delta$. Thus both partial derivatives exist and are continuous in $(b,u)$, and a map with continuous partial Fr\'echet derivatives is $C^1$. The statement for $F_b(w)=(U+w)-T_b(U+w)+G*(b(U+w))$ follows, because $(b,w)\mapsto G*(b(U+w))$ is a bounded bilinear map $L^\infty\times X\to Y$ by Lemma~\ref{lem:Q}(ii) (affine in $w$), and $w\mapsto U+w$, $Y\hookrightarrow X$ are bounded.

\emph{(v).} The set $\{w\in Y:\|w\|_X<\delta_0\}$ is open in $Y$ and convex, and $w\mapsto T_b(U+w)$ is $C^1$ on it with values in $Y$, by (iii) and the continuity of the inclusion $Y\hookrightarrow X$. Its derivative at $w$ is the restriction to $Y$ of $T_b'(U+w)$. Hence $Q_b'(w)=T_b'(U+w)-T_b'(U)$ and, by (iii) and $\|\psi\|_X\le\|\psi\|_Y$,
\[
\|Q_b'(w)\|_{L(Y)}\le\|T_b'(U+w)-T_b'(U)\|_{L(X,Y)}\le C\|w\|_X .
\]
Since $Q_b(0)=0$, integrating $Q_b(w)-Q_b(w')=\int_0^1Q_b'(w'+t(w-w'))(w-w')\dd t$ along the segment (which stays in the ball because $\|w'+t(w-w')\|_X\le\|w\|_X+\|w'\|_X$ and the ball is convex) gives the two estimates, with $\|w'+t(w-w')\|_X\le\|w\|_X+\|w'\|_X$ and $w'=0$ for the first.
\end{proof}

\subsection{The linearization and its principal part}\label{sec:KVlin}

For $b$ as above, set $W=1+b$ and $A_b=A_b(U)=A+I_2*(bU^p)$, $\varrho_b:=A_bU^{p-2}$. The linearization of $F_b$ at $0$ is
\[
L_b\varphi=\varphi+G*(b\varphi)-(p-1)G*\left[W\varrho_b\,\varphi\right]-p\,G*\left[W\,U^{p-1}\,I_2*\left(WU^{p-1}\varphi\right)\right],
\]
and $L_0=I-T_0'(U)=I-\mathbb T$. For $\varphi,\eta\in Y$ one has the symmetry
\[
\langle L_b\varphi,\eta\rangle_{H^1}=\int\nabla\varphi\cdot\nabla\eta+W\varphi\eta-(p-1)\int W\varrho_b\varphi\eta-p\int\left(I_2*(WU^{p-1}\varphi)\right)WU^{p-1}\eta ,
\]
which is symmetric in $(\varphi,\eta)$ (as in Lemma~\ref{lem:grad}).

We use the pointwise bounds, valid for $b$ with $\hat m\le1+b\le\hat m^{-1}$:
\begin{equation}\label{eq:rhob}
|\varrho_b-\varrho|\le\|\varrho\|_\infty\,\|b\|_\infty ,\qquad |\varrho_b-\varrho|\le C\left(\eta+\hat L\varepsilon\right)\quad(b\in\mathfrak B_{\hat m,\hat L}(\eta,\varepsilon)).
\end{equation}
The first holds because $|I_2*(bU^p)|\le\|b\|_\infty A$. For the second, $|I_2*(bU^p)|\le\eta A+C\hat L\varepsilon w_{N-2}$ as in Lemma~\ref{lem:KVres}, and $w_{N-2}U^{p-2}\le C w_{N-2}w_{(p-2)\beta}=C$, because $(p-2)\beta=-(N-2)$.

Put $\mu:=2-p>0$ and
\[
\mathcal D_W\varphi:=\varphi-G*\left[(1-\mu W)\varphi\right].
\]
A direct computation, using $b=W-1$ and $(p-1)W-b=1-\mu W$, gives
\begin{equation}\label{eq:Lsplit}
\begin{aligned}
L_b&=\mathcal D_W-K_{1,b}-K_{2,b},\\
K_{1,b}\varphi&=(p-1)G*\left[W(\varrho_b-1)\varphi\right],\\
K_{2,b}\varphi&=p\,G*\left[WU^{p-1}I_2*(WU^{p-1}\varphi)\right].
\end{aligned}
\end{equation}
Formally $\mathcal D_W=G*\,(-\Delta+\mu W)$: this is the far-field operator, with effective mass $\mu W\ge\mu\hat m>0$.

\begin{lemma}[principal part]\label{lem:KVprinc}
Let $\hat m\in(0,1]$. There is $M_1=M_1(N,p,\hat m)$ such that for every measurable $W$ with $\hat m\le W\le\hat m^{-1}$ the operator $\mathcal D_W$ is an isomorphism of $Y$ and $\|\mathcal D_W^{-1}\|_{L(Y)}\le M_1$.
\end{lemma}
\begin{proof}
\emph{$\mathcal D_W$ is bounded on $Y$.} For $\varphi\in Y$, $|(1-\mu W)\varphi|\le(1+\mu\hat m^{-1})\|\varphi\|_Xw_\beta$, so by Lemma~\ref{lem:reg} $\|G*[(1-\mu W)\varphi]\|_Y\le C\|\varphi\|_X$. (No regularity of $W$ is needed.)

\emph{The weighted resolvent.} Let $h\in X$. Since $X\subset L^2$ and $\mu W\ge\mu\hat m$, the Lax--Milgram theorem gives a unique $\psi=:R_Wh\in H^1$ with $-\Delta\psi+\mu W\psi=h$. We claim
\begin{equation}\label{eq:RW}
\|R_Wh\|_Y\le C\|h\|_X .
\end{equation}
Fix $R_0\ge1$ and put $\zeta(y)=(1+|y|^2/R_0^2)^{-\beta/2}$. Then $\zeta\in H^1$ (since $2\beta>N$), $\zeta$ is smooth, comparable to $w_\beta$ with constants depending on $R_0$, and a direct computation gives $|\Delta\zeta|\le C_1R_0^{-2}\zeta$ with $C_1=\beta(N+\beta+2)$. We choose $R_0^2\ge2C_1/(\mu\hat m)$, so that
\[
(-\Delta+\mu W)\zeta\ge\left(\mu\hat m-C_1R_0^{-2}\right)\zeta\ge\frac12\mu\hat m\,\zeta .
\]
With $\Lambda:=2\|h\|_X/(\mu\hat m\,c)$, where $w_\beta\le\zeta/c$, both $\Lambda\zeta-\psi$ and $\Lambda\zeta+\psi$ are $H^1$ functions with $(-\Delta+\mu W)(\Lambda\zeta\mp\psi)\ge\Lambda\frac12\mu\hat m\zeta-|h|\ge0$. Testing with the negative part and using $\mu W>0$ gives $\Lambda\zeta\mp\psi\ge0$, i.e.\ $|\psi|\le\Lambda\zeta\le C\|h\|_Xw_\beta$. Finally $\int|\nabla\psi|^2\le\int|h||\psi|\le\|h\|_2\|\psi\|_2\le C\|h\|_X^2$. This proves \eqref{eq:RW}.

\emph{Injectivity.} If $\mathcal D_W\varphi=0$ with $\varphi\in Y$, then $\varphi=G*[(1-\mu W)\varphi]\in H^2$ and $(-\Delta+1)\varphi=(1-\mu W)\varphi$, i.e.\ $-\Delta\varphi+\mu W\varphi=0$. Testing with $\varphi$ gives $\varphi=0$.

\emph{Surjectivity and the bound.} Given $f\in Y$, put $h=(1-\mu W)f\in X$, $\psi=R_Wh\in Y$ and $\varphi=f+\psi\in Y$. In $H^{-1}$,
\[
(-\Delta+1)(\varphi-f)=(-\Delta+\mu W)\psi+(1-\mu W)\psi=h+(1-\mu W)\psi=(1-\mu W)\varphi,
\]
so $\varphi-f=G*[(1-\mu W)\varphi]$, i.e.\ $\mathcal D_W\varphi=f$. Moreover $\|\varphi\|_Y\le\|f\|_Y+C\|h\|_X\le M_1\|f\|_Y$.
\end{proof}

\subsection{Uniform invertibility}

Let $\mathcal Y:=Y_\perp\times\R^N$ with the norm $\|(\varphi,\theta)\|=\|\varphi\|_Y+|\theta|$, and define
\[
\mathfrak L_b:\mathcal Y\to Y,\qquad \mathfrak L_b(\varphi,\theta)=L_b\varphi-\sum_i\theta_i\partial_iU .
\]
This is the linearization at $(0,0)$ of the Lyapunov--Schmidt problem $F_b(w)=\sum_ic_i\partial_iU$, $w\in Y_\perp$ (cf.\ Proposition~\ref{prop:natural}); the Lagrange multipliers $c$ play the role of $\Pi^\perp$ of Section~\ref{sec:fixedpoint}.

\begin{lemma}[uniform invertibility]\label{lem:KVU}
Let $\hat m\in(0,1]$ and $\hat L>0$. There are $\eta_1\in(0,\frac12]$, $\varepsilon_1>0$ and $C_U>0$, depending only on $N,p,\hat m,\hat L$, such that for all $\eta\le\eta_1$, $\varepsilon\le\varepsilon_1$ and $b\in\mathfrak B_{\hat m,\hat L}(\eta,\varepsilon)$ the map $\mathfrak L_b$ is an isomorphism and
\[
\|\varphi\|_Y+|\theta|\le C_U\,\left\|L_b\varphi-\textstyle\sum_i\theta_i\partial_iU\right\|_Y\qquad\text{for all }(\varphi,\theta)\in\mathcal Y .
\]
\end{lemma}

\begin{proof}
For $\varphi\in Y$ and $R\ge1$ put
\[
\mathcal N_R(\varphi):=\sup_{|y|\le R}\frac{|\varphi(y)|}{w_\beta(y)},\qquad \omega(R):=\sup_{|y|>R}|\varrho(y)-1| .
\]

\emph{Step 1: operator bounds.} First, by Lemma~\ref{lem:KVQ}(iii), $\|L_b-L_{b'}\|_{L(Y)}\le C_3\|b-b'\|_\infty$ for $b,b'\in\mathcal O(\hat m)$.

Second, for $b\in\mathfrak B_{\hat m,\hat L}(\eta,\varepsilon)$, $\varphi\in Y$ and every $R\ge1$,
\begin{equation}\label{eq:Kbound}
\|K_{1,b}\varphi+K_{2,b}\varphi\|_Y\le C_2\left(\mathcal N_R(\varphi)+\left[\omega(R)+R^{N-p\beta}+\eta+\varepsilon\right]\|\varphi\|_Y\right),
\end{equation}
with $C_2=C_2(N,p,\hat m,\hat L)$ independent of $R$. Proof of \eqref{eq:Kbound}. For $K_{1,b}$: by \eqref{eq:rhob}, $|W(\varrho_b-1)\varphi|\le\hat m^{-1}\left(|\varrho-1|+C(\eta+\hat L\varepsilon)\right)|\varphi|$. On $\{|y|\le R\}$ we bound $|\varrho-1||\varphi|\le\|\varrho-1\|_\infty\mathcal N_Rw_\beta$, and on $\{|y|>R\}$, $|\varrho-1||\varphi|\le\omega(R)\|\varphi\|_Xw_\beta$; Lemma~\ref{lem:reg} gives the claim. For $K_{2,b}$: set $g=WU^{p-1}\varphi$. Then $|g|\mathbf 1_{B_R}\le C\mathcal N_R(\varphi)w_{p\beta}$ and $|g|\mathbf 1_{B_R^c}\le C\|\varphi\|_Xw_{p\beta}\mathbf 1_{B_R^c}$, so by Lemma~\ref{lem:KVtail} (with $s=p\beta>N$)
\[
|I_2*g|\le C\left(\mathcal N_R(\varphi)+R^{N-p\beta}\|\varphi\|_X\right)\,w_{N-2}.
\]
Multiplying by $WU^{p-1}\le C\hat m^{-1}w_{(p-1)\beta}$ and using $(p-1)\beta+N-2=\beta$ gives a bound by $C(\dots)w_\beta$, and Lemma~\ref{lem:reg} again yields the claim.

\emph{Step 2: far-field estimate.} Let $M_1$ be as in Lemma~\ref{lem:KVprinc}. By \eqref{eq:Lsplit},
\[
\varphi=\mathcal D_W^{-1}\left(L_b\varphi+K_{1,b}\varphi+K_{2,b}\varphi\right),
\]
so by \eqref{eq:Kbound}
\[
\|\varphi\|_Y\le M_1\|L_b\varphi\|_Y+M_1C_2\left(\mathcal N_R(\varphi)+\left[\omega(R)+R^{N-p\beta}+\eta+\varepsilon\right]\|\varphi\|_Y\right).
\]
Fix $R_1\ge1$ with $M_1C_2\left(\omega(R_1)+R_1^{N-p\beta}\right)\le\frac18$ and then $\eta_1,\varepsilon_1$ with $M_1C_2(\eta_1+\varepsilon_1)\le\frac18$ (further conditions on $\eta_1$ are imposed below). For $\eta\le\eta_1$ and $\varepsilon\le\varepsilon_1$,
\begin{equation}\label{eq:farfield}
\|\varphi\|_Y\le2M_1\|L_b\varphi\|_Y+2M_1C_2\,\mathcal N_{R_1}(\varphi)\qquad(\varphi\in Y).
\end{equation}

\emph{Step 3: the a priori estimate.} Let $\mathfrak L_0$ be the operator for $b=0$; by Lemma~\ref{lem:Lperp}, $L_0:Y_\perp\to Y_\perp$ is an isomorphism and $Y=Y_\perp\oplus\operatorname{span}\{\partial_iU\}$, so $\mathfrak L_0$ is an isomorphism. Let $C_0:=\|\mathfrak L_0^{-1}\|$ and let $C_3$ be the Lipschitz constant of $b\mapsto L_b$, in the sense of Step~1 for $\|b\|_\infty\le\frac12$. We impose in addition $C_0C_3\eta_1\le\frac12$.

We claim that there are $\varepsilon_1,C_U$ such that the stated inequality holds for all $\eta\le\eta_1$ and $\varepsilon\le\varepsilon_1$. Suppose not. Then for every $n$ (taking $\varepsilon_1=1/n$, $C_U=n$) there are $\varepsilon_n\le1/n$, $\eta_n\in[0,\eta_1]$, $b_n\in\mathfrak B(\eta_n,\varepsilon_n)$ and $(\varphi_n,\theta_n)\in\mathcal Y$ with $\|\varphi_n\|_Y+|\theta_n|=1$ and
\[
f_n:=L_{b_n}\varphi_n-\sum_i\theta_{n,i}\partial_iU\longrightarrow0\quad\text{in }Y .
\]
Passing to a subsequence, $\theta_n\to\theta_\infty$, $\eta_n\to\eta_\infty$ and $b_n(0)\to\beta_\infty$, where $|\beta_\infty|\le\eta_\infty\le\eta_1$ and $\hat m\le1+\beta_\infty\le\hat m^{-1}$. Since $|b_n(y)-b_n(0)|\le\hat L\varepsilon_n|y|\to0$, we get $b_n\to\beta_\infty$ locally uniformly, and $\beta_\infty$ is a \emph{constant}. (This is the only place where the Lipschitz condition in the definition of $\mathfrak B$ is used.)

\emph{Local convergence.} From $L_{b_n}\varphi_n=f_n+\sum_i\theta_{n,i}\partial_iU$ we get
\begin{multline*}
\varphi_n=f_n+\sum_i\theta_{n,i}\partial_iU-G*(b_n\varphi_n)+(p-1)G*[W_n\varrho_{b_n}\varphi_n] + {}\\
{}+p\,G*\left[W_nU^{p-1}I_2*(W_nU^{p-1}\varphi_n)\right].
\end{multline*}
The three convolution terms are bounded in $C^{1,\gamma'}(\R^N)$ by Lemma~\ref{lem:reg}, $\partial_iU\in C^{0,\gamma'}$, and $f_n\to0$ uniformly. By Arzel\`a--Ascoli, after passing to a subsequence, $\varphi_n\to\varphi_\infty$ locally uniformly, and, being bounded in $H^1$, also weakly in $H^1$. Thus $\varphi_\infty\in Y$ and $\langle\varphi_\infty,\partial_jU\rangle_{H^1}=0$. Using dominated convergence (all integrands are dominated by $C\,G(y-y')w_\beta(y')$ or by $C\,I_2(y'-y'')w_{p\beta}(y'')$) and $b_n\to\beta_\infty$, $W_n\to1+\beta_\infty$, $\varrho_{b_n}\to(1+\beta_\infty)\varrho$ pointwise, we pass to the limit and get $\mathfrak L_{\beta_\infty}(\varphi_\infty,\theta_\infty)=0$, i.e.\ $L_{\beta_\infty}\varphi_\infty=\sum_i\theta_{\infty,i}\partial_iU$, where $L_{\beta_\infty}$ is $L_b$ with the constant $b\equiv\beta_\infty$ (for which $W=1+\beta_\infty$, $A_b=(1+\beta_\infty)A$ and $\varrho_b=(1+\beta_\infty)\varrho$). Writing $\mathfrak L_{\beta_\infty}=\mathfrak L_0+(L_{\beta_\infty}-L_0)$ and using Step~1,
\[
\|(\varphi_\infty,\theta_\infty)\|\le C_0\,\|(L_{\beta_\infty}-L_0)\varphi_\infty\|_Y\le C_0C_3\eta_1\|\varphi_\infty\|_Y\le\frac12\|(\varphi_\infty,\theta_\infty)\| ,
\]
so $\varphi_\infty=0$ and $\theta_\infty=0$. In particular $\mathcal N_{R_1}(\varphi_n)\to0$.

\emph{Conclusion of the contradiction.} Since $\theta_n\to0$ and $f_n\to0$, $\|L_{b_n}\varphi_n\|_Y\le\|f_n\|_Y+|\theta_n|\max_i\|\partial_iU\|_Y\to0$. By \eqref{eq:farfield}, $\|\varphi_n\|_Y\to0$. Together with $\theta_n\to0$ this contradicts $\|\varphi_n\|_Y+|\theta_n|=1$. This proves the inequality.

\emph{Step 4: surjectivity.} Let $b\in\mathfrak B(\eta,\varepsilon)$ and $b_t=tb$, $t\in[0,1]$. Then $b_t\in\mathfrak B(\eta,\varepsilon)$ (as $1+tb$ is a convex combination of $1$ and $1+b$). The a priori estimate holds for every $b_t$ with the same constant $C_U$, $t\mapsto\mathfrak L_{b_t}$ is continuous in operator norm by Step~1, and $\mathfrak L_{b_0}=\mathfrak L_0$ is an isomorphism. By the method of continuity, $\mathfrak L_{b_t}$ is invertible for all $t$, in particular for $t=1$.
\end{proof}

\subsection{Reduction at a fixed center}\label{sec:KVred}

\begin{proposition}[reduction, fixed $\xi$]\label{prop:KVred}
Let $\hat m\in(0,1]$ and $\hat L>0$. There are $\eta_2,\varepsilon_2,r_0>0$ and $C_4$, depending only on $N,p,\hat m,\hat L$, such that for all $\eta\le\eta_2$, $\varepsilon\le\varepsilon_2$ and $b\in\mathfrak B_{\hat m,\hat L}(\eta,\varepsilon)$ the problem
\[
F_b(w)=\sum_ic_i\,\partial_iU,\qquad (w,c)\in Y_\perp\times\R^N ,
\]
has a unique solution in $\{\|w\|_Y+|c|\le r_0\}$. It satisfies
\[
\|\tilde w\|_Y+|\tilde c|\le C_4\left(\eta+\chi(\varepsilon)+\varepsilon\right),\qquad U+\tilde w\ge\frac12U .
\]
If $\tilde c=0$, then $U+\tilde w$ is a positive solution of \eqref{eq:KVresc} with $W=1+b$.
\end{proposition}
\begin{proof}
Let $\eta_1,\varepsilon_1,C_U$ be as in Lemma~\ref{lem:KVU} and $C_Q'$, $\delta_0$ as in Lemma~\ref{lem:KVQ}. Since $F_b(w)=F_b(0)+L_bw-Q_b(w)$, the problem is equivalent to
\[
(w,c)=\Gamma(w,c):=\mathfrak L_b^{-1}\left[-F_b(0)+Q_b(w)\right]\quad\text{in }\mathcal Y,
\]
where $\|\mathfrak L_b^{-1}\|\le C_U$. Take $r_0:=\min\{\delta_0,\,c_U/2,\,(4C_UC_Q')^{-1}\}$ and let $B=\{\|w\|_Y+|c|\le r_0\}$. For $(w,c),(w',c')\in B$ (so $\|w\|_X\le\|w\|_Y\le r_0<\delta_0$),
\[
\|\Gamma(w,c)\|\le C_U\|F_b(0)\|_Y+C_UC_Q'r_0^2,
\]
\[
\|\Gamma(w,c)-\Gamma(w',c')\|\le2C_UC_Q'r_0\,\|w-w'\|_Y\le\frac12\|(w,c)-(w',c')\|.
\]
By Lemma~\ref{lem:KVres}, $\|F_b(0)\|_Y\le C_1(\eta+\chi(\varepsilon)+\varepsilon)$, so if $\eta_2,\varepsilon_2$ satisfy $C_UC_1(\eta_2+\chi(\varepsilon_2)+\varepsilon_2)\le r_0/2$ (possible since $\chi(\varepsilon)\to0$; $\chi$ is nondecreasing in $\varepsilon$, so the bound holds for all $\varepsilon\le\varepsilon_2$), then $\Gamma(B)\subset B$ (we use $C_UC_Q'r_0^2\le r_0/4$). Banach's theorem gives a unique fixed point in $B$. The fixed point satisfies $\|x\|\le C_U\|F_b(0)\|+\frac14\|x\|$, whence $\|x\|\le\frac43C_U\|F_b(0)\|$. Finally, $|\tilde w|\le\|\tilde w\|_Xw_\beta\le(\|\tilde w\|_Y/c_U)U\le\frac12U$ because $r_0\le c_U/2$.
\end{proof}

\begin{remark}\label{rem:KVcenter}
With $\eta=0$ (case (a)) the proposition gives $\|\tilde w\|_Y+|\tilde c|\to0$ as $\varepsilon\to0$, uniformly in the potential class, i.e.\ for every $\xi$ and every $V$ satisfying (W). In particular the oscillation of $V$ enters only through $\hat m=m/M$ and $\hat L$.
\end{remark}

\subsection{The reduced equation}\label{sec:KVreduced}

We now sketch how Subsection~\ref{sec:KVred} is used to derive the analogues of Section~\ref{sec:semiclassical}.

\paragraph{Frame.} Normalize $V(\xi_0)=1$ by Lemma~\ref{lem:dil} (and replace $\xi_0$ by $\sqrt{V(\xi_0)}\xi_0$). Choose $r>0$ so small that $\Omega_0:=B_r(\xi_0)$ satisfies $|V-1|\le\eta_*:=\min\{\eta_1,\eta_2\}$ on $\Omega_0$. For $\xi\in\Omega_0$ and $\varepsilon\le\varepsilon_*:=\min\{\varepsilon_1,\varepsilon_2\}$ put $b_\xi=V(\xi+\varepsilon\,\cdot)-1\in\mathfrak B_{\hat m,L}(\eta_*,\varepsilon)$, and let $(\tilde w_\xi,\tilde c_\xi)$ be given by Proposition~\ref{prop:KVred}. Using the \emph{fixed} frame $Y_\perp$ (orthogonal to $\partial_iU$, $\lambda_0=1$) avoids differentiating dilations with respect to $\xi$: we use only $\xi\mapsto b_\xi\in L^\infty$, which is $C^1$ with $\partial_{\xi_j}b_\xi=(\partial_jV)(\xi+\varepsilon\,\cdot)$ by Lemma~\ref{lem:bxi}. Here the global oscillation of $V$ does not matter: only $|V-1|\le\eta_*$ near $\xi_0$ is used, and $\eta_*$ depends on $N,p,m,M,L$.

\paragraph{$C^1$ dependence.} The map $(b,w,c)\mapsto F_b(w)-\sum c_i\partial_iU$ is $C^1$ on its domain, by Lemma~\ref{lem:KVQ}(iv). At the solution the linearization is $\mathfrak L_b-(Q_b'(\tilde w),0)$, which is invertible because $\|Q_b'(\tilde w)\|\le C_Q'\|\tilde w\|_X\le C_Q'r_0\le\frac1{4C_U}$. The implicit function theorem gives $\xi\mapsto(\tilde w_\xi,\tilde c_\xi)$ of class $C^1$.

\paragraph{Natural constraint.} Put $\mathcal J(b,v)=\frac12\int|\nabla v|^2+(1+b)v^2-\frac1{2p}\int\left(I_2*((1+b)v^p)\right)(1+b)v^p$ and $\Phi_\varepsilon(\xi)=\mathcal J(b_\xi,U+\tilde w_\xi)$. One has $\partial_v\mathcal J(b,v)\eta=\langle F_b(v-U),\eta\rangle_{H^1}$ and
\[
\partial_b\mathcal J(b,v)[\hat b]=\frac12\int\hat bv^2-\frac1p\int\left(I_2*((1+b)v^p)\right)\hat b\,v^p .
\]
By the chain rule and $\partial_{\xi_j}\tilde w\in Y_\perp\perp g:=F_b(\tilde w)=\sum\tilde c_i\partial_iU$,
\begin{equation}\label{eq:KVgrad}
\partial_j\Phi_\varepsilon(\xi)=\int(\partial_jV)(\xi+\varepsilon z)\left[\frac12v^2-\frac1p\left(I_2*((1+b_\xi)v^p)\right)v^p\right]\dd z,\qquad v=U+\tilde w_\xi .
\end{equation}
Translation covariance, $\mathcal J(b,v(\cdot+\tau e_j))=\mathcal J(b(\cdot-\tau),v)$, gives $\langle g,\partial_jv\rangle_{H^1}=-\varepsilon\,\partial_j\Phi_\varepsilon$. Expanding the left side as in Proposition~\ref{prop:natSC} yields
\[
(\kappa I-B)\tilde c=-\varepsilon\nabla\Phi_\varepsilon(\xi),\qquad B_{ij}=\langle\partial_{ij}U,\tilde w_\xi\rangle_{H^1},\quad\|B\|\le Nm_3C_*\|\tilde w\|_Y<\kappa
\]
shrinking $r_0$ and $\eta_*$ if necessary. So $\tilde c_\xi=0$ if and only if $\nabla\Phi_\varepsilon(\xi)=0$, and $|\tilde c_\xi|\le C\varepsilon$.

\paragraph{The constant-potential problem.} For $b\equiv\lambda-1$ with $|\lambda-1|\le\eta_*$ one has $\tilde c=0$ and $U+\tilde w=U_\lambda=U(\sqrt\lambda\,\cdot)$. The proof is that of Proposition~\ref{prop:redb}(iii): $\langle g,\partial_jv\rangle_{H^1}=0$ by translation invariance (constant $b$), so $\tilde c=0$; hence $v>0$ solves $-\Delta v+\lambda v=\lambda^2(I_2*v^p)v^{p-1}$, and by Lemma~\ref{lem:dil}, the uniqueness theorem of \cite{XZZ2} and the $O(N)$-equivariance of the fixed point, $v=U_\lambda$.

\paragraph{The limit $\varepsilon\to0$.} Let $\xi_n\to\bar\xi$ and $\varepsilon_n\to0$. As in Lemma~\ref{lem:limitSC} (compactness in $C^{0,\gamma'}_{\mathrm loc}$, weak limit in $Y$, dominated convergence in $T_b$, $b_{\xi_n}\to V(\bar\xi)-1$ locally uniformly and boundedly), every limit of $U+\tilde w_{\xi_n}$ solves the constant-potential problem and lies in the ball of Proposition~\ref{prop:KVred}; by uniqueness it equals $U_{\bar\lambda}$, $\bar\lambda=V(\bar\xi)$. Inserting this into \eqref{eq:KVgrad},
\[
\nabla\Phi_{\varepsilon_n}(\xi_n)\longrightarrow h(\bar\xi):=\nabla V(\bar\xi)\left[\frac12\|U_{\bar\lambda}\|_{L^2}^2-\frac1p\bar\lambda\!\int\left(I_2*U_{\bar\lambda}^p\right)U_{\bar\lambda}^p\right]=-c_*\,V(\bar\xi)^{-N/2}\,\nabla V(\bar\xi),
\]
where $c_*=\frac{(p-1)(N-2)}{4p}\|U\|_{H^1}^2>0$. The precise value of the constant $c_*$ comes from the fact that $h=\nabla(\Upsilon\circ V)$ with the energy $\Upsilon(s)=\frac{p-1}{2p}\|U\|_{H^1}^2\,s^{1-N/2}$ of the frozen profile $U_s$.
%
%
%
The convergence is uniform on compact sets of $\Omega_0$ by the usual contradiction argument.

\begin{theorem}\label{thm:KV}
Let $V$ satisfy (W) and put $\Upsilon(s)=\frac{p-1}{2p}\|U\|_{H^1}^2s^{1-N/2}$, a strictly decreasing function.
\begin{enumerate}[label=(\roman*),nosep]
\item Let $\xi_0$ be an isolated critical point of $V$ with $\deg(\nabla V,B_r(\xi_0),0)\ne0$ for small $r$ (for instance a nondegenerate critical point or a strict local extremum). Then for every small $r>0$ there is $\varepsilon_r>0$ such that for $\varepsilon<\varepsilon_r$ equation \eqref{eq:KV} has a positive solution $u_\varepsilon$ of the form $u_\varepsilon(x)=v_\varepsilon\left(\sqrt{V(\xi_0)}\,(x-\xi_\varepsilon)/\varepsilon\right)$ with $|\xi_\varepsilon-\xi_0|<r$, $v_\varepsilon\ge\frac12U$ and $\|v_\varepsilon-U(\sqrt{V(\xi_\varepsilon)/V(\xi_0)}\,\cdot)\|_{L^2}\to0$.
\item Let $\Omega$ be a bounded open set with $\max_{\overline\Omega}V>\max_{\partial\Omega}V$, and $K=\{x\in\overline\Omega:V(x)=\max_{\overline\Omega}V\}$. For every $\delta>0$ with $K_\delta\subset\Omega$ and every sufficiently small $\varepsilon$ (depending on $\delta$), \eqref{eq:KV} has at least $\operatorname{cat}(K,K_\delta)$ pairwise distinct positive solutions concentrating at $K$ as $\varepsilon\to0$.
\end{enumerate}
\end{theorem}
\begin{proof}[Proof (sketch)]
(i) Normalize $V(\xi_0)=1$ and work on $\Omega_0=B_r(\xi_0)$ as above. The map $\nabla\Phi_\varepsilon$ is continuous and converges uniformly on $\overline{B_r(\xi_0)}$ to $h=-c_*V^{-N/2}\nabla V$. Since $\omega=-c_*V^{-N/2}$ is continuous and negative, $\deg(h,B_r(\xi_0),0)=\deg(\omega\nabla V)=(-1)^N\deg(\nabla V,B_r(\xi_0),0)\neq0$ by the homotopy $((1-t)\omega-t)\nabla V$. Brouwer degree gives $\xi_\varepsilon$ with $\nabla\Phi_\varepsilon(\xi_\varepsilon)=0$, hence $\tilde c=0$ and the solution.

(ii) After a dilation assume $\max_{\overline\Omega}V=1$. Then $V\le1$ on $\overline\Omega$ and $K=\{V=1\}\cap\overline\Omega\subset\Omega$ is compact. Choose $\delta'\in(0,\delta]$ with $V>1-\eta_*$ on $K_{\delta'}$ and put $\Omega'=\operatorname{int}K_{\delta'}$. Then $\partial\Omega'\subset\overline\Omega\setminus K$, so $\max_{\partial\Omega'}V<1$ by compactness, and $K_{\delta''}\subset\Omega'$ for $\delta''<\delta'$. The reduction of this subsection applies on $\Omega'$, and $\Phi_\varepsilon$ is $C^1$-close to $\Upsilon\circ V$ up to an additive constant (the analogue of Corollary~\ref{cor:C1conv}). Since $\Upsilon$ is decreasing, $\Upsilon\circ V$ attains its minimum over $\overline{\Omega'}$ exactly on $K$, and Lemma~\ref{lem:LS} gives $\operatorname{cat}(K,K_{\delta''})$ critical points of $\Phi_\varepsilon$ in $\Omega'$. Because $K_{\delta''}\subset K_\delta$, $\operatorname{cat}(K,K_{\delta''})\ge\operatorname{cat}(K,K_\delta)$. Distinct critical points give distinct solutions by the argument of Lemma~\ref{lem:orbitsep}, which only uses $\|\tilde w_\xi\|_Y\le r_0$ small, $\tilde w_\xi\in Y_\perp$ and $v_\xi\in H^2$.
\end{proof}

Some concluding remarks are in order.

\begin{enumerate}[label=(\arabic*)]
\item The proof of Theorem~\ref{thm:KV} uses the global oscillation of $V$ only through the constants $m$, $M$, $L$ (hence $\hat m$, $\hat L$), never through a smallness condition. The condition $|V-1|\le\eta_*$ is imposed only on a neighbourhood of the point or set where the solutions concentrate, and such a neighbourhood always exists by continuity.
\item The proof of Lemma~\ref{lem:KVU} is where the structure $K=V$ enters, through the cancellation behind Lemma~\ref{lem:KVres} and the coercivity of the principal part $-\Delta+(2-p)W$ in Lemma~\ref{lem:KVprinc}. For a potential $K\ne V$ on the right-hand side, the residual contains a term of the order of $\sup|K/V-1|$ that does not decay at infinity, and one needs this quantity to be small.
\item The constants $\eta_*,\varepsilon_*$ are not explicit. Also, case (ii) of the theorem is stated with the weak requirement $\max_{\overline\Omega}V>\max_{\partial\Omega}V$ of the paper, but it uses only the behaviour of $V$ near $K$.
\end{enumerate}

\end{document}